\documentclass[a4paper,11pt]{amsart}
\usepackage{amsmath,amsthm,amssymb}
\usepackage[mathscr]{eucal}
 \usepackage{cite}
\usepackage{upgreek}
\usepackage[bookmarks,bookmarksnumbered]{hyperref}
\usepackage{enumerate}
\usepackage{amsmath}
\usepackage{mathrsfs}
\usepackage{blindtext}
\usepackage{scrextend}
\usepackage{enumitem}

\usepackage{amsmath}
\newcommand{\Mod}[1]{\ (\text{mod}\ #1)}
\addtokomafont{labelinglabel}{\sffamily}

\numberwithin{equation}{section}

\theoremstyle{plain}
\newtheorem{main}{Theorem}
\newtheorem{mcor}[main]{Corollary}

\newtheorem{theorem}{Theorem}[section]

\newtheorem{lemma}[theorem]{Lemma}
\newtheorem{proposition}[theorem]{Proposition}

\theoremstyle{definition}
\newtheorem{definition}[theorem]{Definition}

\newtheorem{notation}[theorem]{Notation}
\newtheorem{remark}[theorem]{Remark}

\newcommand{\R}{\mathcal{R}}
\renewcommand{\S}{\mathcal{S}}
\newcommand{\T}{\mathcal{T}}

\begin{document}

\title[Prime II$_1$ factors from lattices in products of rank one groups]
{Prime II$_1$ factors arising from irreducible lattices in products of rank one simple Lie groups}

\author[D. Drimbe]{Daniel Drimbe}
\address{Department of Mathematics, University of California San Diego, 9500 Gilman Drive, La Jolla, CA 92093, USA.}
\email{ddrimbe@ucsd.edu}
\thanks{D.D. was partially supported by NSF Career Grant DMS \#1253402.}

\author[D. Hoff]{Daniel Hoff}
\address{Department of Mathematics, University of California Los Angeles, Box 951555, 520 Portola Plaza, Los Angeles, CA 90095-1555, USA.}
\email{hoff@math.ucla.edu}
\thanks{D.H. was partially supported by NSF GRFP Grant DGE \#1144086, and NSF RTG DMS \#1344970.}

\author[A. Ioana]{Adrian Ioana}
\address{Department of Mathematics, University of California San Diego, 9500 Gilman Drive, La Jolla, CA 92093, USA.}
\email{aioana@ucsd.edu}
\thanks{A.I. was supported in part by  NSF Career Grant DMS \#1253402, and a Alfred. P. Sloan Foundation Fellowship.}


\begin{abstract} 
We prove that if $\Gamma$  is an icc irreducible lattice in a product of connected non-compact rank one simple Lie groups with finite center, then the II$_1$ factor $L(\Gamma)$ is prime. In particular, we deduce that the II$_1$ factors associated to the arithmetic groups $\text{PSL}_2(\mathbb Z[\sqrt{d}])$ and $\text{PSL}_2(\mathbb Z[S^{-1}])$ are prime, for any square-free integer $d\geq 2$ with $d\not\equiv 1\Mod{4}$  and any finite non-empty set of primes $S$.
This provides the first examples of prime II$_1$ factors arising from lattices in higher rank semisimple Lie groups. 
More generally, we describe all tensor product decompositions of $L(\Gamma)$ for icc countable groups $\Gamma$ that are measure equivalent to a product of non-elementary hyperbolic groups. In particular, we show that $L(\Gamma)$ is prime, unless $\Gamma$ is a product of infinite groups, in which case we prove a unique prime factorization result for $L(\Gamma)$.

\end{abstract}

\maketitle

\section{Introduction}

\subsection{Background and statement of results} An important theme in operator algebras is the study of tensor product decompositions of II$_1$ factors. A II$_1$ factor $M$ is called {\it prime} if it is not isomorphic to a tensor product of II$_1$ factors. 
In \cite{Po83}, Popa proved that the free groups on uncountably many generators give rise to prime II$_1$ factors.  By using Voiculescu's free probability theory, Ge showed that the free group factors $L(\mathbb F_n)$, $2\leq n\leq\infty$, are also prime, thus providing the first examples of separable prime II$_1$ factors \cite{Ge96}. Ozawa then used subtle C$^*$-algebraic methods to prove that for any icc hyperbolic group $\Gamma$, the II$_1$ factor $L(\Gamma)$ is {\it solid}, that is, the  relative commutant of any diffuse subalgebra is amenable \cite{Oz03}. Since solid non-amenable II$_1$ factors are clearly prime, this recovers the primeness of $L(\mathbb F_n)$. By developing a new technique based on closable derivations, Peterson proved that $L(\Gamma)$ is prime, for any non-amenable icc group $\Gamma$ that admits an unbounded $1$-cocycle into its left regular representation \cite{Pe06}.   Popa then used his powerful deformation/rigidity theory to  give a new proof of solidity of $L(\mathbb F_n)$ \cite{Po06b}.  For additional primeness results, see \cite{Oz04, Po06a, CI08, CH08, Va10b, Bo12, HV12, DI12, CKP14, Ho15}.

A common feature 
of these results is that the groups $\Gamma$ for which $L(\Gamma)$ was proven to be prime have   
 ``rank one" properties, such as hyperbolicity or the existence of certain unbounded (quasi) 1-cocycles.   On the other hand, the primeness question for the ``higher rank" arithmetic groups $\text{PSL}_n(\mathbb Z)$, $n\geq 3$, is notoriously hard and remains open. Moreover, in spite of the remarkable advances made in the study of II$_1$ factors in the last 15 years (see the surveys \cite{Po07,Va10a,Io12b}), little is known about the structure of II$_1$ factors associated to lattices in higher rank semisimple Lie groups.
 In fact, while II$_1$ factors arising from lattices in connected rank one simple Lie groups have already been shown to be prime in \cite{Oz03}, not a single example of a lattice, whose II$_1$ factor is prime, in either a higher rank simple or  semisimple Lie group is known.

Our first main result provides the first examples of lattices in higher rank semisimple Lie groups which give rise to prime II$_1$ factors. More precisely, we prove:

\begin{main}\label{A} 
If $\Gamma$ is an icc irreducible lattice in a product $G=G_1\times...\times G_n$ of $n\geq 1$ connected non-compact rank one simple real Lie groups with finite center, then the II$_1$ factor $L(\Gamma)$ is prime.

More generally, if $\Gamma\in\mathscr L$ is an icc group, then the II$_1$ factor $L(\Gamma)$ is prime.
\end{main}

Before stating a consequence of Theorem \ref{A}, we will first explain the terminology used, give several examples of groups to which Theorem \ref{A} applies, and compare it with a result in the literature.

\begin{definition}\label{L} We denote by $\mathscr L$ the family of countable  groups $\Gamma$ which can be realized as an irreducible lattice in a product $G=G_1\times...\times G_n$ of $n\geq 1$ locally compact second countable groups such that  (i) $G_j$ admits a lattice that is a non-elementary hyperbolic group, for every $1\leq j\leq n$,  (ii) $G_j$ does not admit an open normal compact subgroup, for some $1\leq j\leq n$, and (iii) $\Gamma$ does not contain a non-trivial element which commutes with an open subgroup of $G$.

A subgroup $\Gamma<G$ is called a {\bf lattice} if it is discrete and the homogeneous space $G/\Gamma$ carries a $G$-invariant Borel probability measure.
A lattice $\Gamma<G$ in a product group $G=G_1\times...\times G_n$ is called {\bf irreducible} if its projection onto $\underset{i\not=j}{\times}G_i$ is dense, for every $1\leq j\leq n$.
\end{definition}

\begin{remark}\label{L} Assume that $G_j$ is a connected non-compact simple Lie group (or algebraic group), for every $1\leq j\leq n$. Then  condition (ii) of Definition \ref{L} is satisfied. 
Since any element of $G_j,1\leq j\leq n$, which commutes with an open subgroup is necessarily central, condition (iii) is satisfied by any icc lattice $\Gamma<G$.  Here, we point out that if $G_1,...,G_n$ are of rank one then condition (i) is also satisfied, and provide several examples of countable groups belonging to $\mathscr L$.

\begin{enumerate}

\item 
If $G$ is a connected non-compact rank one simple real Lie group with finite center, then any co-compact lattice $\Lambda<G$ is non-elementary hyperbolic.  
This in particular applies to $G=\text{SL}_2(\mathbb R)$. Moreover, in this case,   $\text{SL}_2(\mathbb Z)$ and the free group $\mathbb F_2$ arise as lattices of $G$.

\item Let $G$ be a rank one simple algebraic group over a locally compact non-archimedean field. 
Then any such group admits a lattice $\Lambda<G$ which is a finitely generated free group (see \cite[Corollaries 4.8 and 4.14]{BK90} and \cite[Theorem 2.1]{Lu91}). In particular, this applies to $G=\text{SL}_2(\mathbb Q_p)$, where $\mathbb Q_p$ denotes the field of $p$-adic numbers for a prime $p$.

\item Let $d\geq 2$ be a square-free integer and $S$ be a finite non-empty set of primes. Denote by $\mathcal O_d$ the ring $\mathbb Z[\frac{1+\sqrt{d}}{2}]$, if $d\equiv 1\Mod 4$, and the ring $\mathbb Z[\sqrt{d}]$, otherwise. Denote by $\mathbb Z[S^{-1}]$ the ring of rational numbers whose denominators have all prime factors from $S$. Then $\text{SL}_2(\mathcal O_d)$ and $\text{SL}_2(\mathbb Z[S^{-1}])$ can be realized as irreducible lattices in $\text{SL}_2(\mathbb R)\times\text{SL}_2(\mathbb R)$ and $\text{SL}_2(\mathbb R)\times(\prod_{p\in S}\text{SL}_2(\mathbb Q_p))$, respectively. Since the same holds if $\text{SL}_2$ is replaced  with $\text{PSL}_2$, it follows that $\text{PSL}_2(\mathcal O_d)$ and $\text{PSL}_2(\mathbb Z[S^{-1}])$ belong to $\mathscr L$.
\end{enumerate}
\end{remark}

\begin{remark} 
Let $G=G_1\times...\times G_n$ be as in the first part of Theorem \ref{A} and $\Gamma<G$  be an irreducible, but not necessarily icc, lattice. Then the center $Z(\Gamma)$ of $\Gamma$ is contained in $Z(G)$  and  $\Gamma/Z(\Gamma)$ is icc.  Thus, $\Gamma/Z(\Gamma)$ is an irreducible icc lattice in $G/Z(G)$. By Remark \ref{L}(1) it follows that  $\Gamma/Z(\Gamma)\in\mathscr L$ and  Theorem \ref{A} implies that the II$_1$ factor
  $L(\Gamma/Z(\Gamma))$ is prime. \end{remark}

\begin{remark}
Let $G=G_1\times...\times G_n$ be as in the first part of Theorem \ref{A}.
Popa and Vaes proved that any lattice $\Gamma<G$ is Cartan-rigid:   any II$_1$ factor $L^{\infty}(X)\rtimes\Gamma$ arising from a free ergodic pmp action of $\Gamma$ has a unique Cartan subalgebra up to unitary conjugacy (see \cite[Theorem 1.3]{PV12}). Moreover, their proof shows that $L(\Gamma)$ does not have a Cartan subalgebra (see \cite[Section 5]{PV12}).  Both the approach of \cite{PV12} and our proof of Theorem \ref{A} use the fact that the lattices in $G$ are measure equivalent to a product of non-elementary hyperbolic groups.  However, unlike the results from \cite{PV12}, the conclusion of Theorem \ref{A} does not hold for arbitrary lattices $\Gamma<G$, as it obviously fails for product lattices $\Gamma=\Gamma_1\times...\times\Gamma_n$, whenever $n\geq 2$ and $\Gamma_i<G_i$ is a lattice for all $1\leq i\leq n$. To prove Theorem \ref{A} we will perform a detailed analysis which shows that if $\Gamma<G$ is any icc lattice such that $L(\Gamma)$ is not prime, then $\Gamma$ is a product group and thus cannot be an irreducible lattice. 
\end{remark}

The following corollary in an immediate consequence of Theorem \ref{A}.

\begin{mcor}\label{B}
Let $d\geq 2$ be a square-free integer and $S$ be a finite non-empty set of primes.

Then  $L(\text{PSL}_2(\mathcal O_d))$ and $L(\text{PSL}_2(\mathbb Z[S^{-1}]))$ are prime II$_1$ factors.
\end{mcor}

\begin{remark}  \cite[Corollary]{CdSS15} implies that $L(\text{PSL}_2(\mathbb Z[\sqrt{2}]))$ is not isomorphic to $L(\Gamma_1\times\Gamma_2)$, for any non-amenable groups $\Gamma_1$ and $\Gamma_2$ in Ozawa's class $\mathcal S$ \cite{Oz04,BO08}.
Corollary \ref{B} strengthens this fact by showing that $L(\text{PSL}_2(\mathbb Z[\sqrt{2}]))$ is prime and hence not isomorphic to $L(\Gamma_1\times\Gamma_2)$, for any non-trivial countable groups $\Gamma_1$ and $\Gamma_2$. 
\end{remark}

Theorem \ref{A} will be deduced from a general result which describes all tensor product decompositions of II$_1$ factors associated to groups that are measure equivalent to products of  hyperbolic groups. Before stating this result, let us recall the notion of measure equivalence due to Gromov \cite{Gr91}, and the construction of amplifications of II$_1$ factors. 

Two countable groups $\Gamma$ and $\Lambda$ are called {\bf measure equivalent} if there exist commuting free measure preserving actions of $\Gamma$ and $\Lambda$ on a standard measure space $(\Omega, m)$, such that the actions of $\Gamma$ and $\Lambda$ each admit a finite measure fundamental domain. 
Natural examples of measure equivalent groups are provided by pairs of lattices $\Gamma,\Lambda$ in an unimodular locally compact second countable group $G$. 
Indeed, endowing $G$ with a Haar measure $m$ and the left and right translation actions of $\Gamma$ and $\Lambda$ shows that $\Gamma$ and $\Lambda$ are measure equivalent.

If $M$ is a II$_1$ factor and $t>0$, then the {\bf amplification} $M^t$ is defined as the isomorphism class of $p(M\overline{\otimes}\mathbb B(\ell^2(\mathbb N))p$, where $p\in M\overline{\otimes}\mathbb B(\ell^2(\mathbb N))$  is a projection satisfying $(\tau\otimes\text{Tr})(p)=t$. Here, $\tau$ and $\text{Tr}$ denote the canonical traces of $M$ and $\mathbb B(\ell^2(\mathbb N))$, respectively. 
Finally, recall that if  $M=P_1\overline{\otimes}P_2$, for some II$_1$ factors $P_1$ and $P_2$, then for every $t>0$ we have a natural identification $M=P_1^t\overline{\otimes}P_2^{1/t}$. 

The following theorem is the main technical result of this paper:

\begin{main}\label{C}
Let $\Gamma$ be a countable icc group and denote $M=L(\Gamma)$. Assume that $\Gamma$ is measure equivalent to a product $\Lambda=\Lambda_1\times...\times\Lambda_n$ of $n\geq 1$ non-elementary hyperbolic groups $\Lambda_1,...,\Lambda_n$. Suppose that $M=P_1\overline{\otimes}P_2$, for some II$_1$ factors $P_1$ and $P_2$.

Then there exist a decomposition  $\Gamma=\Gamma_1\times\Gamma_2$, a partition $S_1\sqcup S_2=\{1,...,n\}$, a decomposition $M=P_1^{t}\overline{\otimes}P_2^{1/t}$, for some $t>0$, and a unitary $u\in M$  such that 
\begin{enumerate}
\item  $P_1^t=uL(\Gamma_1)u^*$ and $P_2^{1/t}=uL(\Gamma_2)u^*$.
\item  $\Gamma_1$ is measure equivalent to $\underset{j\in S_1}{\times}\Lambda_j$ and $\Gamma_2$ is measure equivalent to $\underset{j\in S_2}{\times}\Lambda_j$.

\end{enumerate}
 
\end{main}
In order to put Theorem \ref{C} in a better perspective, we first emphasize a new rigidity phenomenon that Theorem \ref{C} leads to, and then discuss several applications of it. 

 Connes'  classification of injective factors implies that no algebraic information regarding an icc amenable group $\Gamma$  can be recovered from $L(\Gamma)$ \cite{Co76}.  In sharp contrast, Theorem \ref{C} implies that for a natural and wide class of groups $\Gamma$, any 
tensor product decomposition of $L(\Gamma)$ must arise from a direct product decomposition of $\Gamma$. 
This adds to the few known instances where algebraic properties of the von Neumann algebra $L(\Gamma)$ can be promoted to algebraic properties of the group $\Gamma$.  We highlight here two recent developments in this direction: Ioana, Popa and Vaes' discovery of the first classes of ``W$^*$-superrigid'' groups \cite{IPV10} (see  Berbec and Vaes   \cite{BV12} for the only other known examples), and Chifan, de Santiago and Sinclair's ``product rigidity" theorem \cite{CdSS15}.

There are three main applications of Theorem \ref{C}.
First, we use Theorem \ref{C} to deduce Theorem \ref{A}. 
To briefly indicate how this works, let $\Gamma\in\mathscr{L}$ be an icc group. 
Then $\Gamma$ can be realized as an irreducible lattice in a locally compact group $G=G_1\times...\times G_n$ which also admits a product of non-elementary hyperbolic groups  as a lattice.  Assuming that $L(\Gamma)$ is not prime, we apply Theorem \ref{C} to conclude that $\Gamma$ decomposes as a product of infinite groups.  In the case $G_1,\dots,G_n$ are non-compact simple Lie groups with finite center, such a decomposition can be ruled out by appealing to Margulis' normal subgroup theorem (see \cite[Theorem 8.1.1]{Zi84}).
In the general case, we will show that such product decompositions do not exist by using a stronger version of Theorem \ref{C} (see Theorem \ref{general}). 

Secondly, Theorem \ref{C} allows us to prove a unique prime factorization result for tensor products of II$_1$ factors arising from irreducible lattices in products of rank one simple Lie groups.  

\begin{mcor}\label{D}
Let $\Gamma$ be a countable icc group which is measure equivalent to a product
of $n\geq 1$ non-elementary hyperbolic groups.  
Denote $M=L(\Gamma)$.

Then there exists a unique (up to permutation of factors) decomposition $\Gamma=\Gamma_1\times...\times\Gamma_k$, for some $1\leq k\leq n$, such that $L(\Gamma_i)$ is a prime II$_1$ factor, for every $1\leq i\leq k$. 
Moreover, the following hold:

\begin{enumerate}
\item If $M=P_1\overline{\otimes}P_2$, for some II$_1$ factors $P_1, P_2$, then there exist a partition $I_1\sqcup I_2=\{1,...,k\}$ and a decomposition $M=P_1^t\overline{\otimes}P_2^{1/t}$, for some $t>0$, such that $P_1^t=\overline{\otimes}_{i\in I_1}L(\Gamma_i)$ and $P_2^{1/t}=\overline{\otimes}_{i\in I_2}L(\Gamma_i)$, modulo unitary conjugacy in $M$. 
\item If $M=P_1\overline{\otimes}\dots\overline{\otimes}P_m$, for some $m\geq k$ and II$_1$ factors $P_1,...,P_m$, then $m=k$ and there exists a decomposition $M=P_1^{t_1}\overline{\otimes}...\overline{\otimes}P_k^{t_k}$ for some $t_1,...,t_k>0$ with $t_1t_2\dots t_k=1$ such that after permutation of indices and unitary conjugacy we have $L(\Gamma_i)=P_i^{t_i}$, for all $1\leq i\leq k$.
\item In (2), the assumption $m\geq k$ can be omitted if each $P_i$ is assumed to be prime. 
\end{enumerate}
\end{mcor}

Corollary \ref{D} in particular applies if $\Gamma=\Gamma_1\times...\times\Gamma_k$, where $\Gamma_i\in\mathscr L, 1\leq i\leq k$, are icc groups.

The first unique prime factorization results for II$_1$ factors were obtained by Ozawa and Popa in their pioneering work \cite{OP03}.
More precisely, \cite{OP03} established conclusions (1)-(3) of Corollary~\ref{D} for $M = L(\Gamma_1)\overline{\otimes}...\overline{\otimes}L(\Gamma_k)$ whenever $\Gamma_i$, $1 \le i \le k$, are icc non-amenable groups which are either hyperbolic or discrete subgroups (in particular, lattices) of connected simple Lie groups of rank one.  In the meantime, several other unique prime factorization results have been obtained in \cite{Pe06, CS11, SW11,Is14, CKP14, HI15, Ho15, Is16}. Corollary \ref{D} is the first unique prime factorization result that applies to II$_1$ factors coming from irreducible lattices in certain higher rank semisimple Lie groups.
It implies in particular that if $\Gamma$ is any irreducible lattice in a product of $n\geq 1$ connected non-compact rank one simple Lie groups with finite center and $\Gamma_0=\Gamma/Z(\Gamma)$, then the II$_1$ factors $$L(\Gamma_0)^{\overline{\otimes}k}=\underbrace{L(\Gamma_0)\overline{\otimes}...\overline{\otimes}L(\Gamma_0)}_{\text{$k$ times}},\;\; \;k\geq 1,$$ are pairwise non-isomorphic. This generalizes the case $n=1$ obtained in \cite{CH88}, for lattices in the simplectic groups $\text{Sp}(m,1)$, and in \cite{OP03}, for lattices in arbitrary connected non-compact rank one simple Lie groups with finite center.

Our last application of Theorem \ref{C} relates to prime factorization for measure equivalence.  
In \cite{MS02}, Monod and Shalom proved a series of striking rigidity results for orbit and measure equivalence. In particular, they also studied groups $\Gamma$ which are measure equivalent to a product $\Lambda=\Lambda_1\times...\times\Lambda_n$ of non-elementary hyperbolic groups (more generally, of groups in the class $\mathcal C_{\text{reg}}$). In this context, they proved a prime factorization result: if $\Gamma=\Gamma_1\times...\times\Gamma_m$ is itself a product group and all  the groups involved are torsion-free,  then $m\leq n$, and if $m=n$ then, after permutation of the indices, $\Gamma_i$ is measure equivalent to $\Lambda_i$, for $1\leq i\leq n$ (see \cite[Theorem 1.16]{MS02} and \cite[Theorem 3]{Sa09}).  Theorem \ref{C} recovers and strengthens this result in the case $\Gamma$ is icc and $\Lambda_i$ are hyperbolic. More precisely, it implies that if instead of assuming that $\Gamma$ is a product of $m$ infinite groups, one merely requires that $L(\Gamma)$ is a tensor product of $m$ II$_1$ factors, then $m\leq n$, and if $m=n$, then there exists a unique product decomposition $\Gamma=\Gamma_1\times...\times\Gamma_m$ such that the above conclusion holds. 

\subsection{Comments on the proof of Theorem \ref{C}} Since all of our main results are deduced from Theorem \ref{C}, we outline briefly and informally its method of proof. 
Let $\Gamma$ be an icc group which is measure equivalent to a product $\Lambda=\Lambda_1\times...\times\Lambda_n$ of non-elementary hyperbolic groups.  By \cite{Fu99}, $\Gamma$ and $\Lambda$ must have stably orbit equivalent actions. To simplify notation, assume that $\Gamma$ and $\Lambda$ admit in fact orbit equivalent actions, i.e. there exist free ergodic pmp actions of $\Gamma$ and $\Lambda$ on a probability space $(X,\mu)$ whose orbits are equal, almost everywhere.
Denote $M=L^{\infty}(X)\rtimes\Gamma=L^{\infty}(X)\rtimes\Lambda$.

 Our goal is to classify all tensor product decompositions $L(\Gamma)=P_1\overline{\otimes}P_2$.
 To achieve this goal, we use a combination of techniques from Popa's deformation/rigidity theory.  
 
 First, we use repeatedly the relative strong solidity property of hyperbolic groups (see Section \ref{SS: RSS}) established in the breakthrough work \cite{PV11,PV12}, to conclude the existence of a partition $S_1\sqcup S_2=\{1,\dots,n\}$ such that letting $\Lambda_{S_i} = \underset{j\in S_i}{\times}\Lambda_j$ for $i \in \{1, 2\}$, we have
 \begin{align} \label{E: RSSimbed}
 P_1 \prec L^{\infty}(X)\rtimes \Lambda_{S_1}, \quad\text{ and }\quad P_2 \prec L^{\infty}(X)\rtimes \Lambda_{S_2},
 \end{align}
 where $P\prec Q$ denotes that a corner of $P$ embeds into $Q$ inside the ambient algebra, in the sense of Popa \cite{Po03}. For simplicity below, we will write $P\sim Q$ to indicate that $Pp'\prec Q$ and $Qq'\prec P$, for all non-zero projections $p'$ and $q'$ in the relative commutants of $P$ and $Q$.

To see the importance of \eqref{E: RSSimbed}, note that for each $i$, we have $P_i \subset L(\Gamma) \subset L^\infty(X) \rtimes \Gamma$, and in this sense $P_i$ is ``far away" from $L^\infty(X)$. This remains true after passing through the intertwining in \eqref{E: RSSimbed}, and hence one thinks of the image of $P_i$ as being not far from $L(\Lambda_{S_i})$ in $L^{\infty}(X)\rtimes \Lambda_{S_i}$.
The critical consequence of \eqref{E: RSSimbed} is then that it allows one to show that each $P_i$ inherits a weaker form of the relative strong solidity present in $L(\Lambda_{S_i})$.

In particular, if we follow \cite{PV09} and consider the comultiplication $*$-homomorphism $\Delta:M\rightarrow M\overline{\otimes}L(\Gamma)$  given by $\Delta(au_g)=au_g\otimes u_g$, for all $a\in L^{\infty}(X), g\in\Gamma$, it allows us to conclude that
\begin{equation}\label{intro1}
\Delta(L^{\infty}(X)\rtimes \Lambda_{S_1})\prec M\overline{\otimes}P_1, \quad \text{ and } \quad \Delta(L^{\infty}(X)\rtimes \Lambda_{S_2})\prec M\overline{\otimes}P_2.
\end{equation}
 
 This is achieved in the first part of Section 5.
 The conclusion \eqref{intro1} enables us to then make crucial use of an ultrapower technique from \cite{Io11} (see Section 4) in combination with the transfer of rigidity principle from \cite{PV09} to find subgroups $\Sigma_1, \Sigma_2 <\Gamma$ such that
 \begin{equation}\label{intro2}
 L^{\infty}(X)\rtimes\Sigma_1\sim L^{\infty}(X)\rtimes\Lambda_{S_1}, \quad \text{ and } \quad L^{\infty}(X)\rtimes\Sigma_2\sim L^{\infty}(X)\rtimes\Lambda_{S_2};
 \end{equation}
 \begin{equation}
 \label{intro3}P_1\sim L(\Sigma_1), \quad \text{ and } \quad P_2\sim L(\Sigma_2).
 \end{equation}
 
 This is achieved in the second part of Section 5.
 
 We then use \eqref{intro2} to deduce that $\Sigma_i$ is measure equivalent to $\Lambda_{S_i}$, for all $i\in\{1,2\}$ (see Section 3). 
 
Finally, inspired by results in \cite{CdSS15}, we show that \eqref{intro3} implies that, after replacing $\Sigma_i$ with a commensurable subgroup $\Gamma_i<\Gamma$ we have $\Gamma=\Gamma_1\times\Gamma_2$ with $P_i=L(\Gamma_i)$ for all $i\in\{1,2\}$, modulo unitary conjugacy and amplification (see Section 6).
This altogether proves Theorem \ref{C}.
 
\subsection{Organization of the paper}  Besides the introduction and a section of preliminaries, this paper has five other sections:
Sections 3-6 are devoted to the different ingredients of the proof of Theorem \ref{C}, as explained above. In Section 7, we finalize the proof of Theorem \ref{C} and derive the rest of our main results.

\subsection{Acknowledgment} We are grateful to the referee for many comments that helped improve the exposition.

\section{Preliminaries}
\subsection{Terminology} We fix notation regarding tracial von Neumann algebras, countable groups, and countable pmp equivalence relations.

A {\it tracial von Neumann algebra} 
$(M,\tau)$ is a von Neumann algebra $M$ equipped with a faithful normal tracial state $\tau$. We denote by $L^2(M)$ the completion of $M$ with respect to the norm $\|x\|_2=\sqrt{\tau(x^*x)}$ and consider the standard representation $M\subset\mathbb B(L^2(M))$. 
Unless stated otherwise, we will always assume that $M$ is separable, i.e.  $L^2(M)$ is a separable Hilbert space.
For a set $\mathcal S\subset\mathbb B(L^2(M))$, we denote by $\mathcal S'$ its commutant. If $\mathcal S$ is closed under adjoint, then by von Neumann's double commutant theorem, $\mathcal S''=(\mathcal S')'$ is exactly the von Neumann algebra generated by $\mathcal S$.
We denote by $\mathcal U(M)$ the group of unitary elements of $M$, by $(M)_1=\{x\in M\;|\;\|x\|\leq 1\}$ the unit ball of $M$, and by $\mathcal Z(M)=M\cap M'$ the center of $M$.

Let $P\subset M$ be a von Neumann subalgebra, which we  will always assume to be unital. We denote by $e_P:L^2(M)\rightarrow L^2(P)$ the orthogonal projection onto $L^2(P)$, by $E_P:M\rightarrow P$ the conditional expectation onto $P$,  and by $\mathcal N_{M}(P)=\{u\in\mathcal U(M)\;|\;uPu^*=P\}$
the {\it normalizer of $P$ in $M$}. The subalgebra $P\subset M$ is called {\it regular} if $\mathcal N_M(P)''=M$.
{\it Jones' basic construction} of the inclusion $P\subset M$ is defined as the von Neumann subalgebra of $\mathbb B(L^2(M))$ generated by $M$ and $e_P$, and is denoted by $\langle M,e_P\rangle$. If $J:L^2(M)\rightarrow L^2(M)$ denotes the involution given by $J(x)=x^*$, for every $x\in M$, then $\langle M,e_P\rangle= (JPJ)'\cap \mathbb B(L^2(M))$.

For a countable group $\Gamma$,  its left regular representation $u:\Gamma\rightarrow\mathcal U(\ell^2(\Gamma))$ is given by $u_g(\delta_h)=\delta_{gh}$, where $\{\delta_h \,|\, h\in\Gamma\}$ denotes the usual orthonormal basis of $\ell^2(\Gamma)$. The weak operator closure of $\{u_g \,|\, g\in\Gamma\}$ is a tracial von Neumann algebra which is called the {\it group von Neumann algebra of $\Gamma$}, and denoted by $L(\Gamma)$. This algebra is a II$_1$ factor if and only if $\Gamma$ is {\it icc}, i.e. every non-trivial conjugacy class of $\Gamma$ is infinite. Let $S, T\subset\Gamma$ be two subsets. We denote by $\langle S\rangle$ the group generated by $S$, and by $C_S(T)=\{g\in S \,|\, gh=hg,\;\text{for all $h\in T$}\}$ the {\it centralizer of $T$ in $S$}.

For a pmp action $\Gamma\curvearrowright (X,\mu)$ of a countable group $\Gamma$ on a standard probability space $(X,\mu)$, we denote by $\mathcal R(\Gamma\curvearrowright X)=\{(x,y)\in X\times X \,|\, \Gamma\cdot x=\Gamma\cdot y\}$ the associated orbit equivalence relation. For a countable pmp equivalence relation $\mathcal R$ on $(X,\mu)$ and a measurable subset $Y\subset X$, we denote by $\mathcal R|_{Y}=\mathcal R\cap (Y\times Y)$ the restriction of $\mathcal R$ to $Y$.  For every $x\in X$, $[x]_{\mathcal R}$ denotes its equivalence class. We  denote by $[[\mathcal R]]$ the set of partially defined measurable isomorphisms $\theta:Y=\text{dom}(\theta)\rightarrow Z=\text{ran}(\theta)$ between measurable subsets $Y,Z\subset X$ which satisfy $(\theta(x),x)\in\mathcal R$, for almost every $x\in Y$.  The group of measurable isomorphisms $\theta:X\rightarrow X$ which satisfy $(\theta(x),x)\in\mathcal R$, for almost every $x\in X$, is called the {\it full group} of $\mathcal R$ and denoted by $[\mathcal R]$.

Finally, two pmp actions $\Gamma\curvearrowright (X,\mu)$ and $\Lambda\curvearrowright (Y,\nu)$ are called {\it stably orbit equivalent (SOE)} if there exist
non-negligible measurable subsets $X_0\subset X$ and $Y_0\subset Y$, and a measure preserving isomorphism $\theta:(X_0,\mu(X_0)^{-1}\mu|_{X_0})\rightarrow (Y_0,\nu(Y_0)^{-1}\nu|_{Y_0})$ such that $(\theta\times\theta)(\mathcal R(\Gamma\curvearrowright X)|_{X_0})=\mathcal R(\Lambda\curvearrowright Y)|_{Y_0}$.
If this holds for $X_0=X$ and $Y_0=Y$, the actions are called {\it orbit equivalent (OE)}.

\subsection {Intertwining-by-bimodules} We next recall from  \cite [Theorem 2.1 and Corollary 2.3]{Po03} the powerful {\it intertwining-by-bimodules} technique of Popa.

\begin {theorem}[\!\!\cite{Po03}]\label{corner} Let $(M,\tau)$ be a tracial von Neumann algebra and $P\subset pMp, Q\subset qMq$ be von Neumann subalgebras. Let 
 $\mathcal U\subset\mathcal U(P)$ be a subgroup such that $\mathcal U''=P$.

Then the following are equivalent:

\begin{itemize}

\item There exist projections $p_0\in P, q_0\in Q$, a $*$-homomorphism $\theta:p_0Pp_0\rightarrow q_0Qq_0$  and a non-zero partial isometry $v\in q_0Mp_0$ such that $\theta(x)v=vx$, for all $x\in p_0Pp_0$.

\item There is no sequence $u_n\in\mathcal U$ satisfying $\|E_Q(x^*u_ny)\|_2\rightarrow 0$, for all $x,y\in pMq$.
\end{itemize}

If one of these equivalent conditions holds true,  then we write $P\prec_{M}Q$, and say that {\it a corner of $P$ embeds into $Q$ inside $M$.}
If $Pp'\prec_{M}Q$ for any non-zero projection $p'\in P'\cap pMp$, then we write $P\prec^{s}_{M}Q$.
\end{theorem}

\begin{remark}
In the context of Theorem \ref{corner}, let $(\tilde M,\tilde{\tau})$ be a tracial von Neumann algebra such that $M\subset\tilde M$ and $\tilde{\tau}_{|M}=\tau$. If $P\prec_{M}Q$, then clearly $P\prec_{\tilde M}Q$. But the following fact, which we will use throughout the paper, also holds: if $P\prec^s_MQ$, then $P\prec^s_{\tilde M}Q$. To see this, assume that $P\prec^s_MQ$ and let $p'\in P'\cap\tilde M$ be a non-zero projection. Let $p''\in P'\cap M$ be the support projection of $E_M(p')$. Since $Pp''\prec_{M}Q$, we can find projections $p\in P, q\in Q$, a $*$-homomorphism $\theta:pPpp''\rightarrow qQq$, and a non-zero partial isometry $v\in qMpp''$ such that $\theta(x)v=vx$, for all $x\in pPpp''$. Let $\tilde\theta:pPpp'\rightarrow qQq$ be the $*$-homomorphism given by $\tilde\theta(xp')=\theta(xp'')$, for all $x\in pPp$, and put $\tilde v=vp'$. Then $\tilde\theta$ is well-defined and $\tilde\theta(y)\tilde v=\tilde v y$, for all $y\in pPpp'$. Since $\tilde v\not=0$, we get that $Pp'\prec_{\tilde M}Q$.
\end{remark}

\begin{remark}\label{corner2} Let $P$ and $Q_i\subset M_i$, $1\leq i\leq m$, be tracial von Neumann algebras. Let  $\mathcal U\subset\mathcal U(P)$ a subgroup such that $\mathcal U''=P$, and $\pi_i:P\rightarrow M_i$, $1\leq i\leq m$, be  trace preserving $*$-homomorphisms. Assume that there exist $\delta>0$ and $a_i,b_i\in M_i$, $1\leq i\leq m$, such that $\sum_{i=1}^m\|E_{Q_i}(a_i\pi_i(u)b_i)\|_2^2\geq\delta$, for all $u\in\mathcal U.$ Then $\pi_i(P)\prec_{M_i}Q_i$, for some $1\leq i\leq m$.
Indeed, the above inequality implies that a corner of the von Neumann algebra generated by $\{\oplus_{i=1}^m\pi_i(u) \,|\, u\in\mathcal U\}$ embeds into $\oplus_{i=1}^mQ_i$ inside $\oplus_{i=1}^{m}M_i$, which implies the desired conclusion (see \cite[proof of Theorem 4.3]{IPP05} for details).

\end{remark}

We continue with the following two lemmas, in which we collect several elementary facts.

\begin{lemma}\label{facts1}
Let $(M,\tau)$ be a tracial von Neumann algebra and $P\subset pMp, Q\subset qMq, R\subset rMr$ be von Neumann subalgebras. 
Then the following hold:

\begin{enumerate}
\item \emph{\cite{Va08}} Assume that $P\prec_{M}Q$ and $Q\prec_{M}^sR$. Then $P\prec_{M}R$.
\item Assume that $Pz\prec_MQ$, for any non-zero projection $z\in \mathcal N_{pMp}(P)'\cap pMp\subset\mathcal Z(P'\cap pMp)$. Then $P\prec_M^sQ$.
\item Assume that $P\prec_{M}Q$.  Then there is a non-zero projection $z\in\mathcal N_{pMp}(P)'\cap pMp$ such that $Pz\prec^s_{M}Q$.
\item Assume that $P\prec_{M}Q$. Then there is a non-zero projection $z\in\mathcal N_{qMq}(Q)'\cap qMq$ such that $P\prec_MQq'$, for every non-zero projection $q'\in Q'\cap M$ with $q'\leq z$.
\end{enumerate}
\end{lemma}

{\it Proof.} (1) This part is precisely \cite[Lemma 3.7]{Va08}.

(2) $\&$ (3) Using Zorn's lemma and a maximality argument, we can find a projection $z\in P'\cap pMp$ such that $Pz\prec_M^sQ$ and $P(p-z)\nprec_MQ$. 
We claim that $z\in\mathcal N_{pMp}(P)'\cap pMp$. This claim clearly implies both (2) and (3).

Let us first show that
 $z\in\mathcal Z(P'\cap pMp)$. Otherwise, we can find non-zero equivalent projections $p_1,p_2\in P'\cap pMp$ satisfying $p_1\leq z, p_2\leq p-z$. Let $u\in\mathcal U(P'\cap pMp)$ such that $up_1u^*=p_2$. Then $uPp_1u^*=Pp_2$, which contradicts that $Pp_1\prec_MQ$, while $Pp_2\nprec_MQ$. This shows that indeed $z\in\mathcal Z(P'\cap pMp)$. Now, if $u\in\mathcal N_{pMp}(P)$, then $uzu^*\in\mathcal Z(P'\cap pMp)$ and $Puzu^*=uPzu^*\prec_M^sQ$. The maximality property of $z$ forces $uzu^*\leq z$, hence $uzu^*=z$. This proves the claim.

(4) Let $p_0\in P, q_0\in Q$ be projections, $\theta:p_0Pp_0\rightarrow q_0Qq_0$ a $*$-homomorphism, and  $v\in q_0Mp_0$ a non-zero partial isometry such that $\theta(x)v=vx$, for all $x\in p_0Pp_0$. Let $r\in Q'\cap qMq$ be the support projection of $E_{Q'\cap qMq}(vv^*)$. 
Let $r'\in Q'\cap qMq$ be a non-zero projection with $r'\leq r$. Let $\psi:p_0Pp_0\rightarrow q_0r'(Qr')q_0r'$ be given by $\psi(x)=\theta(x)r'$ and $w=r'v\in q_0r'Mp_0$. Then $\psi(x)w=wx$, for all $x\in p_0Pp_0$. Since $E_{Q'\cap qMq}(wv^*)=r'E_{Q'\cap qMq}(vv^*)\not=0$, we get that $w\not=0$, hence $P\prec_{M}Qr'$. 

Let $z'\in\mathcal Z(Q'\cap qMq)$ be the central support of $r$, and put $z=\vee_{u\in\mathcal N_{qMq}(Q)}uz'u^*\in\mathcal N_{qMq}(Q)'\cap qMq$. If $q'\in Q'\cap qMq$ is a non-zero projection with $q'\leq z$,  we can find $u\in\mathcal N_{qMq}(Q)$ such that $q'uz'u^*\not=0$. This implies the existence of non-zero equivalent projections $q'', r'\in Q'\cap qMq$ such that $q''\leq q'$ and $r'\leq uru^*$. As $u^*r'u\leq r$, we get that $P\prec_{M}Qu^*r'u=u^*(Qr')u$, hence $P\prec_{M}Qr'$. 
This implies that $P\prec_{M}Qq''$ and since $q''\leq q'$, we derive that $P\prec_{M}Qq'$. This finishes the proof.
\hfill$\blacksquare$

\begin{lemma}\label{groups}
Let $\Lambda<\Gamma$ be a countable groups. 
\begin{enumerate}
\item If $L(\Gamma)\prec_{L(\Gamma)}L(\Lambda)$, then $\Lambda$ has finite index in $\Gamma$.
\item If $\Lambda$ has finite index in $\Gamma$, then  $L(\Gamma)\prec^s_{L(\Gamma)}L(\Lambda)$.
\end{enumerate}
\end{lemma}

{\it Proof.} (1) Assume that $L(\Gamma)\prec_{L(\Gamma)}L(\Lambda)$. Thus,  the $L(\Gamma)$-$L(\Gamma)$-bimodule $L^2(\langle L(\Gamma),e_{L(\Lambda)}\rangle)$ contains a non-zero $L(\Gamma)$-central vector (see \cite[Theorem 2.1]{Po03}).  Therefore, the unitary representation $\pi:\Gamma\rightarrow\mathcal U(L^2(\langle L(\Gamma),e_{L(\Lambda)}\rangle))$ given by $\pi(g)(\xi)=u_g\xi u_g^*$, has a non-zero invariant vector.  As $\pi$ is isomorphic to a subrepresentation of the representation $\Gamma\curvearrowright\oplus_{k\in\Gamma}\ell^2(\Gamma/k\Lambda k^{-1})$, we deduce that $\ell^2(\Gamma/k\Lambda k^{-1})$ contains a non-zero $\Gamma$-invariant vector, for some $k\in\Gamma$. This implies that $k\Lambda k^{-1}$ and hence $\Lambda$ has finite index in $\Gamma$.

(2) Assume that $[\Gamma:\Lambda]<\infty$. Let $g_1,...,g_m\in\Gamma$ such that $\Gamma$ is the disjoint union of $\{g_i\Lambda\}_{i=1}^m$. Fix any non-zero projection $p\in L(\Gamma)' \cap L(\Gamma) = \mathcal Z(L(\Gamma))$. Then $0<\|p\|_2^2=\|up\|_2^2=\sum_{i=1}^m\|E_{L(\Lambda)}(u_{g_i}^*up)\|_2^2$, for every $u\in\mathcal U(L(\Gamma))$. This shows that $L(\Gamma)p\prec_{L(\Gamma)}L(\Lambda)$, and the conclusion follows.
\hfill$\blacksquare$

\subsection {Relative amenability} A tracial von Neumann algebra $(M,\tau)$ is called {\it amenable} if there exists a positive linear functional $\varphi:\mathbb B(L^2(M))\rightarrow\mathbb C$ such that $\varphi_{|M}=\tau$ and $\varphi$ is $M$-{\it central}, in the following sense:  $\varphi(xT)=\varphi(Tx)$, for all $x\in M$ and $T\in\mathbb B(L^2(M))$. By 
Connes' celebrated classification of injective factors \cite{Co76}, $M$ is amenable iff it is approximately finite dimensional.

Throughout the paper, we make extensive use of the notion of relative amenability introduced by Ozawa and Popa. 
Let $p\in M$ be a projection, and
 $P\subset pMp, Q\subset M$ be von Neumann subalgebras. 
Following \cite[Section 2.2]{OP07} we say that $P$ is  {\it amenable relative to $Q$ inside $M$} if there exists a positive linear functional $\varphi:p\langle M,e_Q\rangle p\rightarrow\mathbb C$ such that $\varphi_{|pMp}=\tau$ and $\varphi$ is $P$-central.

{\bf Convention.} Whenever the ambient algebra $(M,\tau)$ is clear from the context, we will write $P\prec Q$ instead of $P\prec_{M}Q$. We will also say that {\it $P$ is amenable relative to $Q$} instead of {\it $P$ is amenable relative to $Q$ inside $M$}.

We continue with two lemmas containing several elementary facts regarding relative amenability.

\begin{lemma}\label{facts2}
Let $(M,\tau)$ be a tracial von Neumann algebra,  and $P\subset pMp, Q\subset M$ be von Neumann subalgebras. Then the following hold:

\begin{enumerate}
\item Assume that $P$ is amenable relative to $Q$. Then $Pp'$ is amenable relative to $Q$, for every projection $p'\in P'\cap pMp$.
\item Assume that $p_0Pp_0p'$ is amenable relative to $Q$, for some projections $p_0\in P, p'\in P'\cap pMp$. Let $z$ be the smallest projection belonging to $\mathcal N_{pMp}(P)'\cap pMp$  such that $p_0p'\leq z$. Then $Pz$ is amenable relative to $Q$.

\item Assume that $P\prec^s_{M}Q$. Then $P$ is amenable relative to $Q$.
\end{enumerate}
\end{lemma}

{\it Proof.} (1)  Let $\varphi:p\langle M,e_Q\rangle p\rightarrow\mathbb C$ be a $P$-central positive linear functional such that $\varphi_{|pMp}=\tau$. Then the restriction of $\varphi$ to $p'\langle M,e_Q\rangle p'$ witnesses that $Pp'$ is amenable relative to $Q$.

(2) Let $p''\in\mathcal Z(P'\cap pMp)$ be the smallest projection such that $p_0p'\leq p''$. By \cite[Remark 2.2]{Io12a}, $Pp''$ is 
amenable relative to $Q$. Since $z=\vee_{u\in\mathcal N_{pMp}(P)}up''u^*$, we can find  $p_n\in\mathcal Z(P'\cap pMp)p''$ and $u_n\in\mathcal N_{pMp}(P)$ such that $z=\sum_n u_np_nu_n^*$. Since $Pz\subset\bigoplus_n u_nPp_nu_n^*$ and $Pp_n$ is amenable relative to $Q$ for every $n$ by part (1),  it follows that $Pz$ is amenable relative to $Q$.

(3) If $P$ is not amenable relative to $Q$, then there is a non-zero projection $z\in\mathcal Z(P'\cap pMp)$ such that $Pz'$ is not amenable relative to $Q$, for any non-zero projection $z'\in\mathcal Z(P'\cap pMp)z$. Since $Pz\prec_{M}Q$, we can find projections 
 $p_0\in P$, $q_0\in Q$, a $*$-homomorphism $\theta:p_0Pp_0z\rightarrow q_0Qq_0$, and a non-zero partial isometry $v\in q_0Mp_0z$ such that $\theta(x)v=vx$, for all $x\in p_0Pp_0z$. Then $v^*v=p_0p'$, for a projection $p'\in (P'\cap pMp)z$, and $vv^*\in\theta(p_0Pp_0z)'\cap q_0Qq_0$. Since $\theta(p_0Pp_0z)\subset q_0Qq_0$, by part (1), $\theta(p_0Pp_0)vv^*$ is amenable relative to $Q$. Since $\theta(p_0Pp_0)vv^*$ is unitarily conjugate to $p_0Pp_0p'$, the latter algebra is also amenable relative to $Q$. By part (2), we can find a projection $z'\in\mathcal Z(P'\cap pMp)$  such that $p_0p'\leq z'$ and $Pz'$ is amenable relative to $Q$. This contradicts the definition of $z$.
\hfill$\blacksquare$

\begin{lemma}\label{union}
Let $(M,\tau)$ be a tracial von Neumann algebra and $Q\subset M$ a von Neumann subalgebra. Let $(P_i)_{i\in I}\subset pMp$ be an increasing net of von Neumann subalgebras, and denote $P=(\cup_{i\in I}P_i)''$. 

If $P_i$ is amenable relative to $Q$, for every $i\in I$, then $P$ is amenable relative to $Q$.
\end{lemma}

{\it Proof.} Let $\lim_{i}$ denote a state on $\ell^{\infty}(I)$ which extends the usual limit.
For every $i\in I$, let $\varphi_i:p\langle M,e_Q\rangle p\rightarrow\mathbb C$ be a $P_i$-central positive linear functional such that ${\varphi_i}_{|pMp}=\tau$. We define $\varphi:p\langle M,e_Q\rangle p\rightarrow\mathbb C$ by letting 
$\varphi(T)=\lim_i\varphi_i(T)$, for every $T\in p\langle M,e_Q\rangle p$.

Then $\varphi$ is a positive linear functional and $\varphi_{|pMp}=\tau$. Moreover, $\varphi$ is $P_i$-central, for every $i\in I$. To see this, let $x\in P_i$, for some $i\in I$, and $T\in p\langle M,e_Q\rangle p$. If $j\in I$ satisfies $j\geq i$, then $P_i\subset P_j$ and thus $x\in P_j$. Hence, we have $\varphi_j(xT)=\varphi_j(Tx)$, for every $j\geq i$, which implies that $\varphi(xT)=\varphi(Tx)$.

Let $A\subset P$ be the set of all $x\in P$ such that  $\varphi(xT)=\varphi(Tx)$, for every $T\in p\langle M,e_Q\rangle p$. By the above, $A$ contains $\cup_{i\in I}P_i$. On the other hand, the Cauchy-Schwarz inequality implies that $|\varphi(xT)|\leq\sqrt{\varphi(x^*x)\varphi(TT^*)}\leq \|x\|_2 \|T\|$ and similarly that $|\varphi(Tx)|\leq \|x\|_2\|T\|$, for all $x\in pMp$ and $T\in p\langle M,e_Q\rangle p$. This implies that $A$ is closed in $\|.\|_2$. Hence, $A=P$ and  thus $\varphi$ is $P$-central. 
\hfill$\blacksquare$

We next record the following useful result:

\begin{lemma}\label{PV11}
Let $(M,\tau)$ be a tracial von Neumann algebra and $Q_1, Q_2\subset M$ von Neumann subalgebras which form a commuting square, i.e. $E_{Q_1}\circ E_{Q_2}=E_{Q_2}\circ E_{Q_1}$. Assume that $Q_1$ is regular in $M$.
Let $P\subset pMp$ be a von Neumann subalgebra. Then the following hold:

\begin{enumerate} 

\item \emph{\cite{PV11}} If $P$ is amenable relative to $Q_1$ and $Q_2$, then $P$ is amenable relative to $Q_1\cap Q_2$. 
\item If $P\prec_M^s Q_1$ and $P\prec_M^s Q_2$, then $P\prec_{M}^sQ_1\cap Q_2$.
\end{enumerate} 
\end{lemma}

{\it Proof.}
Part (1) is precisely \cite[Proposition 2.7]{PV11}. Part (2) follows easily by adapting the proof of \cite[Proposition 2.7]{PV11}. For completeness we include a proof, using the notation therein. 

Assume that $P\prec_M^s Q_1$ and $P\prec_M^s Q_2$. Let $p'\in P'\cap pMp$ be a non-zero projection. We will prove the conclusion of part (2) by showing that $Pp'\prec_{M}Q_1\cap Q_2$.
To this end, for $i\in\{1,2\}$, we let $\text{Tr}_i:\langle M,e_{Q_i}\rangle\rightarrow\mathbb C$ be the canonical semifinite trace given by $\text{Tr}_i(xe_{Q_i}y)=\tau(xy)$, for all $x,y\in M$. Let  $\mathcal T_i:L^1(\langle M,e_{Q_i})\rangle\rightarrow L^1(M)$ given by $\tau(\mathcal T_i(T)x)=\text{Tr}_i(Tx)$, for all $T\in L^1(\langle M,e_{Q_i}\rangle)$ and $x\in M$.

Since $P\prec_M^sQ_1$, we get that $Pp'\prec_MQ_1$. By \cite[Theorem 2.1]{Po03} we can find a non-zero projection $e_1\in (Pp')'\cap p'\langle M,e_{Q_1}\rangle p'$ such that $\text{Tr}_1(e_1)<\infty$. Let $p''\in M$ be the support projection of $\mathcal T_1(e_1)$. Then $p''\in (Pp')'\cap p'Mp'=p'(P'\cap pMp)p'$. Since $P\prec_M^sQ_2$, we get that $Pp''\prec_MQ_2$.  By \cite[Theorem 2.1]{Po03} we can find a non-zero projection $e_2\in (Pp'')'\cap p''\langle M,e_{Q_2}\rangle p''$ with $\text{Tr}_2(e_2)<\infty$.

Next, consider the $M$-$M$-bimodule $\mathcal H=L^2(\langle M,e_{Q_1}\rangle)\otimes_{M}L^2(\langle M,e_{Q_2})\rangle$ and put $\xi=e_1\otimes_Me_2\in\mathcal H$. Then $x\xi=\xi x$, for all $x\in Pp''$. Moreover, since $p''e_2=e_2\not=0$ and $p''$ is the support projection of $\mathcal T_1(e_1)$ we have $\mathcal T_1(e_1)^{{1}/{2}}e_2\not=0$. Since $\|\xi\|^2=\langle e_1\otimes_Me_2,e_1\otimes_Me_2\rangle=\langle\mathcal T_1(e_1)e_2,e_2\rangle=\|\mathcal T_1(e_1)^{{1}/{2}}e_2\|^2$, we deduce that $\xi\not=0$. 

Now, by the last part of the proof of \cite[Proposition 2.7]{PV11}, the $M$-$M$-bimodule $\mathcal H$ is contained in a multiple of $_ML^2(\langle M,e_Q\rangle)_M$, where $Q=Q_1\cap Q_2$. Since $0\not=\xi=p''\xi p''$, we derive the existence of a non-zero vector $\eta\in p''L^2(\langle M,e_Q\rangle)p''$ such that $x\eta=\eta x$, for all $x\in Pp''$. Then $\zeta=\eta^*\eta\in L^1(\langle M,e_Q\rangle)_{+}$ satisfies $0\not=\zeta=p''\zeta p''$ and $x\zeta=\zeta x$, for all $x\in Pp''$.
Let $t>0$ such that the spectral projection $f={\bf 1}_{[t,\infty)}(\zeta)$ is non-zero. Then $f\in (Pp'')'\cap p''\langle M,e_Q\rangle p''$ and  since $tf\leq \zeta$, we get that $\text{Tr}(f)\leq\text{Tr}(\zeta)/{t}<\infty$, where $\text{Tr}:\langle M,e_Q\rangle\rightarrow\mathbb C$ denotes the canonical semifinite trace. 
By \cite[Theorem 2.1]{Po03} we conclude that $Pp''\prec_MQ$ and hence that $Pp'\prec_{M}Q$, as desired.
\hfill$\blacksquare$

For the last result of this subsection, assume the following context: let $\Gamma\curvearrowright (X,\mu)$, $\Lambda\curvearrowright (Y,\nu)$ be stably orbit equivalent free ergodic pmp actions. Thus, there exists $\ell\ge 1$ such that we can view $X$ as a subset of $Y\times \mathbb Z/\ell\mathbb Z$  satisfying $\mathcal R(\Gamma\curvearrowright X)=\mathcal R(\Lambda\times\mathbb Z/\ell\mathbb Z\curvearrowright Y\times\mathbb Z/\ell\mathbb Z)|_{X}$, where $\mathbb Z/\ell\mathbb Z$ acts on itself by addition. Hence,
$L^\infty(X)\rtimes\Gamma=pMp$, where $M=L^{\infty}(Y\times\mathbb Z/\ell\mathbb Z)\rtimes(\Lambda\times\mathbb Z/\ell\mathbb Z)$ and $p=1_X$.  If $B=L^\infty(Y)\otimes \mathbb M_{\ell}(\mathbb C)$, then we identify $M=B\rtimes\Lambda$, where $\Lambda$ acts trivially on $\mathbb M_{\ell}(\mathbb C)$.
 
\begin{lemma}\label{coamenable}
Let $\Sigma$ be a subgroup of $\Lambda$ such that $L(\Gamma)$ is amenable relative to $B\rtimes\Sigma$ inside $M$.

 Then $\Sigma$ is co-amenable in $\Lambda$, i.e. $\ell^{\infty}(\Lambda/\Sigma)$ admits a left $\Lambda$-invariant state.
\end{lemma}

{\it Proof.} 
Assume first that $(\nu\times c)(X)\geq 1$, where $c$ denotes the counting measure on $\mathbb Z/\ell\mathbb Z$. Then by using the ergodicity of the actions, we may assume that the inclusion $X\subset Y\times\mathbb Z/\ell\mathbb Z$ satisfies  $Y\subset X$, where $Y$ denotes its copy $Y\times\{0\}\subset Y\times\mathbb Z/\ell\mathbb Z$. Thus, we have $q=1_{Y}\leq p=1_X$ and $qBq=L^{\infty}(Y)$.
Put $A=L^\infty(X)$ and note that $A\rtimes\Gamma=pMp$ and $L^\infty(Y)\rtimes\Lambda=q(A\rtimes\Gamma)q.$ We also denote by $\{u_g\}_{g\in\Gamma}\subset A\rtimes\Gamma$ and $\{v_h\}_{h\in\Lambda}\subset B\rtimes\Lambda$ the canonical unitaries implementing the actions of $\Gamma$ and $\Lambda$ on $A$ and $B$, respectively.    We end this paragraph by observing that $\{v_h q\}_{h\in \Lambda}\subset L^{\infty}(Y)\rtimes\Lambda$ are precisely the canonical unitaries which implement the action of $\Lambda$ on $L^\infty(Y)$.

Next, since $L(\Gamma)$ is amenable relative to $B\rtimes\Sigma$ inside $M$, there exists a positive linear functional $\varphi:p\langle M,e_{B\rtimes\Sigma} \rangle p\to \mathbb C$ which is $L(\Gamma)$-central and satisfies $\varphi_{|pMp}=\tau.$

Let $\mathcal D\subset q\langle M, e_{B\rtimes\Sigma}\rangle q$ be the von Neumann algebra generated by $\{v_h q e_{B\rtimes\Sigma} \;qv_h^*\}_{h\in\Lambda}$.
If $h\in\Lambda\setminus\Sigma$, then $e_{B\rtimes\Sigma} v_h q e_{B\rtimes\Sigma}=E_{B\rtimes\Sigma}(v_h)qe_{B\rtimes\Sigma}=0$. On the other hand, if $h\in\Sigma$, then $v_hq\in B\rtimes\Sigma$, thus $v_hqe_{B\rtimes\Sigma}=e_{B\rtimes\Sigma} v_hq$. Let $S\subset\Lambda$ be a complete set of representatives for $\Lambda/\Sigma$.
The above observations imply that the formula $\pi (f)=\sum_{h\in S}f(h\Sigma) v_{h}q\;e_{B\rtimes\Sigma}\;qv_h^*$ defines a $*$-isomorphism $\pi:\ell^{\infty}(\Lambda/\Sigma)\rightarrow \mathcal D$. Moreover, we have $\pi(k\cdot f)=v_kq\pi(f)qv_k^*$, for every $k\in\Lambda$ and $f\in\ell^{\infty}(\Lambda/\Sigma)$.

 Now, we claim that $\varphi(v_kqTqv_k^*)=\varphi(T)$, for all $k\in\Lambda$ and $T\in\mathcal D$.
Since $v_kq\in\mathcal N_{q(A\rtimes\Gamma)q}(Aq)$, we can find mutually orthogonal projections $a_g\in Aq$, $g\in\Gamma$,  such that $v_kq=\sum_{g\in\Gamma}u_g a_g$ and $\sum_{g\in\Gamma}a_g=q$, where both series converge in $\|.\|_2$.  Note that $Aq=L^{\infty}(Y)$ commutes with $\mathcal D$, hence $a_g$ commutes with $\mathcal D$, for every $g\in\Gamma$. Moreover, $a_gqv_k^*=a_gu_g^*$, for every $g\in\Gamma$.
Also, the Cauchy-Schwarz inequality implies that $|\varphi(xV))|, |\varphi(Vx)|\leq \|x\|_2\|V\|$, for every $x\in pMp$ and $V\in p\langle M,e_{B\rtimes\Sigma}\rangle p$. By combining these facts with the fact that $\varphi$ is $L(\Gamma)$-central we obtain that
$$
\begin{array}{rcl}
\varphi(v_kqTqv_k^*)&=&\sum_{g\in\Gamma}\varphi(u_ga_gTqv_k^*) \\
&=& \sum_{g\in\Gamma}\varphi(u_gTa_gqv_k^*)\\
&{=}& \sum_{g\in\Gamma}\varphi(u_gTa_gu_g^*) \\
&=&\sum_{g\in\Gamma}\varphi(Ta_g)\\
&{=}&\varphi(T).
\end{array}
$$
 
It is now clear that  the positive linear functional $\varphi\circ\pi:\ell^{\infty}(\Lambda/\Sigma)\rightarrow\mathbb C$ is $\Lambda$-left invariant, which implies that $\Sigma$ is co-amenable in $\Lambda$. This finishes the proof in the case $(\nu\times c)(X)\geq 1$.

In general, let $r\geq 1$ such that $r(\nu\times c)(X)\geq 1$, and put $X_1=X\times\mathbb Z/r\mathbb Z$, $\Gamma_1=\Gamma\times\mathbb Z/r\mathbb Z$. Then $X_1\subset Y\times\mathbb Z/\ell\mathbb Z\times\mathbb Z/r\mathbb Z$ and if we consider a bijection $\mathbb Z/\ell\mathbb Z\times\mathbb Z/r\mathbb Z\equiv\mathbb Z/\ell r\mathbb Z$, we have that $\mathcal R(\Gamma_1\curvearrowright X_1)=\mathcal R(\Lambda\times\mathbb Z/\ell r\mathbb Z\curvearrowright Y\times\mathbb Z/\ell r\mathbb Z)|_{X_1}$. Denote $M_1=L^{\infty}(Y\times\mathbb Z/\ell r\mathbb Z)\rtimes(\Lambda\times\mathbb Z/\ell r\mathbb Z)$, $B_1=L^{\infty}(Y)\rtimes\mathbb M_{\ell r}(\mathbb C)$, and identify $M_1=B_1\rtimes\Lambda$. 
Since the inclusion $B_1\rtimes\Sigma\subset M_1$ is naturally identified to the inclusion $(B\rtimes\Sigma)\otimes\mathbb M_r(\mathbb C)\subset M\otimes\mathbb M_r(\mathbb C)$, we get that 
$L^{\infty}(\mathbb Z/r\mathbb Z)\rtimes\Gamma_1\equiv L(\Gamma)\otimes\mathbb M_r(\mathbb C)$ is amenable relative to $B_1\rtimes\Sigma$ inside $M_1$. Thus, $L(\Gamma_1)$ is amenable relative to $B_1\rtimes\Sigma$ inside $M_1$. Since $(\nu\times c)(X_1)\geq 1$, we can apply the above and derive that $\Sigma$ is co-amenable in $\Lambda$.
\hfill$\blacksquare$

\subsection{Relatively strongly solid groups}\label{SS: RSS}
In his breakthrough work  \cite{Oz03}, Ozawa proved that II$_1$ factors arising from non-elementary hyperbolic groups $\Gamma$ (e.g. $\Gamma=\mathbb F_n, 2\leq n\leq\infty$) are {\it solid}: if $P_1,P_2\subset L(\Gamma)$ are commuting von Neumann subalgebras, then either $P_1$ is not diffuse or $P_2$ is amenable. In the last ten years, this result has been generalized and strengthened in many ways. Remarkably, Ozawa and Popa proved that if $\Gamma=\mathbb F_n, 2\leq n\leq\infty$, then $L(\Gamma)$ {\it strongly solid}: the normalizer $\mathcal N_{L(\Gamma)}(P)''$ is amenable, for any diffuse amenable von Neumann subalgebra $P\subset L(\Gamma)$.  Chifan and Sinclair extended this to cover all non-elementary hyperbolic groups $\Gamma$ \cite{CS11}. 

Most recently, a breakthrough was made by Popa and Vaes who proved that non-abelian free groups and, more generally, non-elementary hyperbolic groups $\Gamma$ are relatively strong solid \cite{PV11,PV12}. Following \cite[Definition 2.7]{CIK13}, we say that a countable non-amenable group $\Gamma$ is {\it relatively strongly solid} and write $\Gamma\in\mathcal C_{\text{rss}}$ if for any trace preserving action $\Gamma\curvearrowright Q$ on a tracial von Neumann algebra $(Q,\tau)$ the following holds: if $M=Q\rtimes\Gamma$ and $P\subset pMp$ is any von Neumann subalgebra which is amenable relative to $Q$, then either $P\prec_{M}Q$ or the normalizer $\mathcal N_{pMp}(P)''$ is amenable relative to $Q$. Note that $\mathcal C_{\text{rss}}$ more generally contains all weakly amenable groups that either admit a proper $1$-cocycle into an orthogonal representation weakly contained in the left regular representation \cite[Theorem 1.6]{PV11}, or are bi-exact \cite[Theorem 1.4]{PV12}.
 
We will use repeatedly the following consequence of belonging to $\mathcal C_{\text{rss}}$ (see \cite[Lemma 5.2]{KV15}).

\begin{lemma}[\!\!{\cite{KV15}}]\label{KV15}
Let $\Gamma$ be a group in $\mathcal C_{\text{rss}}$, and $M=Q\rtimes\Gamma$, where $\Gamma\curvearrowright Q$ is a trace preserving action on a tracial von Neumann algebra. 
Let $P_1, P_2\subset M$ be commuting von Neumann subalgebras.

Then either $P_1\prec_{M}Q$ or $P_2$ is amenable relative to $Q$.
\end{lemma}

\section{From intertwining to measure equivalence}
The main goal of this section is to establish the following proposition, which provides the tool used to deduce the measure equivalence in part (2) of Theorem \ref{C}:

\begin{proposition}\label{P: twine2ME}
Let $\R$ be a countable pmp equivalence relation on $(X, \mu)$, and $Y, Z \subset X$ be positive measure subsets. Suppose that $\R|_Y = \R(\Gamma_1 \times \Gamma_2 \curvearrowright Y)$ and $\R|_Z \ge \R(\Lambda \curvearrowright Z)$ for free measure preserving actions of countable groups $\Gamma_1, \Gamma_2$, and $\Lambda$.
Assume that
\begin{enumerate}[label=(\roman*)]
\item $L^\infty(Y) \rtimes \Gamma_1 \prec_{L(\R)} L^\infty(Z) \rtimes \Lambda$, and 
\item $L^\infty(Z) \rtimes \Lambda \prec^s_{L(\R)} L^\infty(Y) \rtimes \Gamma_1$.
\end{enumerate}
Then $\Gamma_1$ and $\Lambda$ are measure equivalent. 
\end{proposition}

Throughout this section, all subsets of probability spaces that we consider are assumed measurable.

In order to prove Proposition \ref{P: twine2ME}, we first establish a series of lemmas in subsections \ref{SS: index}-\ref{SS: twineR}. The proof of Proposition \ref{P: twine2ME} is then given in Subsection \ref{SS: propproof}.

\subsection{Essentially finite index subequivalence relations}\label{SS: index}

Consider an inclusion of countable pmp equivalence relations $\T \le \mathcal{R}$ on $(X, \mu)$. 
Decompose $X = \bigsqcup_{N \in \{1, 2, \dots, \aleph_0\}} X_N$, where the $X_N$ are the $\R$-invariant sets defined by
\begin{align}\label{E: X_N}
X_N = \{x \in X \;|\; \text{$[x]_\R$ is the union of $N$ $\T$-classes}\} \quad\text{for}\quad N = 1, 2, \dots, \aleph_0.
\end{align}

If $\mu(X_\infty) = 0$, we say that the inclusion $\T \le \R$ has {\it essentially finite index}. If in fact there exists $k \ge 1$ such that $\mu(X_N) = 0$ for all $N > k$, the inclusion is said to have {\it bounded index}. 

We will use the following basic fact, whose proof we include for the sake of completeness.

\begin{lemma}\label{L: biintersect}
Let $\S, \T \le \R$ be inclusions of pmp countable equivalence relations and suppose that $\S \le \R$ has essentially finite (respectively, bounded) index. 

Then $\S \cap \T \le \T$ has essentially finite (respectively, bounded) index. 
\end{lemma}

{\it Proof.} Note that if $C_{\S}$ is an $\S$-class and $C_{\T}$ is a $\T$-class, then $C_{\S} \cap C_{\T}$ is either empty or equal to the $(\S \cap \T)$-class of any of its elements. Hence for $x \in X$, if $\mathcal{C}_x$ denotes the set of $\S$-classes in $[x]_{\R}$, we have that \begin{align*} [x]_{\T} = [x]_{\T} \cap [x]_{\R} = \bigsqcup_{C \in \mathcal{C}_x} ([x]_{\T} \cap C) \end{align*} is the union of at most $|\mathcal{C}_x|$ $(\S \cap \T)$-classes. If $\S \le \R$ is essentially finite (resp. bounded) index, then $|\mathcal{C}_x| < \infty$ (resp. there is $k \ge 1$ such that $|\mathcal{C}_x| < k$) for a.e. $x \in X$, and hence $\S \cap \T \le \T$ has essentially finite (resp. bounded) index. \hfill$\blacksquare$ \smallskip

The product structure $\Gamma_1 \times \Gamma_2$ assumed in Proposition \ref{P: twine2ME} will be exploited via the following lemma:
 
\begin{lemma}\label{L: prod2bi}
Let $\R = \R(\Gamma_1 \times \Gamma_2 \curvearrowright X)$ for a pmp action of the product of countable groups $\Gamma_1, \Gamma_2$ on $(X, \mu)$. Let $\T = \R(\Gamma_1 \curvearrowright X)$, $Y \subset X$ a positive measure subset, $\T_0 \le \T|_Y$ a subequivalence relation, and $\theta \in [[\R]]$ with $Y = \operatorname{dom}(\theta)$ such that $(\theta \times \theta)(\T_0) \le \T|_{\theta(Y)}$. Assume that  $\T_0 \le \T|_Y$ has essentially finite (respectively, bounded) index.
  
Then there is a sequence of $\T_0$-invariant positive measure $Y_n \subset Y$ with $Y = \bigcup_{n = 1}^\infty Y_n$ such that 
 $(\theta \times \theta)(\T_0|_{Y_n}) \le \T|_{\theta(Y_n)}$ has essentially finite (respectively, bounded) index for each $n \ge 1$. 
\end{lemma}
{\it Proof.}
Enumerate $\Gamma_2 = \{s_n\}_{n = 1}^\infty$ and let 
\begin{align*}
Y_{n} = \{x \in Y\;|\; \text{there exists}\; h_1(x) \in \Gamma_1 \text{ such that } \theta(x) = h_1(x)s_nx\}
\end{align*} 
Then $Y = \bigcup_{n = 1}^\infty Y_n$, and each $Y_n$ is $\T_0$-invariant, for if $x \in Y_{n}$ and $(x, x') \in \T_0$, then $(\theta(x), \theta(x')) \in \T$, so there is $k_1 \in \Gamma_1$ such that 
$\theta(x') = k_1\theta(x) = k_1h_1(x)s_nx$.

Now for any $x, x' \in Y_{n}$ such that $(\theta(x), \theta(x')) \in \T$, there is $k_1 \in \Gamma_1$ such that 
$\theta(x') = k_1\theta(x) = k_1h_1(x)s_nx$,
and on the other hand,
$\theta(x') = h_1(x')s_nx'$ and so we conclude that $x' = h_1(x')^{-1}k_1h_1(x)x$, giving $(x, x') \in \T$. 
Thus we have
\begin{align*}
\T_0|_{Y_n} \le (\theta^{-1} \times \theta^{-1})(\T|_{\theta(Y_n)}) \le \T|_{Y_n},
\end{align*}
and since $\T_0|_{Y_n} \le \T|_{Y_n}$ has bounded index, so too does $\T_0|_{Y_n} \le (\theta^{-1} \times \theta^{-1})(\T|_{\theta(Y_n)})$ and its image under $\theta$, as desired. 
\hfill$\blacksquare$

\subsection{Realizing subequivalence relations as restrictions}\label{SS: sub2action}
We recall in this section a useful construction appearing in \cite{IKT08}.
Consider as above an inclusion of countable pmp equivalence relations $\S \le \mathcal{R}$ on $(X, \mu)$ and the decomposition $X = \bigsqcup_{N \in \{1, 2, \dots, \aleph_0\}} X_N$ defined by \eqref{E: X_N}. 
For each $N$, let $\{C_n^{(N)}\}_{0 \le n < N}$ be a sequence of choice functions, i.e. a sequence of Borel functions $C^{(N)}_n: X_N \to X_N$ such that for each $x \in X_N$ the sequence $\{C_n^{(N)}(x)\}_{n = 0}^{N-1}$ contains exactly one element of each $\S$-class contained in $[x]_\R$. We take $C_0^{(N)} = \operatorname{Id}_X$.

Each $(x, y) \in \R|_{X_N}$ gives rise to a permutation $\pi_N(x, y) \in S_N$ defined by
\begin{align*}
m = \pi_N(x, y)(n) \iff (C_m^{(N)}(x), C_n^{(N)}(y)) \in \S,
\end{align*}
and the map $\pi_N: \R|_{X_N} \to S_N$ is called the {\it index cocycle} associated to these choice functions (see \cite{FSZ89}).  

Let 
\begin{align}
(\tilde X, \lambda) = \bigsqcup_{N \in \{1,2, \dots, \aleph_0\}} \left(X_N \times \{0, \dots, N-1\},\; \mu \otimes c\right)
\end{align}
where $c$ denotes the counting measure. 
In the case of an essentially finite index inclusion $\S \le \R$, we may instead endow the space $\tilde X$ with an $\tilde \R$-invariant probability measure $\tilde \mu$ by normalizing the counting measure:
\begin{align}\label{E: prob}
(\tilde X, \tilde \mu) = \bigsqcup_{N \in \{1, 2, \dots, \aleph_0\}} \left(X_N \times \{0, \dots, N-1\},\; \mu \otimes \frac{c}{N}\right)
\end{align}

We define a measurable equivalence relation $\tilde \R$ on $\tilde X$ by
\begin{align}
((x, n), (y, m)) \in \tilde \R \iff (x, y) \in \R \text{ and } n = \pi_N(x, y)(m), 
\end{align}
for $x, y \in X_N$, $n, m \in \{0, \dots, N-1\}$.
For $x, y \in X_N$ we have $((x, 0), (y, 0)) \in \tilde \R$ if and only if $(x, y) \in \R$ and $\pi_N(x, y)(0) = 0$ which occurs exactly when $(x, y) = (C_0^{(N)}(x), C_0^{(N)}(y)) \in \S$. Thus 
\begin{align}\label{E: restrict}
\S = \tilde \R|_{X \times \{0\}}
\end{align}

Now let $p: \tilde X \to X$ be the projection map $p(x, n) = x$. 
Any element $\phi \in [[\R]]$ gives rise to $\tilde \phi \in [[\tilde \R]]$ defined by 
\begin{align}
\tilde \phi:{} p^{-1}(\operatorname{dom} \phi) &\to p^{-1}(\operatorname{ran} \phi) \notag\\
(x, n) &\mapsto (\phi(x), \pi_N(\phi(x), x)(n)) \quad\text{for}\quad x \in X_N,\; n \in \{0, \dots, N-1\}, \label{E: Rrep}
\end{align}
such that $\tilde\phi\tilde\psi = \widetilde{\phi\psi}$ for $\phi, \psi \in [[\R]]$. 
In particular, if $\R = \R(\Gamma \curvearrowright X)$ is given by the free pmp action of a countable group $\Gamma$, then $\tilde \R$ is given by the free measure preserving action $\Gamma \curvearrowright \tilde X$ defined by
\begin{align}\label{E: upaction}
g \cdot (x, n) = (gx, \pi_N(gx, x)(n)) \quad\text{for}\quad x \in X_N,\; n \in \{0, \dots, N-1\}.
\end{align}

\subsection{A stable orbit equivalence-type characterization of measure equivalence}\label{SS: MEchar}

The main purpose of this subsection is to prove Lemma \ref{L: MEchar}, which allows one to deduce that countable groups $\Gamma$ and $\Lambda$ admit SOE free pmp actions (and hence are ME) from a seemingly weaker condition. 
We begin with the following general ergodic-theoretic lemma, whose proof we include for the sake of completeness.  

\begin{lemma}\label{L: saturate}
Let $\R$ be a countable pmp equivalence relation on $(X, \mu)$ and $E \subset X$ a positive measure subset. 

Then there exist a positive measure subset $E_0 \subset E$ and $\phi_0 = \operatorname{id}_{E_0}, \phi_1, \dots, \phi_k \in [[\R]]$, such that
\begin{enumerate}
\item $\operatorname{dom}(\phi_i) = E_0$ for $i = 0, \dots, k$,
\item $\operatorname{ran}(\phi_i) \cap \operatorname{ran}(\phi_j) = \emptyset$ for $i \ne j$, and
\item $Y = \bigsqcup_{i = 0}^k \operatorname{ran}(\phi_i)$ is $\R$-invariant. 
\end{enumerate}
\end{lemma}

{\it Proof.}
Let $Z$ be the set of $\R$-ergodic invariant probability measures on $X$, and let $\pi: X \to Z$ denote the ergodic decomposition of $\mu$ with respect to $\R$ (see \cite[Theorem 18.5]{KM04}). 
Thus, if we denote  $\nu=\pi_*\mu$, then $\mu=\int_Z m\; d\nu(m)$. Consider the natural embedding $L^2(Z)\ni f\mapsto f\circ\pi\in L^2(X)$ and denote by $e:L^2(X)\rightarrow L^2(Z)$ the orthogonal projection, noting that for $f \in L^2(X)$, we have that $e(f)$ is given by $e(f)(x) = \int_X f(y) d(\pi(x))(y)$ for almost every $x \in X$. 

Since $\mu(E) = \int_Z m(E) \;d\nu$ is positive, the set $Z_1=\{m\in Z|m(E)>0\}$ has positive measure. 
Since each $m \in Z_1$ is $\R$-ergodic, there is a positive measure subset $Z_0 \subset Z_1$ such that either $(X, m)$ is non-atomic for all $m \in Z_0$ or such that there is an integer $k \ge 0$ with $m$ supported on $k+1$ atoms each of measure $\frac{1}{k+1}$ for all $m \in Z_0$.
In any case, we can find an integer $k \ge 0$ and a measurable subset $E_0 \subset E \cap \pi^{-1}(Z_0)$ with $m(E_0) = \frac{1}{k+1}$ for all $m \in Z_0$. Moreover, we can then find measurable subsets $E_1,...,E_k\subset\pi^{-1}(Z_0)$ with $\pi^{-1}(Z_0)=E_0\cup E_1\cup...\cup E_k$ such that $E_i\cap E_j=\emptyset$ for all $0\leq i<j\leq k$, and 
$m(E_i)=\frac{1}{k+1}$ for all $0\leq i\leq k$ and $m\in Z_0$.

Since $\pi^{-1}(Z_0)$ is $\mathcal R$-invariant, in order to get the conclusion, it suffices to prove the following claim:

{\bf Claim.} Let $A,B\subset X$ be measurable sets satisfying $m(A)=m(B)$, for almost every $m\in Z$. Then there is $\theta\in[[\mathcal R]]$ such that dom$(\theta)=A$ and ran$(\theta)=B$. 

To this end, let $\{A_j\}_{j\in J}$ and $\{B_j\}_{j\in J}$ be maximal families of disjoint non-negligible measurable subsets of $A$ and $B$ such that for every $j\in J$ we can find $\theta_j\in [[\mathcal R]]$ with dom$(\theta_j)=A_j$ and ran$(\theta_j)=B_j$. Since $\sum_{j\in J}\mu(A_j)\leq\mu(A)\leq 1$, we deduce that $J$ is countable.
In particular, the sets $A'=\cup_{j\in J}A_j$, $B'=\cup_{j\in J}B_j$, $A''=A\setminus A'$, and $B''=B\setminus B'$ are measurable. 

Our goal is to show that $\mu(A'')=\mu(B'')=0$. Granting this, $\theta\in[[\mathcal R]]$ given by $\theta(x)=\theta_j(x)$ for all $x\in A_j$ and $j\in J$ satisfies $\theta(A')=B'$, and since $\mu(A\setminus A')=\mu(B\setminus B')=0$, the claim follows.

 Assume by contradiction that $\mu(A'')=\mu(B'')>0$.
For any $m\in Z$ and $j\in J$, since $m$ is $\mathcal R$-invariant and $B_j = \theta_j(A_j)$, we have $m(B_j)=m(A_j)$. Together with the assumption made on $A$ and $B$, this implies that $m(A'')=m(B'')$, for almost every $m\in Z$.

Let us show that there is $\rho\in [\mathcal R]$ such that $\mu(\rho(A'')\cap B'')>0$. Otherwise, we would have that $\int_{B''}{\bf 1}_{A''}\circ\rho \;\text{d}\mu=0$, for all $\rho\in[\mathcal R]$. Thus, if $\mathcal K\subset L^2(X,\mu)$ denotes the $\|.\|_2$-closure of the convex hull of $\{{\bf 1}_{A''}\circ\rho|\rho\in[\mathcal R]\}$, then $\int_{B''}f\;\text{d}\mu=0$, for every $f\in\mathcal K$. If $f\in\mathcal K$ denotes the element of minimal $\|.\|_2$, then $f$ is $\mathcal R$-invariant, hence $f=e(f)$. 
Moreover, since $e({\bf 1}_{A''}\circ\rho)=e({\bf 1}_{A''})$, for all $\rho\in[\mathcal R]$, we conclude that $f=e({\bf 1}_{A''})\geq 0$.
This and the condition $\int_{B''}f\;\text{d}\mu=0$ imply that $\pi(x)(A'')=f(x)=0$, for almost every $x\in B''$. Thus, $\pi(x)(B'')=0$, for almost every $x\in B''$, contradicting our assumption that $\mu(B'')>0$.

Finally, let $\tilde A=A''\cap\rho^{-1}(B''), \tilde B=\rho(A'')\cap B''$, and $\tilde\theta\in[[\mathcal R]]$ be the restriction of $\rho$ to $\tilde A$. Since $\mu(\tilde A)=\mu(\tilde B)>0$, $\tilde\theta(\tilde A)=\tilde B$, and $\tilde A\cap A'=\tilde B\cap B'=\emptyset$, this contradicts the maximality of the families $\{A_j\}_{j\in J}$ and $\{B_j\}_{j\in J}$, and finishes the proof of the claim.
\hfill$\blacksquare$ \smallskip

\begin{lemma}\label{L: fi2restrict}
Let $\R = \R(\Gamma \curvearrowright X)$ for a free pmp action of a countable group $\Gamma$ and let $E \subset X$ be a positive measure subset.

Then there exists a positive measure subset $E_0 \subset E$ with the following property: for any essentially finite index subequivalence relation $\T \le \R|_{E_0}$ there is a free pmp action $\Gamma~\curvearrowright~(\tilde X, \tilde \mu)$ such that $\T \cong \R(\Gamma \curvearrowright \tilde X)|_{\tilde E_0}$ for some measurable subset $\tilde E_0 \subset \tilde X$.
\end{lemma}
{\it Proof.}
Let $E_0 \subset E$, $Y \subset X$ and $\phi_0, \dots, \phi_k \in [[\R]]$ be as in the conclusion of Lemma \ref{L: saturate}. 
Let $\S = \R(\Gamma \curvearrowright X)|_Y$ and note that since $Y$ is $\Gamma$-invariant we have $\S = \R(\Gamma \curvearrowright Y)$ with $\Gamma$ acting freely. 
Define a subequivalence relation $\S_0 \le \S$ by $\S_0 = \bigsqcup_{j = 0}^k (\phi_j \times \phi_j)(\T)$.

Then for $x \in Y$, 
\begin{align*}
[x]_\S 
= \bigsqcup_{j = 0}^k [x]_\S \cap \phi_j(E_0) 
= \bigsqcup_{j = 0}^k \phi_j([x]_\S \cap E_0)
\end{align*}
and as $[x]_\S \cap E_0$ is the union of finitely many $\T$-classes, we see that $[x]_\S$ is the union of finitely many $\S_0$-classes. 
Thus, $\S_0 \le \S$ is an essentially finite index inclusion. 
Let $\Gamma \curvearrowright (\tilde X, \tilde \mu)$ be the free pmp action arising from this inclusion as in \eqref{E: prob} and \eqref{E: upaction}. 
Then by \eqref{E: restrict} we have $\S_0 \cong \R(\Gamma \curvearrowright \tilde X)|_{Y \times \{0\}}$ and so $\T = \S_0|_{E_0} \cong \R(\Gamma \curvearrowright \tilde X)|_{E_0 \times \{0\}}$ as desired. \hfill$\blacksquare$

\begin{lemma}\label{L: MEchar}
Let $\Gamma \curvearrowright (X, \mu)$ and $\Lambda \curvearrowright (Y, \nu)$ be free pmp actions of countable groups. Suppose there are positive measure subsets $E \subset X$, $F \subset Y$ and essentially finite index subequivalence relations $\T \le \R(\Gamma \curvearrowright X)|_E$ and $\S \le \R(\Lambda \curvearrowright Y)|_F$ with $\T \cong \S$.

Then $\Gamma$ and $\Lambda$ admit SOE free pmp actions (and hence are measure equivalent).  
\end{lemma}

{\it Proof.}
Applying Lemma \ref{L: fi2restrict} we find a free pmp action $\Gamma \curvearrowright (\tilde X, \tilde \mu)$ and positive measure subsets $E_0 \subset E$, $\tilde E_0 \subset \tilde X$ such that $\T|_{E_0} \cong \R(\Gamma \curvearrowright \tilde X)|_{\tilde E_0}$. Let $\Psi: \tilde E_0 \to E_0$ denote the measure space isomorphism implementing this identification. 

Since $\T \cong \S$, let $\theta: E \to F$ be a measure space isomorphism such that $(\theta \times \theta)(\T) = \S$.
Let $F_0 = \theta(E_0)$ and again apply Lemma \ref{L: fi2restrict} to find a positive measure subset $F_1 \subset F_0$ and a free pmp action $\Lambda \curvearrowright (\tilde Y, \tilde \nu)$ such that $\S|_{F_1} \cong \R(\Lambda \curvearrowright \tilde Y)|_{\tilde F_1}$ for some measurable $\tilde F_1 \subset \tilde Y$. 

Letting $E_1 = \theta^{-1}(F_1)$ and $\tilde E_1 = \Psi^{-1}(E_1)$ we see that 
\begin{align*}
\R(\Gamma \curvearrowright \tilde X)|_{\tilde E_1} \cong \T|_{E_1} \cong \S|_{F_1} \cong \R(\Lambda \curvearrowright \tilde Y)|_{\tilde F_1},
\end{align*} 
giving the desired stable orbit equivalence. \hfill$\blacksquare$

\subsection{Intertwining subequivalence relations}\label{SS: twineR}
We will need the techniques of \cite{Io11} which give the analogue of Popa's intertwining in the setting of countable pmp equivalence relations. 
Consider an inclusion of countable pmp equivalence relations $\mathcal{S} \le \mathcal{R}$ on $(X, \mu)$ such that each $\R$-class contains infinitely many $\S$-classes. 
For a positive measure subset $E \subset X$, the formula \eqref{E: Rrep} gives rise to a unitary representation $\rho: [\R|_{E}] \to \mathcal{U}(L^2(E \times \mathbb{Z}_{\ge 0}))$ defined by
\begin{align*}
[\rho(\theta)\xi](x, n) = \xi(\tilde \theta^{-1}(x, n)) \quad\text{for}\quad \xi \in L^2(E \times \mathbb{Z}_{\ge 0}).
\end{align*}
For $\xi \in L^2(E \times \mathbb{Z}_{\ge 0})$, denote $S(\xi) = \{x \in E\;|\;\sup_i |\xi(x, i)| \ne 0\}$. For further reference, we note that if $\xi$ is $\rho([\mathcal T])$-invariant, for some subequivalence relation $\mathcal T\leq\mathcal R|_{E}$, then $S(\xi)$ is $\mathcal T$-invariant.

Following \cite{IKT08} we define a function $\varphi_{\mathcal{S}}: [[\mathcal{R}]] \to [0, 1]$ by 
\begin{align*}
\varphi_{\mathcal{S}}(\theta) 
= \mu(\{x \in \operatorname{dom}(\theta) \;|\;(\theta(x), x) \in \mathcal{S}\}). 
\end{align*}

The following result established in \cite{Io11} shows the connection between $\varphi_{\S}$, Popa's intertwining, and intertwining of subequivalence relations: 

\begin{lemma}[\!\!{\cite[Lemmas 1.7 and 1.8]{Io11}}]\label{L: twine2R}
Let $E \subset X$ be a positive measure subset and $\T \le \R|_{E}$ a subequivalence relation. Then the following are equivalent:
\begin{enumerate}
\item $L(\T) \prec_{L(\R)} L(\S)$.
\item There is no sequence $\{\theta_n\}_{n = 1}^\infty \subset [\mathcal{T}]$ such that $\varphi_{\mathcal{S}}(\psi \theta_n \psi') \to 0$ for all $\psi, \psi' \in [\mathcal{R}]$.
\item  There is a non-zero $\rho([\T])$-invariant vector $\eta \in L^2(E \times \mathbb{Z}_{\ge 0})$.  
Moreover, in this case,  
 there is a subequivalence relation $\T_0 \le \T$ such that for any positive measure subset  $E_0 \subset S(\eta)$ there is a positive measure subset $Y \subset E_0$ and $\theta \in [[\R]]$, $\theta: Y \to Z$, satisfying
\begin{enumerate}
\item $\T_0|_{Y} \le \T|_{Y}$ has bounded index, and
\item $(\theta \times \theta)(\T_0|_{Y}) \le \S|_{Z}$. 
\end{enumerate}
\end{enumerate}

\end{lemma}

 In order to exploit strong intertwining $L(\T) \prec^s_{L(\R)} L(\S)$, we will use the following lemma:

\begin{lemma}\label{L: strongtwine}
Let $\T \le \mathcal{R}$ be a subequivalence relation such that for all $\T$-invariant subsets $E \subset X$ of positive measure there is no sequence $\{\theta_n\}_{n = 1}^\infty \subset [\mathcal{T}|_{E}]$ such that $\varphi_{\mathcal{S}}(\psi \theta_n \psi') \to 0$ for all $\psi, \psi' \in [\mathcal{R}]$.

Then there is a non-zero $\rho([\T])$-invariant vector $\eta \in L^2(X \times \mathbb{Z}_{\ge 0})$ such that $\mu(S(\eta)) = 1$.  
\end{lemma}
{\it Proof.}
Let $\mathcal F$ be the set of families $\{\eta_i\}_{i\in I}\subset L^2(X\times\mathbb Z_{\ge 0})$ of $\rho([\mathcal T])$-invariant vectors
 which satisfy $S(\eta_i)\cap S(\eta_j)=\emptyset$, for all $i\not= j$, and $\|\eta_i\|_2=\sqrt{\mu(S(\eta_i))}>0$, for all $i$. 
 By Zorn's lemma, we can find a family $\{\eta_i\}_{i\in I}\in\mathcal F$ that is maximal with respect to inclusion.   
 
 We claim that $\sum_{i\in I}\mu(S(\eta_i))=1$. Indeed, otherwise $E=X\setminus(\cup_{i\in I}S(\eta_i))$  would be a $\mathcal T$-invariant set of positive measure. By applying Lemma \ref{L: twine2R} $(2)\Rightarrow (3)$ we find a non-zero $\rho(\mathcal T)$-invariant vector $\xi\in L^2(X\times\mathbb Z_{\ge 0})$ with $S(\xi)\subset E$. But then the family $\{\eta_i\}_{i\in I}\cup\{\frac{\sqrt{\mu(S(\xi))}}{\|\xi\|_2}\xi\}$ also belongs to $\mathcal F$, which contradicts the maximality of $\{\eta_i\}_{i\in I}$, and thus proves the claim.
 
 It is now clear $\eta=\sum_{i\in I}\eta_i\in L^2(X\times\mathbb Z_{\ge 0})$ is a $\rho([\T])$-invariant unit vector with $\mu(S(\eta))=1$.
\hfill$\blacksquare$

We can now prove the intertwining lemma to be used in the proof of Proposition \ref{P: twine2ME}. A countable pmp equivalence relation $\T$ on $(Y,\nu)$ is called {\it aperiodic} 
if $[y]_{\T}$ is infinite for almost every $y\in Y$.
\begin{lemma}\label{L: twineR}
Let $\R$ be a countable pmp equivalence relation on $(X, \mu)$, $Y, Z \subset X$ subsets of positive measure, and $\T \le \R|_Y$, $\S \le \R|_Z$ subequivalence relations with $\T$ aperiodic. 

If $L(\T) \prec_{L(\R)} L(\S)$, then there is a subequivalence relation $\T_0 \le \T$, subsets of positive measure $Y_1 \subset Y$, $Z_1 \subset Z$, and $\theta \in [[\R]]$, $\theta: Y_1 \to Z_1$, such that
\begin{enumerate}
\item $\T_0|_{Y_1} \le \T|_{Y_1}$ has bounded index, and
\item $(\theta \times \theta)(\T_0|_{Y_1}) \le \S|_{Z_1}$. 
\end{enumerate}

 If we assume moreover that $L(\T) \prec_{L(\R)}^s L(\S)$, then for any positive measure $Y_0 \subset Y$, the subset $Y_1$ above can be taken with $Y_1 \subset Y_0$. 
\end{lemma}
{\it Proof.}
Let $\S' = \S \sqcup \{(x, x) \;|\; x \in X \setminus Z\}$ and note that $L(\T) \prec_{L(\R)} L(\S)$ implies $L(\T) \prec_{L(\R)} L(\S')$. 
Then by Lemma \ref{L: twine2R}, we can find $\T_0 \le \T$, positive measure subsets $Y_1 \subset Y$, $Z_1 \subset X$ and $\theta \in [[\R]]$, $\theta: Y_1 \to Z_1$, such that conclusions (1) and (2) hold. 
Since $\T$ is aperiodic, conclusion~(1) forces the $\T_0|_{Y_1}$-class of almost every $x \in Y_1$ to be infinite, and so conclusion~(2) forces $\mu(Z_1 \cap Z) = \mu(Z_1)$, and so we may indeed take $Z_1 \subset Z$. 

The moreover conclusion follows because Lemma \ref{L: strongtwine} allows us to apply the moreover assertion of Lemma \ref{L: twine2R} with $E_0$ a positive measure subset of $Y_0$.
\hfill$\blacksquare$

\subsection{Proof of Proposition \ref{P: twine2ME}}\label{SS: propproof} 
Let $\T = \R(\Gamma_1 \curvearrowright Y)$ and $\S = \R(\Lambda \curvearrowright Z)$. 
By assumption {\it(i)} and Lemma \ref{L: twineR}, there is a subequivalence relation $\T_0 \le \T$, positive measure subsets $Y_1 \subset Y$, $Z_1 \subset Z$, and $\theta \in [[\R]]$, $\theta: Y_1 \to Z_1$, such that $\T_0|_{Y_1} \le \T|_{Y_1}$ has bounded index and $(\theta \times \theta)(\T_0|_{Y_1}) \le \S|_{Z_1}$.
 
Similarly, by assumption {\it(ii)} and Lemma \ref{L: twineR}, there is a subequivalence relation $\S_0 \le \S$, positive measure subsets $Z_2 \subset Z$, $Y_2 \subset Y$, and $\phi \in [[\R]]$, $\phi: Z_2 \to Y_2$, such that $\S_0|_{Z_2} \le \S|_{Z_2}$ has bounded index and $(\phi \times \phi)(\S_0|_{Z_2}) \le \T|_{Y_2}$. Moreover, by Lemma \ref{L: twineR} we can take $Z_2 \subset Z_1$, since {\it(ii)} assumes strong intertwining.

Define $Y_2' \subset Y_1$ by $Y_2' = \theta^{-1}(Z_2)$ and let $\psi=\phi \circ \theta: Y_2' \to Y_2$.  
Since $\S_0|_{Z_2} \le \S|_{Z_2}$ has bounded index and $(\theta \times \theta)(\T_0|_{Y_2'}) \le \S|_{Z_2}$, Lemma \ref{L: biintersect} gives that 
\begin{align*}
\S_0|_{Z_2} \cap (\theta \times \theta)(\T_0|_{Y_2'}) \le (\theta \times \theta)(\T_0|_{Y_2'})
\end{align*} 
has bounded index. Letting
\begin{align}
\T_{00} = (\theta^{-1} \times \theta^{-1})(\S_0|_{Z_2} \cap (\theta \times \theta)(\T_0|_{Y_2'})) = (\theta^{-1} \times \theta^{-1})(\S_0|_{Z_2}) \cap \T_0|_{Y_2'}
\end{align}
we see that $\T_{00} \le \T_0|_{Y_2'}$ has bounded index and therefore so to does $\T_{00} \le \T|_{Y_2'}$.
Moreover, $(\theta \times \theta)(\T_{00}) \le \S_0|_{Z_2}$ and so 
\begin{align}
(\psi\times\psi)(\T_{00})\le\T|_{Y_2}. \label{E: inc}
\end{align}
As $Y_2', Y_2 \subset Y$, we may regard $\psi$ in $[[\R|_{Y}]] = [[\R(\Gamma_1 \times \Gamma_2 \curvearrowright Y)]]$ and apply Lemma \ref{L: prod2bi} to find positive measure $Y_3' \subset Y_2'$ such that for $Y_3 = \psi(Y_3')$, the inclusion \eqref{E: inc} has bounded index when restricted to $Y_3$, i.e. 
$
(\psi\times\psi)(\T_{00}|_{Y_3'})\le\T|_{Y_3}
$
has bounded index. 
Let $Z_3 = \theta(Y_3')$. Then because 
\begin{align*}
(\psi\times\psi)(\T_{00}|_{Y_3'}) \le (\phi \times \phi)(\S_0|_{Z_3}) \le \T|_{Y_3},
\end{align*}
we conclude that $(\phi \times \phi)(\S_0|_{Z_3}) \le \T|_{Y_3}$ has bounded index. 

Thus, $\S_0|_{Z_3}$ is a subequivalence relation of $\S|_{Z_3} = \R(\Lambda \curvearrowright Z)|_{Z_3}$ with bounded index whose isomorphic image $(\phi \times \phi)(\S_0|_{Z_3})$ has bounded index in $\T|_{Y_3} = \R(\Gamma_1 \curvearrowright Y)|_{Y_3}$. 
An application of Lemma \ref{L: MEchar} finishes the proof. 
\hfill$\blacksquare$

\bibliographystyle{\string~/Dropbox/Research/References/myamsalpha}
\bibliography{\string~/Dropbox/Research/References/References}

\section{Transfer of commutation from subalgebras to subgroups}

In this section we prove the following result which will be crucial in the proof of Theorem \ref{C}. This result is an immediate consequence of the ``ultrapower technique" developed in \cite{Io11}, being essentially contained in the proof of \cite[Theorem 3.1]{Io11} (see also \cite[Theorem 3.3]{CdSS15} and \cite[Lemma 5.6]{KV15}). Nevertheless, for completeness, we include a detailed proof.

\begin{theorem}[\!\!{\cite{Io11}}]\label{ultraproduct}
Let $M$ be a II$_1$ factor and $p\in M$ a projection such that $pMp=A\rtimes\Gamma$, where $\Gamma\curvearrowright A$ is a trace preserving action on a tracial von Neumann algebra.  Let $\Delta:M\rightarrow M\overline{\otimes}\text{L}(\Gamma)$ be a $*$-homomorphism which satisfies $\Delta(a)=a\otimes 1$ and $\Delta(u_g)=u_g\otimes u_g$, for all $a\in A$ and $g\in\Gamma$. Assume that $P\subset L(\Gamma)$ and $Q\subset qMq$ are von Neumann subalgebras such that $\Delta(Q)\prec_{M\overline{\otimes}L(\Gamma)}M\overline{\otimes}P$.
 
 Then there exists a decreasing sequence of subgroups $\Omega_k<\Gamma$ such that 
 
 \begin{enumerate}
 \item $Q\prec_{M}A\rtimes\Omega_k$, for all $k\geq 1$, and
 \item $P'\cap L(\Gamma)\prec_{L(\Gamma)}L(\cup_{k\geq 1}C_{\Gamma}(\Omega_k))$.
 \end{enumerate}

\end{theorem}

Throughout this section, we assume the setting of Theorem \ref{ultraproduct}.  Since $M$ is a II$_1$ factor, after replacing $Q$ with a unitary conjugate of one its corners, we may clearly assume that $q\leq p$.

In preparation for the proof of Theorem \ref{ultraproduct}, let us introduce some notation. We denote by $\mathcal G$ the family of all subgroups $\Sigma<\Gamma$ such that $Q\nprec_{M}A\rtimes\Sigma$. We may assume that $\mathcal G$ is non-empty. Indeed, if $\mathcal G=\emptyset$, then $Q\prec_{M}A$, and thus the conclusion holds with $\Omega_k=\{e\}$, for every $k\geq 1$.

We say that a set $S\subset\Gamma$ is {\it small relative to $\mathcal G$} if $S\subset\cup_{i=1}^mb_i\Sigma_ic_i$, for some $b_i,c_i\in\Gamma$ and $\Sigma_i\in\mathcal G$.  We denote by $I$ the family of subsets of $\Gamma$ that are small relative to $\mathcal G$. We order $I$ by inclusion and consider a cofinal ultrafilter $\mathcal V$ on $I$. Thus, $\{S'\in I|S'\supset S\}$ belongs to $\mathcal V$, for every $S\in I$.

\begin{lemma}\label{transfer}
We can find a finite set $F\subset L(\Gamma)$ and $\delta>0$ such that the following holds: for any $S\subset\Gamma$ which is small relative to $\mathcal G$, there exists $g\in\Gamma\setminus S$ such that $\sum_{\alpha,\beta\in F}\|E_P(\alpha u_{g}\beta)\|_2^2\geq\delta$.
\end{lemma}

{\it Proof.} The proof uses the ``transfer of rigidity" principle from \cite{PV09} (see the proof of \cite[Lemma 3.2]{PV09}).
Since $\Delta(Q)\prec_{M\overline{\otimes}L(\Gamma)}M\overline{\otimes}P$, we can find $F\subset (L(\Gamma))_1$ finite and $\kappa>0$ such that \begin{equation}\label{1}\sum_{\alpha,\beta\in F}\|E_{M\overline{\otimes}P}((1\otimes\alpha)\Delta(u)(1\otimes\beta))\|_2^2\geq\kappa,\;\;\;\text{for every $u\in\mathcal U(Q)$}.\end{equation}
Put $\delta=\frac{\kappa}{2\|q\|_2^2}$. Let $S\subset\Gamma$ be small relative to $\mathcal G$. Thus, $S\subset\cup_{i=1}^mb_i\Sigma_ic_i$, for some $b_i,c_i\in\Gamma$ and $\Sigma_i\in\mathcal G$. For $g\in\Gamma$, we denote $\varphi(g)=\sum_{\alpha,\beta\in F}\|E_P(\alpha u_{g}\beta)\|_2^2$. Since $F\subset (L(\Gamma))_1$, we have that $\varphi(g)\leq |F|^2$, for every $g\in\Gamma$.
Our goal is to show the existence of $g\in\Gamma\setminus S$ such that $\varphi(g)\geq\delta$.

Since $Q\nprec_{M}A\rtimes\Sigma_i$, for every $i\in\{1,...,m\}$, by Remark \ref{corner2} we can find $u\in\mathcal U(Q)$ such that \begin{equation}\label{2}\|E_{A\rtimes\Sigma_i}(u_{b_i}^*uu_{c_i}^*)\|_2^2\leq\frac{\kappa}{2m|F|^2},\;\;\;\text{for every $1\leq i\leq m$}.\end{equation}
Since $u\in Q\subset qMq\subset q(A\rtimes\Gamma)q$, we can write $u=\sum_{g\in\Gamma}a_gu_g$, where $a_g\in A$. By using \ref{2} we get that
\begin{equation}\label{3}
\sum_{g\in S}\varphi(g)\|a_g\|_2^2\leq|F|^2\sum_{g\in S}\|a_g\|_2^2\leq |F|^2\sum_{i=1}^m\|E_{A\rtimes\Sigma_i}(u_{b_i}^*uu_{c_i}^*)\|_2^2\leq\frac{\kappa}{2}.
\end{equation}
On the other hand, since $\Delta(u)=\sum_{g\in\Gamma}a_gu_g\otimes u_g$, equation \ref{1}  rewrites as $\sum_{g\in\Gamma}\varphi(g)\|a_g\|_2^2\geq\kappa$. In combination with \ref{3} this gives that $\sum_{g\in\Gamma\setminus S}\varphi(g)\|a_g\|_2^2\geq\frac{\kappa}{2}$. Since $\sum_{g\in\Gamma\setminus S}\|a_g\|_2^2\leq\|u\|_2^2=\|q\|_2^2$, it follows that we can find $g\in\Gamma\setminus S$ such that $\varphi(g)\geq\delta$, as claimed.
\hfill$\blacksquare$

{\bf Proof of Theorem \ref{ultraproduct}.} Denote $N=L(\Gamma)$. By Lemma \ref{transfer}, for every $S\in I$ we can find $g_S\in\Gamma\setminus S$ such that $\sum_{\alpha,\beta\in F}\|E_P(\alpha u_{g_S}\beta)\|_2^2\geq\delta$. We put $g=(g_S)_{S\in I}\in\Gamma^{\mathcal V}$ and consider the canonical inclusions $\Gamma\subset\Gamma^{\mathcal V}\subset\mathcal U(N^{\mathcal V})$. 
We define $\Sigma=\Gamma\cap g\Gamma g^{-1}$ and claim that $P'\cap N\prec_{N}L(\Sigma)$.

Assume by contradiction that this is false. By Theorem \ref{corner}, we can find a sequence $u_n\in\mathcal U(P'\cap N)$ such that $\|E_{L(\Sigma)}(xu_ny)\|_2\rightarrow 0$, for every $x,y\in N$. We denote by $\mathcal K\subset L^2(N^{\mathcal V})$ the closed linear span of $Nu_gN$, and by $e$ the orthogonal projection from $L^2(N^{\mathcal V})$ onto $\mathcal K$. 

Let us show that $\langle u_n\xi u_n^*,\eta\rangle\rightarrow 0$, for every $\xi,\eta\in\mathcal K$. To prove this, it suffices to show that $\langle u_nxu_gy u_n^*,x'u_gy'\rangle\rightarrow 0$, for every $x,y\in N$. But this is clear since 
$\|E_{L(\Sigma)}(x'^*u_nx)\|_2\rightarrow 0$ and
\begin{align*}\langle u_nxu_gy u_n^*,x'u_gy'\rangle=&\tau(u_g^*(x'^*u_nx)u_g(yu_n^*y'^*))=\tau(E_N(u_g^*(x'^*u_nx)u_g)yu_n^*y'^*)\\=&\tau(E_N(u_g^*E_{L(\Sigma)}(x'^*u_nx)u_g)yu_n^*y'^*).\end{align*}

Next, since $\sum_{\alpha,\beta\in F}\|E_{P^{\mathcal V}}(\alpha u_g\beta)\|_2^2=\lim\limits_{S\rightarrow\mathcal V}(\sum_{\alpha,\beta\in F}\|E_P(\alpha u_{g_S}\beta)\|_2^2)\geq\delta$, we can find $\alpha,\beta\in F$ such that $E_{P^{\mathcal V}}(\alpha u_g\beta)\not=0$. 
Thus, $\|E_{P^{\mathcal V}}(\alpha u_g\beta)-\alpha u_g\beta\|_2<\|\alpha u_g\beta\|_2$.
Since $\alpha u_g\beta\in\mathcal K$, we get that $\|e(E_{P^{\mathcal V}}(\alpha u_g\beta))-\alpha u_g\beta\|_2<\|\alpha u_g\beta\|_2$.
This implies that
 $\xi=e(E_{P^{\mathcal V}}(\alpha u_g\beta))\in\mathcal K$ is non-zero.
 On the other hand, as $e$ is $N$-$N$-bimodular and $u_n\in P'\cap N$, we get that $u_n\xi u_n^*=e(u_nE_{P^{\mathcal V}}(\alpha u_g\beta)u_n^*)=\xi$
 and therefore $\langle u_n\xi u_n^*,\xi\rangle=\|\xi\|_2^2>0$, for every $n$. This contradicts the previous paragraph and thus proves the claim.
 
Now, enumerate $\Sigma=\{\sigma_j\}_{j\geq 1}$. If $\sigma\in\Gamma$, then $\sigma$ belongs to $\Sigma$ if and only if $\sigma$ commutes with $\{g_Sg_{S'}^{-1}|S,S'\in T\}$, for some $T\in\mathcal V$. In particular, for every $j\geq 1$, we can find $T_j\in\mathcal V$ such that $\sigma_j$ commutes with $\{g_Sg_{S'}^{-1}|S,S'\in T_j\}$. For $k\geq 1$, define $W_k=\cap_{j=1}^kT_j$ and $\Omega_k=\langle g_Sg_{S'}^{-1}|S,S'\in W_k\rangle$. 
 Then $W_k\in\mathcal V$ and $\Omega_k\supset\Omega_{k+1}$. Since $\sigma_1,...,\sigma_k\in C_{\Gamma}(\Omega_k)$, we deduce that $\Sigma=\cup_{k\geq 1}C_{\Gamma}(\Omega_k).$
 
To finish the proof, it suffices to show that if $W\in\mathcal V$, then $\Omega=\langle g_Sg_{S'}^{-1}|S,S'\in W\rangle$ does not belong to $\mathcal G$.
 Indeed, this implies that $\Omega_k\not\in\mathcal G$ and hence that $Q\prec_{M}A\rtimes\Omega_k$, for every $k\geq 1$.
 Assume by contradiction that $\Omega\in\mathcal G$. Fix $S'\in W$. Then $g_S\in\Omega g_{S'}$, for every $S\in W$. Since $\Omega\in\mathcal G$, the set $\Omega g_{S'}\subset\Gamma$ is small relative to $\mathcal G$. Since $\mathcal V$ is cofinal and $W\in\mathcal V$, we get that $W\cap\{S''\in I|S''\supset\Omega g_{S'}\}$ belongs to $\mathcal V$, and hence is non-empty. Let $S''\in W$ such that $S''\supset\Omega g_{S'}$. But then we get that $g_S\in S''$, for every $S\in W$. Taking $S=S''$, this contradicts the fact that $g_{S''}\in\Gamma\setminus S''$.
\hfill$\blacksquare$

\section{Groups measure equivalent to products of hyperbolic groups \\ and  tensor decompositions}\label{S: tech}
The proof of Theorem \ref{C} is divided between this and the next section. 
Before stating the main result of this section, we need to introduce some notation. 

\begin{notation}\label{ME} Let $\Gamma$ be an icc group which is measure equivalent to a product $\Lambda=\Lambda_1\times...\times\Lambda_n$ of $n\geq 1$  groups belonging to $\mathcal C_{\text{rss}}$. By \cite[Lemma 3.2]{Fu99}, $\Gamma$ and $\Lambda$ admit stably orbit equivalent free ergodic pmp actions. 
We may thus find a free ergodic pmp action $\Lambda\curvearrowright (Y,\nu)$ and $\ell\geq 1$, such that the following holds: consider the product action $\Lambda\times\mathbb Z/\ell\mathbb Z\curvearrowright (Y\times\mathbb Z/\ell\mathbb Z,\nu\times c)$, where $\mathbb Z/\ell\mathbb Z$ acts on itself by addition and $c$ denotes the counting measure on $\mathbb Z/\ell\mathbb Z$.  Then there is a non-negligible measurable set $X\subset Y\times\mathbb Z/\ell\mathbb Z$ and a free ergodic measure preserving action $\Gamma\curvearrowright X$ such that $$\mathcal R(\Gamma\curvearrowright X)=\mathcal R(\Lambda\times\mathbb Z/\ell\mathbb Z\curvearrowright Y\times\mathbb Z/\ell\mathbb Z)|_{X}.$$

We put $A=L^{\infty}(X), M=L^{\infty}(Y\times\mathbb Z/\ell\mathbb Z)\rtimes(\Lambda\times\mathbb Z/\ell\mathbb Z)$, $p=1_X$, and note that
$A\rtimes\Gamma=pMp$. We identify $L^{\infty}(\mathbb Z/\ell\mathbb Z)\rtimes\mathbb Z/\ell\mathbb Z=\mathbb M_{\ell}(\mathbb C)$, and use this identification to write $M=B\rtimes\Lambda$, where $B=L^{\infty}(Y)\otimes\mathbb M_{\ell}(\mathbb C)$ and $\Lambda$ acts trivially on $\mathbb M_{\ell}(\mathbb C)$. 
We let $\{u_g\}_{g\in\Gamma}\subset A\rtimes\Gamma$ and $\{v_h\}_{h\in\Lambda}\subset B\rtimes\Lambda$ denote the canonical unitaries 
implementing the actions of $\Gamma$ and $\Lambda$ on $A$ and $B$, respectively.

For a set $T\subset\{1,2,...,n\}$, we denote $\Lambda_T=\underset{j\in T}{\times}\Lambda_j$ and let $\widehat T = \{1,2,...,n\} \setminus T$.

We define a $*$-homomorphism $\Delta:M\rightarrow M\overline{\otimes}L(\Gamma)$ as follows \cite{PV09}.
Let $k\geq\tau(p)^{-1}$ be an integer, where $\tau$ denotes the trace of $M$.
Let $\widetilde{\Delta}:\mathbb M_k(pMp)\rightarrow \mathbb M_k(pMp)\overline{\otimes}L(\Gamma)$ be the  $*$-homomorphism given by $\widetilde{\Delta}(a)=a\otimes 1$ and 
$\widetilde{\Delta}(u_g)=u_g\otimes u_g$, for all $a\in \mathbb M_k(A)$ and $g\in\Gamma$. 
 Let $q\in\mathbb M_k(A)$ be a projection satisfying $(\text{Tr}\otimes\tau)(q)=1$ and $e_{11}\otimes p\leq q$, where $\text{Tr}$ is the non-normalized trace of $\mathbb M_k$. We fix an identification  $\alpha:M\rightarrow q\mathbb M_k(pMp)q$ which satisfies $\alpha(x)=e_{11}\otimes x$, for all $x\in pMp$.
 Since $\widetilde{\Delta}(q)=q\otimes 1$, we have that $\widetilde{\Delta}(q\mathbb M_k(pMp)q)\subset q\mathbb M_k(pMp)q\overline{\otimes}L(\Gamma)$. 
 
 Finally, we put $\Delta=(\alpha^{-1}\otimes\text{id})\circ\widetilde{\Delta}\circ\alpha:M\rightarrow M\overline{\otimes}L(\Gamma)$. Then one checks that
$$\text{$\Delta(a)=a\otimes 1$ and $\Delta(u_g)=u_g\otimes u_g$, for every $a\in A$ and $g\in\Gamma$}.$$ 
  
  For further reference, we also record two facts. Firstly, if $\Gamma$ is icc, then $\Delta(M)'\cap M\overline{\otimes}L(\Gamma)=\mathbb C$.
  Indeed, if $\Gamma$ is icc, it is easy to see that $\widetilde{\Delta}(\mathbb M_k(pMp))'\cap \mathbb M_k(pMp)\overline{\otimes}L(\Gamma)=\mathbb C$, which  gives the fact.
The second fact goes back to \cite[Proposition 7.2.4]{IPV10}. In the more general context needed below, it is due to \cite[Proposition 2.4]{KV15}).
\end{notation}

\begin{lemma}[\!\!{\cite{KV15}}]\label{L: comult}
If $N \subset M$ has no amenable direct summand, then $\Delta(N)p'$ is non-amenable relative to $M \overline{\otimes} \mathbb{C}$ inside $M \overline{\otimes} L(\Gamma)$ for any non-zero projection $p' \in \Delta(N)' \cap M \overline{\otimes} L(\Gamma)$. 
\end{lemma}

The following is the main result of this section:

\begin{theorem}\label{tensor}
Assume that $L(\Gamma)=P_1\overline{\otimes}P_2$, where $P_1,P_2\subset L(\Gamma)$ are II$_1$ factors.

Then there are subgroups $\Sigma_1,\Sigma_2<\Gamma$ and a partition $S_1\sqcup S_2=\{1,...,n\}$ such that for all $i\in\{1,2\}$,
\begin{enumerate}
\item $P_i\prec^s_{L(\Gamma)}L(\Sigma_i)$, $L(\Sigma_i)\prec^s_{L(\Gamma)}P_i$, 
\item $A\rtimes\Sigma_i\prec^s_{M}B\rtimes\Lambda_{S_i}$, $B\rtimes\Lambda_{S_i}\prec^s_{M} A\rtimes\Sigma_i$, and
\item $\Sigma_i$ is measure equivalent to $\Lambda_{S_i}$.
\end{enumerate}
\end{theorem}
 
  The rest of the section is devoted to the proof of Theorem \ref{tensor}.
We assume throughout the notation from \ref{ME} and that $L(\Gamma)=P_1\overline{\otimes}P_2$.

\subsection{Outline of proof of Theorem \ref{tensor}} The proof of Theorem \ref{tensor} is divided between five steps, which we now briefly outline in order to facilitate reading.

\begin{labeling}{steps}
\item[\bf{Step 1.}]  {\it There is a partition $T_1\sqcup T_2=\{1,...,n\}$ such that $P_i\prec^s_M B\rtimes\Lambda_{T_i}$, for all $i\in\{1,2\}$.}
This conclusion will be obtained in Proposition \ref{step1} by using that $\Lambda_j\in\mathcal C_{\text{rss}}$, for all $1\leq j\leq n$.

\vskip 0.1in
\item[\bf{Step 2.}] {\it There is a  partition  $S_1\sqcup S_2=\{1,...,n\}$ such that $\Delta(B\rtimes\Lambda_{S_i})\prec_{M\overline{\otimes}L(\Gamma)} M\overline{\otimes}P_i$,  for all $i\in\{1,2\}$. }
This conclusion will be obtained in Proposition \ref{step2} by using that $\Lambda_j\in\mathcal C_{\text{rss}}$, for all $j$, and the embeddings $\varphi_i:P_i\rightarrow\mathbb M_{m_i}(B\rtimes\Lambda_{T_i})$ (for some $m_i\geq 1$) provided by {\bf Step 1}.

\vskip 0.1in
\item[\bf{Step 3.}] {\it There is a decreasing sequence of subgroups $\Omega_k<\Gamma$ such that
$B\rtimes\Lambda_{S_1}\prec_M A\rtimes\Omega_k$, for all $k\geq 1$, and 
 $P_2\prec_{L(\Gamma)} L(\cup_{k\geq 1}C_{\Gamma}(\Omega_k)).$} This is an immediate consequence of {\bf Step 2} and Theorem \ref{ultraproduct}; see Lemma \ref{ultra}.

\vskip 0.1in
\item[\bf{Step 4.}] {\it There is a subgroup $\Sigma_1<\Gamma$ such that  $B\rtimes\Lambda_{S_1}\prec^s_{M} A\rtimes\Sigma_1$,  $A\rtimes\Sigma_1\prec^s_{M}B\rtimes\Lambda_{S_1}$, $P_1\prec^s_{L(\Gamma)}L(\Sigma_1)$, and
 $L(\Sigma_1)\prec^s_{L(\Gamma)}P_1$.}
Specifically, Lemma \ref{Sigma_i} will show that $\Sigma_1=\Omega_k$ works, for $k$ large. A key part is showing that $L(\Omega_k)\prec_{L(\Gamma)}P_1$, for large $k$; see
Lemma \ref{Omega_k}.
This uses again that $\Lambda_j\in\mathcal C_{\text{rss}}$ for all $j$ and the embeddings $\varphi_i:P_i\rightarrow\mathbb M_{m_i}(B\rtimes\Lambda_{T_i})$ for $i\in\{1,2\}$. {\it Similarly, there is a subgroup $\Sigma_2<\Gamma$ with analogous properties.}
\vskip 0.1in
\item[\bf{Step 5.}]  {\it $\Sigma_i$ is measure equivalent to $\Lambda_{S_i}$, for every $i\in\{1,2\}$.} This will follow readily by combining the result of {\bf Step 4} with Proposition \ref{P: twine2ME}. 

\end{labeling}

\begin{remark} Since {\bf Steps 1-3} suffice in order to deduce Corollary \ref{B}, we include its proof right after {\bf Step 3}.
\end{remark}
\subsection{Step 1}

\begin{proposition}\label{step1} There is a partition $T_1\sqcup T_2=\{1,...,n\}$ such that $P_i\prec^s_M B\rtimes\Lambda_{T_i}$, for all $i\in\{1,2\}$. Moreover, 
if $P_i$ is amenable relative to $B\rtimes\Lambda_T$, for some $i\in\{1,2\}$ and $T\subset\{1,..,n\}$, then $T\supset T_i$.
\end{proposition}

{\it Proof.}
For $t\in\{1,...,n\}$, denote by $\hat t$ the set $\{1,...,n\}\setminus\{t\}.$
For $i\in\{1,2\}$, let $T_i\subset \{1,...,n\}$ be a minimal set with respect to inclusion such that $P_i$ is amenable relative to $B\rtimes\Lambda_{T_i}$. 

We claim that $P_2\prec^s_{M} B\rtimes\Lambda_{\{1,...,n\}\setminus T_1}.$ This is immediate if $T_1 = \emptyset$.\footnote{In fact, since each $P_i$ is type {\rm II}$_1$ and $B$ is type {\rm I}, after the proposition is proved, the conclusion that $P_i\prec^s_M B\rtimes\Lambda_{T_i}$ for all $i\in\{1,2\}$ will imply that $T_1$ and $T_2$ are nonempty.} Otherwise consider any $t \in T_1$. Since $\Lambda_{t}\in \mathcal {C}_{rss}$, Lemma \ref{KV15} implies that $P_1$ is amenable relative to $B\rtimes \Lambda_{\hat t}$ or $P_2\prec_M B\rtimes\Lambda_{\hat t}.$ Using Lemma \ref{PV11}(1) and the minimality of $T_1$, it follows that $P_2\prec_M B\rtimes\Lambda_{\hat t}.$ Since $\Gamma$ is icc and the action $\Gamma \curvearrowright X$ is ergodic, we have that $(\mathcal{N}_{pMp}(P_2))'\cap pMp\subset L(\Gamma)'\cap pMp=\mathbb C p$. Lemma \ref{facts1}(3) implies that $P_2\prec^s_{M} B\rtimes\Lambda_{\hat t}$. Since this holds for all $t\in T_1$,  Lemma \ref{PV11}(2) implies that $P_2\prec^s_{M} B\rtimes\Lambda_{\{1,...,n\}\setminus T_1}$ as claimed. Using the minimality of $T_2$ and Lemma \ref{facts2}(3), we get that $T_1\cap T_2=\emptyset$.
In a similar way we obtain that $P_1\prec^s_{M} B\rtimes\Lambda_{\{1,...,n\}\setminus T_2}$. 

The remaining part of the proof is to prove that $T_1\cup T_2=\{1,...,n\}$. 
We claim that $L(\Gamma)$ is not amenable relative to $B\rtimes \Lambda_{T}$ inside $M$, for any proper set $T\subsetneq \{1,...,n\}.$
Otherwise,  Lemma \ref{coamenable} would imply that $\Lambda_T<\Lambda$ is co-amenable, for some $T\subsetneq \{1,...,n\}.$ This would further give that $\Lambda_{\{1,...,n\}\setminus T}$ is amenable, which contradicts the fact that $\Lambda_j$ is non-amenable, for every $1\leq j\leq n$.

Next, fixing any $i\in\{1,2\}$, we claim that $P_i\prec^s_M B\rtimes \Lambda_{T_i}$. This is immediate if $T_i = \{1,...,n\}$; otherwise consider any $t\notin T_i$.
Then $P_i$ is amenable relative to $B\rtimes \Lambda_{\hat t}$ and since $\Lambda_t\in \mathcal C_{rss}$, we must have either $P_i\prec_M B\rtimes \Lambda_{\hat t}$ or $\mathcal N_{pMp}(P_i)''$ amenable relative to $B\rtimes \Lambda_{\hat t}$. Since $L(\Gamma)\subset \mathcal N_{pMp}(P_i)''$, the previous paragraph implies that $P_i\prec_M B\rtimes \Lambda_{\hat t}$, for all $t\notin T_i.$ As above, we get that $P_i\prec^s_M B\rtimes \Lambda_{\hat t}$, for all $t\notin T_i.$ Lemma \ref{PV11}(2) implies now that $P_i\prec^s_M B\rtimes \Lambda_{T_i}$, as claimed.

Thus, in particular $P_i\prec^s_M B\rtimes \Lambda_{T_1\cup T_2}$, for all $i\in \{1,2\}$. Applying \cite[Lemma 2.3]{BV12} implies that $L(\Gamma)\prec_M B\rtimes\Lambda_{T_1\cup T_2}$. As above, Lemma \ref{facts1}(3) implies that $L(\Gamma)\prec^s_M B\rtimes\Lambda_{T_1\cup T_2}$. Applying \cite[Lemma 2.3]{BV12} once again gives that $A\rtimes\Gamma\prec_M B\rtimes \Lambda_{T_1\cup T_2}.$ If there exists $t\in \{1,...,n\}\setminus( T_1\cup T_2)$, then we would get that $L(\Lambda_t) \prec_M B\rtimes \Lambda_{\hat t},$ which contradicts that $\Lambda_t$ is infinite.
Thus, $T_1\cup T_2=\{1,...,n\}$. The moreover assertion follows from the minimality of $T_1$ and $T_2$ using again Lemma \ref{PV11}(1).
\hfill$\blacksquare$

\subsection{Step 2}
Towards the second step of the proof of Theorem \ref{tensor}, we now prove that the each intertwining $P_i\prec^s_M B\rtimes\Lambda_{T_i}$ from Proposition \ref{step1} allows us to deduce that $P_i$ itself has a weaker form of relative solidity present in $B\rtimes\Lambda_{T_i}$. More precisely:

\begin{lemma}\label{L: *k}
Let $P = P_i$ and $k = |T_i|$ for some $i \in \{1, 2\}$. 
Then for any tracial von Neumann algebra $M_0$, any projection $q \in \tilde M = M_0 \overline{\otimes} P$, and any commuting subalgebras $Q_0, \dots, Q_k \subset q\tilde Mq$ we have either
\begin{enumerate}
\item $Q_0 \prec_{\tilde M}^s M_0$, or
\item $Q_jq'$ is amenable relative to $M_0$ inside $\tilde M$, for some $j \in \{1, \dots, k\}$ and some non-zero projection $q' \in Q_j' \cap q\tilde Mq$. 
\end{enumerate}
\end{lemma}
\noindent{\it Proof.} Assume that $Q_jq'$ is not amenable relative to $M_0$ inside $\tilde M$ for any $j \in \{1, \dots, k\}$ and non-zero projection $q' \in Q_j' \cap q\tilde Mq$. We first note that in order to prove the lemma, it suffices to show the conclusion $Q_0\prec_{\tilde M} M_0$.
Indeed, if this is known, then for any $z \in \mathcal{N}_{q\tilde Mq}(Q_0)' \cap q\tilde Mq \subset (\bigcup_{j=0}^kQ_j)' \cap q\tilde Mq$, applying the result to the commuting subalgebras $\{Q_jz\}_{j=0}^k\subset z\tilde Mz$ (noting that $Q_jq'$ is not amenable relative to $M_0$, for all $j\in\{1, \dots, k\}$ and any non-zero projection $q'\in (Q_jz)'\cap z\tilde Mz$), we conclude that $Q_0z \prec_{\tilde M} M_0$ and so by Lemma \ref{facts1}(2), $Q_0 \prec_{\tilde M}^s M_0$ as desired.

For an integer $m \ge 1$, let $e_{11}\in\mathbb M_{m}(\mathbb C)$ denote the matrix unit corresponding the $(1,1)$ entry and view $M$ as a non-unital subalgebra of $\mathbb M_{m}(M)$ via the embedding $x\mapsto x\otimes e_{11}$.
By Proposition \ref{step1} we have that $P \prec^s_M B \rtimes \Lambda_T$ for some $T \subset \{1, \dots, n\}$ with $|T| = k$. 
Hence we have for some $m \geq 1$ a not necessarily unital $*$-homomorphism $\varphi:P \rightarrow \mathbb M_{m}(B\rtimes\Lambda_T)$ and a non-zero partial isometry $v\in \mathbb M_{m,1}(M)p$ such that $\varphi(x)v=vx$, for every $x\in P$.
We define $e = \varphi(p)$, $\mathcal B =\mathbb M_{m}(B)$, and $\mathcal M = \mathbb M_{m}(B\rtimes\Lambda_T) \subset \mathbb{M}_m(M)$ and write canonically $\mathcal M = \mathcal B\rtimes\Lambda_T$. 
Moreover, we may assume that $E_{\mathcal M}(vv^*)\geq ce$, for some $c>0$. 

Replacing $\varphi$ by ${\rm id} \otimes \varphi$ we extend to $\varphi: M_0 \overline{\otimes} P \to M_0 \overline{\otimes} \mathbb M_m(M)$. Note that $\varphi(M_0 \overline{\otimes} P) \subset M_0 \overline{\otimes} \mathcal{M}$ and that $\varphi(x)v=vx$, for every $x\in M_0 \overline{\otimes} P$.
Let $f = \varphi(q)$ and $Q = (\bigcup_{j = 0}^k \varphi(Q_j))'' \subset f(M_0 \overline{\otimes} \mathcal{M})f$. 

\vskip 0.05in
\noindent {\bf Claim 1.} To prove that $Q_0 \prec_{\tilde M} M_0$, it is enough to show that $\varphi(Q_0) \prec_{M_0 \overline{\otimes} \mathcal{M}}^s M_0 \overline{\otimes} \mathcal{B}$.

\vskip 0.05in
\noindent {\it Proof of Claim 1.} Assume by contradiction that $\varphi(Q_0) \prec_{M_0 \overline{\otimes} \mathcal{M}}^s M_0 \overline{\otimes} \mathcal{B}$ and $Q_0 \nprec_{\tilde M} M_0$. Since $Q_0\prec_{M_0\overline{\otimes}\mathbb M_m(M)}\varphi(Q_0)$, Lemma \ref{facts1}(1) implies that $Q_0\prec_{M_0\overline{\otimes}\mathbb M_m(M)}M_0\overline{\otimes}\mathcal B$. From this we get that $Q_0\prec_{M_0\overline{\otimes}pMp}M_0\overline{\otimes}A$. 
 On the other hand, since $Q_0\nprec_{\tilde M}M_0$, by Theorem \ref{corner} we can find a sequence $u_n\in\mathcal U(Q_0)$ satisfying $\|E_{M_0}(xu_ny)\|_2\rightarrow 0$, for all $x,y\in\tilde M$. Let us show that $\|E_{M_0\overline{\otimes}A}(xu_ny)\|_2\rightarrow 0$, for all $x,y\in M_0\overline{\otimes}pMp$. This assertion will give a contradiction, and thus prove the claim.

To prove the assertion, recalling that $pMp=A\rtimes\Gamma$, it suffices to treat the case $x=1$ and $y\in L(\Gamma)$.
But then since $u_n\in Q_0$ and $Q_0\subset\tilde M\subset M_0\overline{\otimes} L(\Gamma)$ we get that $u_ny\in M_0\overline{\otimes}L(\Gamma)$ and thus $E_{M_0\overline{\otimes}A}(u_ny)=E_{M_0}(u_ny)=E_{M_0}(u_nE_{\tilde M}(y)).$
As $\|E_{M_0}(u_nE_{\tilde M}(y))\|_2\rightarrow 0$,  the claim is proven.
\hfill$\square$

\vskip 0.05in
\noindent {\bf Claim 2.} $\varphi(Q_j)q'$ is not amenable relative to $M_0 \overline{\otimes} \mathcal{B}$ inside $M_0 \overline{\otimes} \mathcal{M}$ for any $j \in \{1, \dots, k\}$ and any non-zero projection $q' \in  \varphi(Q_j)' \cap f(M_0 \overline{\otimes} \mathcal{M})f$.

\vskip 0.05in
\noindent {\it Proof of Claim 2.} 
Suppose the claim is false. Since $\mathcal{B}$ is amenable,  by \cite[Proposition 2.4(3)]{OP07}, we would conclude that there is $j \in \{1, \dots, k\}$ such that $\varphi(Q_j)q'$ is amenable relative to $M_0$ inside $M_0 \overline{\otimes} \mathcal M$ for some non-zero projection $q' \in\varphi(Q_j)' \cap f(M_0 \overline{\otimes} \mathcal M)f$. Thus, by Lemma \ref{facts2}(2), there is a projection $z\in\mathcal Z(\varphi(Q_j)'\cap f(M_0\overline{\otimes}\mathbb M_m(M))f)$ such that $q'\leq z$ and $\varphi(Q_j)z$ is amenable relative to $M_0$ inside $M_0\overline{\otimes}\mathbb M_m(M).$
Since  $E_{M_0\overline{\otimes}\mathcal M}(vv^*)\geq ce$, we get that
$v^*q'v\not=0$. Hence we deduce that $z'=v^*zv\in Q_j'\cap q(M_0\overline{\otimes}M)q$ is a non-zero projection such that $Q_jz'$ is amenable relative to $M_0$ inside $M_0\overline{\otimes}\mathbb M_m(M)$, and hence inside $M_0\overline{\otimes}pMp$.

Thus, we can find a $Q_jz'$-central positive linear functional $\psi:z'\langle M_0\overline{\otimes}pMp,e_{M_0}\rangle z'\rightarrow\mathbb C$ such that $\psi_{|z'(M_0\overline{\otimes}pMp)z'}=\tau$. The formula $\Psi(T)=\psi(z'Tz')$ defines a $Q_j$-central positive linear functional $\Psi:\langle M_0\overline{\otimes} pMp,e_{M_0}\rangle \rightarrow\mathbb C$ such that $\Psi(x)=\tau(xz')$, for any $x\in M_0\overline{\otimes}pMp$. 

Note that $L^2(pMp)\cong L^2(P)\otimes\ell^2$, as left $P$-modules. Thus, we can find a unitary operator $U:L^2(pMp)\rightarrow L^2(P)\otimes\ell^2$ such that $U(x\xi)=xU(\xi)$, for any $x\in P$ and $\xi\in L^2(pMp)$.
Let $V=\text{id}_{L^2(M_0)}\otimes U:L^2(M_0\overline{\otimes}pMp)\rightarrow L^2(M_0\overline{\otimes}P)\otimes\ell^2$ and $\theta:\mathbb B(L^2(M_0\overline{\otimes}P))\rightarrow\mathbb B(L^2(M_0\overline{\otimes}pMp))$ be the $*$-homomorphism given by $\theta(T)=V^*(T\otimes\text{id}_{\ell^2})V.$ Then $\theta(\langle M_0\overline{\otimes}P,e_{M_0}\rangle)\subset\langle M_0\overline{\otimes}pMp,e_{M_0}\rangle$ and $\theta(x)=x$, for every $x\in M_0\overline{\otimes}P$. Thus, if $\tilde\Psi:\langle M_0\overline{\otimes}P,e_{M_0}\rangle\rightarrow\mathbb C$ is given by $\tilde{\Psi}(T)=\Psi(\theta(T))$, then $\tilde{\Psi}$ is $Q_j$-central and satisfies $\tilde{\Psi}(x)=\tau(xz')$, for every $x\in M_0\overline{\otimes}P$. If we let $z''\in Q_j'\cap q(M_0\overline{\otimes}P)q$ be the support projection of $E_{M_0\overline{\otimes}P}(z')$, then \cite[Theorem 2.1]{OP07} implies that $Q_jz''$ is amenable relative to $M_0$ inside $\tilde M=M_0\overline{\otimes}P$, which is a contradiction. 
\hfill$\square$

For $j \in \{1, \dots, k\}$ and $S \subset T$, let $q_{j, S}$ be the maximal projection in $\mathcal{Z}(Q' \cap f(M_0 \overline{\otimes} \mathcal{M})f)$ such that $\varphi(Q_j)q_{j, S}$ is amenable relative to $M_0 \overline{\otimes} (\mathcal{B} \rtimes \Lambda_{S})$. 
Noting that $S' \subset S$ implies $q_{j, S'} \le q_{j, S}$, set 
\begin{align}\label{E: q'}
z_{j, S} = q_{j, S} - \bigvee_{S' \subsetneq S} q_{j, S'},
\end{align}
so that $z_{j, S}z_{j, S'} = 0$ whenever $S \ne S'$ by Lemma \ref{PV11}(1). 
Since $q_{j, T} = f$ it follows that if we let $\mathcal{F}_j = \{S \subset T | z_{j, S} \ne 0\}$, then $\sum_{S \in \mathcal{F}_j} z_{j, S} = f$ with the summands being mutually orthogonal.

\vskip 0.05in
\noindent {\bf Claim 3.} If $j \ne j'$ and $S \in \mathcal{F}_j$, $S' \in \mathcal{F}_{j'}$ with $z_{\scriptscriptstyle j, S}z_{\scriptscriptstyle j'\!\!, S'} \ne 0$, then $S \cap S' = \emptyset$. 

\vskip 0.05in
\noindent {\it Proof of Claim 3.} 
For any $\ell \in S$ and any nonzero projection $z \le z_{j, S}$, $z \in \mathcal{Z}(Q' \cap f(M_0 \overline{\otimes} \mathcal{M})f)$, we must have $\varphi(Q_j)z$ non-amenable relative to $M_0 \overline{\otimes} (\mathcal{B}\rtimes\Lambda_{T \setminus \{\ell\}})$. Otherwise, using Lemma \ref{PV11}(1) would give $\varphi(Q_j)z$ is amenable relative to $M_0 \overline{\otimes} (\mathcal{B} \rtimes \Lambda_{S \setminus \{\ell\}})$ implying $z \le q_{j, S \setminus \{\ell\}} \le 1-z_{j, S}$ (this last inequality coming from equation \eqref{E: q'}). Thus, decomposing $M_0 \overline{\otimes} \mathcal{M} = (M_0 \overline{\otimes} (\mathcal{B} \rtimes \Lambda_{T \setminus \{\ell\}})) \rtimes \Lambda_\ell$ and using that $\Lambda_{\ell}\in\mathcal C_{\text{rss}}$ and Lemma \ref{KV15} we conclude that 
\[\varphi(Q_{j'})z \prec_{M_0 \overline{\otimes} \mathcal{M}} M_0 \overline{\otimes} (\mathcal{B} \rtimes  \Lambda_{T \setminus \{\ell\}}).\]
Since 
\[
\mathcal{N}_{z_{j, S}(M_0 \overline{\otimes} \mathcal{M})z_{j, S}}(\varphi(Q_{j'})z_{j, S})' \cap z_{j, S}(M_0 \overline{\otimes} \mathcal{M})z_{j, S} \subset \mathcal{Z}((Qz_{j, S})' \cap z_{j, S}(M_0 \overline{\otimes} \mathcal{M})z_{j, S}),
\]
it follows by Lemma \ref{facts1}(2) that 
$\varphi(Q_{j'})z_{j, S} \prec^s_{M_0 \overline{\otimes} \mathcal{M}} M_0 \overline{\otimes} (\mathcal{B} \rtimes  \Lambda_{T \setminus \{\ell\}})$. 
Applying Lemma \ref{PV11}(2) to intersect over $\ell \in S$, we find that 
$\varphi(Q_{j'})z_{j, S} \prec^s_{M_0 \overline{\otimes} \mathcal{M}} M_0 \overline{\otimes} (\mathcal{B} \rtimes  \Lambda_{T \setminus S})$. 
Lemma \ref{facts2}(3) then implies that 
$\varphi(Q_{j'})z_{j, S}$ is amenable relative to $M_0 \overline{\otimes} (\mathcal{B} \rtimes  \Lambda_{T \setminus S})$. Hence $z_{j, S} \le q_{j'\!\!, \,T \setminus S}$, and so
\begin{align*}
0 < z_{\scriptscriptstyle j'\!\!, S'}z_{\scriptscriptstyle j, S} \le z_{\scriptscriptstyle j'\!\!, S'}q_{\scriptscriptstyle j'\!\!, T \setminus S} \le z_{\scriptscriptstyle j'\!\!, S'}q_{\scriptscriptstyle j'\!\!, S' \cap (T \setminus S)} 
\end{align*}
which forces $S' \cap (T \setminus S) = S'$ (that is, $S \cap S' = \emptyset$), since otherwise $q_{j'\!\!, \,S' \cap (T \setminus S)} \le 1 - z_{j'\!\!, S'}$ by equation \eqref{E: q'}. 
\hfill$\square$

\vskip 0.05in
\noindent {\bf Claim 4.} For each $\ell \in T$ we have $\bigvee \{z_{j, S} | \ell \in S, 1 \le j \le k, S \in \mathcal{F}_j\} = f$. 

\vskip 0.05in
\noindent {\it Proof of Claim 4.} To prove the claim it suffices to show that for any non-zero projection $q' \in \mathcal{Z}(Q' \cap f(M_0 \overline{\otimes} \mathcal{M})f)$, we have $\bigcup \{S \in \mathcal{F}_{j} | 1 \le j \le k, z_{\scriptscriptstyle j, S}q' \ne 0\} = T$. Indeed, assuming this condition, let $\ell \in T$ and put $f'=\bigvee \{z_{j, S} | \ell \in S, 1 \le j \le k, S \in \mathcal{F}_j\}$. Then $q'=f-f'$ satisfies 
$z_{j,S}q'=0$, for every $1\leq j\leq k$ and $S\in\mathcal{F}_j$ such that $\ell\in S$.
The assumed condition forces $q'=0$ and hence $f'=f$.

For each $1 \le j \le k$, using the fact that $\sum_{S \in \mathcal{F}_{j}} z_{j, S} = f$, pick (recursively) some $S_{j} \in \mathcal{F}_{j}$ such that $z_{j, S_{j}}q' \ne 0$ and $z_{j, S_{j}}z_{j'\!\!, S_{j'}} \ne 0$ for all $j' \le j$. 
Then using Claim 3 we have 
\begin{align*}
|\bigcup \{S \in \mathcal{F}_{j} | 1 \le j \le k, z_{j, S}q' \ne 0\}| \ge \sum_{j = 1}^k |S_j|.
\end{align*}
By Claim 2 we have $|S| > 0$ for all $S \in \mathcal{F}_j$, $j \in \{1, \dots, k\}$, so each of the $k = |T|$ terms in the above sum is positive. Thus $|\bigcup \{S \in \mathcal{F}_{j} | 1 \le j \le k, z_{j, S}q' \ne 0\}| = |T|$ and the claim follows.\footnote{This type of reasoning also implies that $|S| = 1$ for any $S \in \mathcal{F}_j, j \in \{1, \dots, k\}$, but we will not need this.} 
\hfill$\square$

\vskip 0.05in
\noindent {\bf Claim 5.} $\varphi(Q_0) \prec^s_{M_0 \overline{\otimes} \mathcal{M}} M_0 \overline{\otimes} (\mathcal{B} \rtimes  \Lambda_{T \setminus \{\ell\}})$ for each $\ell \in T$.

\vskip 0.05in
\noindent {\it Proof of Claim 5.} 
Fix $\ell \in T$. By Lemma \ref{facts1}(2) it is enough to show that 
\[\varphi(Q_0)z \prec_{M_0 \overline{\otimes} \mathcal{M}} M_0 \overline{\otimes} (\mathcal{B} \rtimes  \Lambda_{T \setminus \{\ell\}})\]
for any $z \in \mathcal{N}_{f(M_0\overline{\otimes}\mathcal M)f}(\varphi(Q_0))' \cap f(M_0 \overline{\otimes} \mathcal{M})f \subset \mathcal{Z}(Q' \cap f(M_0 \overline{\otimes} \mathcal{M})f)$. 
Fix any such $z$ and note that by Claim 4 we can find $j \in \{1, \dots, k\}$ and $S \in \mathcal{F}_j$ such that $\ell \in S$ and $zz_{j, S} \ne 0$. 
It follows that  $\varphi(Q_j)z$ is not amenable relative to $M_0 \overline{\otimes} (\mathcal{B} \rtimes \Lambda_{T \setminus \{\ell\}})$, otherwise Lemma \ref{PV11}(1) would give $\varphi(Q_j)zz_{j, S}$ amenable relative to $M_0 \overline{\otimes} (\mathcal{B} \rtimes \Lambda_{S \setminus \{\ell\}})$ implying $zz_{j, S} \le q_{j, S \setminus \{\ell\}} \le 1-z_{j, S}$ (this last inequality coming from equation \eqref{E: q'}). Decomposing $M_0 \overline{\otimes} \mathcal{M} = (M_0 \overline{\otimes} (\mathcal{B} \rtimes \Lambda_{T \setminus \{\ell\}})) \rtimes \Lambda_\ell$ and using that $\Lambda_\ell\in\mathcal C_{\text{rss}}$ and Lemma \ref{KV15} we conclude that $\varphi(Q_0)z \prec_{M_0 \overline{\otimes} \mathcal{M}} M_0 \overline{\otimes} (\mathcal{B} \rtimes  \Lambda_{T \setminus \{\ell\}})$, as desired. 
\hfill$\square$

Note that the subalgebras $\{M_0 \overline{\otimes} (\mathcal{B} \rtimes  \Lambda_{T \setminus \{\ell\}})\}_{\ell\in T}$ pairwise form commuting squares, are each regular in $M_0 \overline{\otimes} \mathcal{M}$, and have $\bigcap_{\ell \in T} M_0 \overline{\otimes} (\mathcal{B} \rtimes  \Lambda_{T \setminus \{\ell\}}) = M_0 \overline{\otimes} \mathcal{B}$. Hence Claim 5 together with Lemma \ref{PV11}(2) implies that $\varphi(Q_0) \prec^s_{M_0 \overline{\otimes} \mathcal{M}} M_0 \overline{\otimes} \mathcal{B}$. By Claim 1, this concludes the proof of the lemma. \hfill$\blacksquare$

\begin{proposition}\label{step2} There is a partition  $S_1\sqcup S_2=\{1,...,n\}$ such that $\Delta(B\rtimes\Lambda_{S_i})\prec_{M\overline{\otimes}L(\Gamma)} M\overline{\otimes}P_i$,  for all $i\in\{1,2\}$. \end{proposition}
\noindent{\it Proof.}  
Set $\tilde M = M \overline{\otimes}L(\Gamma) = M\overline{\otimes}P_1\overline{\otimes}P_2$, for $T \subset \{1, \dots, n\}$ let $Q_T = \Delta(L(\Lambda_T))$, and define $Q = (\bigcup_{j = 1}^n Q_j)'' = \Delta(L(\Lambda))$. For $i \in \{1, 2\}$, let $\hat i$ denote the element in $\{1, 2\} \setminus\{i\}$.  

\vskip 0.05in
\noindent{\bf Claim 1.} There are $i \in \{1, 2\}$, $S_i \subset \{1, \dots, n\}$ with $|S_i| = |T_i|$, and a non-zero projection $q \in \mathcal{Z}(Q' \cap \tilde M)$ such that $Q_jq'$ is not amenable relative to $M \overline{\otimes} P_{\widehat i}$ for all $j \in S_i$ and any non-zero projection $q' \in \mathcal{Z}((Qq)' \cap q\tilde Mq)$. 

\vskip 0.05in
\noindent{\it Proof of Claim 1.} For $j \in \{1, \dots, n\}$, $i \in \{1, 2\}$, let $q_{j, i}$ be the maximal projection in $\mathcal{Z}(Q' \cap \tilde M)$ such that $Q_jq_{j, i}$ is amenable relative to $M \overline{\otimes} P_{\widehat i}$ inside $\tilde M$. Then $Q_jq'$ is non-amenable relative to $M \overline{\otimes} P_{\widehat i}$ for all projections $q' \in \mathcal{Z}(Q' \cap \tilde M)$ with $q' \le 1 - q_{j, i}$, so it suffices to find $S_i \subset \{1, \dots, n\}$ with $|S_i| \ge |T_i|$ and $\bigwedge_{j \in S_i} (1 - q_{j, i}) \ne 0$.   
 Note that for each $j$ we have $Q_jq_{j, 1}q_{j, 2} = \Delta(L(\Lambda_j))q_{j, 1}q_{j, 2}$ amenable relative to $M$ by Lemma \ref{PV11}(1) and hence Lemma \ref{L: comult} forces $q_{j, 1}q_{j, 2} = 0$. 
 
 Let $S_1 \subset \{1, \dots, n\}$ be a maximal subset satisfying $q_1 = \bigwedge_{j \in S_1} (1 - q_{j, 1}) \ne 0$. If $|S_1| \ge |T_1|$ the claim holds with $i = 1$ and we are done. Otherwise, $S_2 = \widehat S_1$ will have $|S_2| \ge |T_2|$ and by the maximality of $S_1$, for any $j \in S_2$ we have $q_1 \le q_{j, 1} \le 1-q_{j, 2}$ and hence $\bigwedge_{j \in S_2} (1 - q_{j, 2}) \ge q_1 \ne 0$ so that the claim holds with $i = 2$. 
\hfill$\square$

\vskip 0.05in
For ease of notation, we assume without loss of generality that Claim 1 holds for $i = 1$. Set $S_2 = \hat S_1$. 

\vskip 0.05in
\noindent{\bf Claim 2.} $\Delta(L(\Lambda_{S_i})) = Q_{S_i} \prec_{\tilde M} M \overline{\otimes} P_i$ for all $i \in \{1, 2\}$. 

\vskip 0.05in
\noindent{\it Proof of Claim 2.}
We apply Lemma \ref{L: *k} with $M_0 = M \overline{\otimes} P_{2}$ to the commuting subalgebras $Q_{S_2}q, \{Q_jq\}_{j \in S_1} \subset q\tilde M q$. Alternative (2) of Lemma \ref{L: *k} cannot hold, for if there were $j \in S_1$ and a non-zero projection $q' \in (Q_jq)' \cap q\tilde Mq$ with $Q_jq'$ amenable relative to $M \overline{\otimes} P_{2}$,  Lemma \ref{facts2}(2) would give a projection $q'' \in \mathcal{N}_{q\tilde Mq}(Q_jq)' \cap q\tilde M q \subset \mathcal{Z}((Qq)' \cap q\tilde M q)$ with $q' \le q''$ (so $q'' \ne 0$) and $Q_jq''$ amenable relative to $M \overline{\otimes} P_{2}$, contradicting Claim 1. Thus Lemma \ref{L: *k} gives that $Q_{S_2}q \prec^s_{\tilde M} M \overline{\otimes} P_{2}$. This implies that $Q_{S_2} \prec_{\tilde M} M \overline{\otimes} P_{2}$, and that $Q_{S_2}q$ is amenable relative to $M \overline{\otimes} P_{2}$ by Lemma \ref{facts2}(3).

Hence for all $j \in S_2$ we have $Q_jq$ amenable relative to $M \overline{\otimes} P_{2}$. It follows that $Q_jq'$ is not amenable relative to $M \overline{\otimes} P_{1}$ for any $j \in S_2$ and non-zero projection $q' \in \mathcal{Z}((Qq)' \cap q\tilde Mq)$. Otherwise, Lemma \ref{PV11}(1) would give $Q_jq' = \Delta(L(\Lambda_j))$ amenable relative to $M$, contradicting Lemma \ref{L: comult}. We then apply Lemma \ref{L: *k} with $M_0 = M \overline{\otimes} P_{1}$ to the commuting subalgebras $Q_{S_1}q, \{Q_jq\}_{j \in S_2} \subset q\tilde M q$, and as before we conclude that $Q_{S_1}q \prec^s_{\tilde M} M \overline{\otimes} P_{1}$ and hence $Q_{S_1} \prec_{\tilde M} M \overline{\otimes} P_{1}$, establishing Claim 2. \hfill$\square$

\vskip 0.05in
We now finish the proof of the proposition. For any $i \in \{1, 2\}$, since $\mathcal{U}(\Delta(B \rtimes \Lambda_{S_i}))$ is generated by $\{\Delta(bu) : b \in \mathcal{U}(B), u \in \mathcal{U}(L(\Lambda_{S_i}))\}$ if we did not have $\Delta(B \rtimes \Lambda_{S_i}) \prec_{\tilde M} M \overline{\otimes} P_i$ there would be sequences $\{b_n\} \subset \mathcal{U}(B), \{u_n\} \subset \mathcal{U}(L(\Lambda_{S_i}))$ such that $\|E_{M \overline{\otimes} P_i}(x\Delta(b_nu_n)y)\|_2~\to~0$ for all $x, y \in \tilde M$. But then for any $x, y \in P_{\widehat i}$, using the fact that $\Delta(B) \subset M\overline{\otimes} P_i$ we would have 
\begin{align*}
\|E_{M \overline{\otimes} P_i}(x\Delta(u_n)y)\|_2
= \|\Delta(b_n)E_{M \overline{\otimes} P_i}(x\Delta(u_n)y)\|_2
= \|E_{M \overline{\otimes} P_i}(x\Delta(b_nu_n)y)\|_2 \to 0.
\end{align*}
Since $\tilde M = M \overline{\otimes} P_i \overline{\otimes} P_{\widehat i}$ it would further follow that $\|E_{M \overline{\otimes} P_i}(x\Delta(u_n)y)\|_2 \to 0$ for all $x, y \in \tilde M$, which would contradict Claim 2.  Hence we must have $\Delta(B \rtimes \Lambda_{S_i}) \prec_{\tilde M} M \overline{\otimes} P_i$ as desired. 
\hfill$\blacksquare$

\subsection{Step 3} Next, by combining {\bf Step 2} and Theorem \ref{ultraproduct}, we obtain:

 \begin{lemma}\label{ultra} We can find a decreasing sequence of subgroups $\Omega_k<\Gamma$ such that
\begin{itemize}
\item$B\rtimes\Lambda_{S_1}\prec_M A\rtimes\Omega_k$, for all $k\geq 1$, and 
\item $P_2\prec_{L(\Gamma)} L(\cup_{k\geq 1}C_{\Gamma}(\Omega_k)).$
\end{itemize}
 \end{lemma}
 
 {\it Proof.} 
 By Proposition \ref{step2} we have that $\Delta(B\rtimes\Lambda_{S_1})\prec_{M\overline{\otimes}L(\Gamma)}M\overline{\otimes}P_1$. Since $P_2\subset P_1'\cap L(\Gamma)$, the conclusion follows from
 Theorem \ref{ultraproduct}.
  \hfill$\blacksquare$

  \subsection{Proof of Corollary \ref{B}} Let $\Gamma=\text{PSL}_2(R)$, where either $R=\mathcal O_d$, for a square-free integer $d\geq 2$, or $R=\mathbb Z[S^{-1}]$, for a non-empty set of primes $S$. 
  Then the centralizer $C_{\Gamma}(g)$ of any non-trivial element $g\in\Gamma\setminus\{e\}$ is solvable, hence amenable. This follows from the following fact which can be derived by using for instance the Jordan normal form of matrices: if $A\in\text{SL}_2(\mathbb R)\setminus\{\pm I\}$, then the group $\{B\in \text{SL}_2(\mathbb R)|AB=\pm BA\}$ is solvable.
    In particular, we deduce that $\Gamma$ is icc and does not contain two commuting non-amenable subgroups.
  
  Assume by contradiction that $L(\Gamma)$ is not prime and write $L(\Gamma)=P_1\overline{\otimes}P_2$. Since $\Gamma$ is non-amenable, we may assume without loss of generality that $P_2$ is non-amenable.
  Since $\Gamma\in\mathscr L$ by Remark \ref{L}, $\Gamma$ is measure equivalent to a product $\Lambda=\Lambda_1\times...\times\Lambda_n$ of $n\geq 1$ non-elementary hyperbolic groups (where $n=2$, if $R=\mathcal O_d$, and $n=|S|+1$, if $R=\mathbb Z[S^{-1}]$). Since non-elementary hyperbolic groups are in class $\mathcal C_{\text{rss}}$ by \cite{PV12}, we are in the setting of \ref{ME}. 
  Thus, we may find a decreasing sequence of subgroups $\Omega_k<\Gamma$ satisfying Lemma \ref{ultra}. Since $\Lambda_i$ is non-amenable, for every $1\leq i\leq n$, and $P_2$ is non-amenable, it follows that for large enough $k$ we have that both $\Omega_k$ and $C_{\Gamma}(\Omega_k)$ are non-amenable. This contradicts the previous paragraph.
  \hfill$\blacksquare$

 \subsection{Step 4} This step is divided between two lemmas. We start with the following:
 
\begin{lemma}\label{Omega_k}  Let $\Omega_k$ be the decreasing sequence of subgroups of $\Gamma$ provided by Lemma \ref{ultra}.

 Then for any large enough $k\geq 1$ we have that $L(\Omega_k)\prec_{L(\Gamma)}P_1$.
 \end{lemma}

{\it Proof.} 
Let $i\in\{1,2\}$. By Proposition \ref{step1}, $P_i\prec_{M}B\rtimes\Lambda_{T_i}$. We can thus find  a not necessarily unital $*$-homomorphism $\varphi_i:P_i\rightarrow\mathbb M_{m_i}(M)$ and a non-zero partial isometry $v_i\in M_{m_i,1}(M)p$ such that $\varphi_i(x)v_i=v_ix$, for every $x\in P_i$, 
and $\varphi(P_i)\subset\mathcal M_i$, where
$\mathcal M_i=\mathbb M_{m_i}(B\rtimes\Lambda_{T_i})$, for some $m_i\geq 1$. Here, we view $P_i\subset M$ as non-unital subalgebras of $\mathbb M_{m_i}(M)$ via the embedding $x\mapsto x\otimes e_{11}$, where $e_{11}\in\mathbb M_{m_i}(\mathbb C)$ is the matrix unit corresponding the $(1,1)$ entry.
Moreover, we may assume that $E_{\mathcal M_i}(v_iv_i^*)\geq c_i\varphi_i(1)$, for some $c_i>0$.
We define $\mathcal B_i=\mathbb M_{m_i}(B)$
and write canonically $\mathcal M_i=\mathcal B_i\rtimes\Lambda_{T_i}$. 

We claim that  $\varphi_i(P_i)p'$ is not amenable relative to $\mathcal B_i\rtimes\Lambda_{T_i\setminus\{j\}}$ inside $\mathbb M_{m_i}(M)$, for any $j\in T_i$ and any non-zero projection $p'\in \varphi_i(P_i)'\cap \varphi_i(1)\mathbb M_{m_i}(M)\varphi_i(1)$ with $p'\leq v_iv_i^*$.  Otherwise, it would follow that  $P_iv_i^*p'v_i$ is amenable relative to $\mathcal B_i\rtimes\Lambda_{T_i\setminus\{j\}}$ inside $\mathbb M_{m_i}(M)$.
Note that  $v_i^*p'v_i$ is a non-zero projection in $P_i'\cap (p\otimes e_{11})\mathbb M_{m_i}(M)(p\otimes e_{11})=P_i'\cap pMp$. But, recalling that $pMp=A\rtimes\Gamma$, $P_1\overline{\otimes}P_2=L(\Gamma)$, $\Gamma$ is icc, and the action $\Gamma\curvearrowright A$ is ergodic, we get that $\mathcal N_{pMp}(P_i)'\cap pMp\subset L(\Gamma)'\cap A\rtimes\Gamma=\mathbb Cp$.
Thus, by Lemma \ref{facts2}(2), we would get that
$P_i$ amenable relative to $\mathcal B_i\rtimes\Lambda_{T_i\setminus\{j\}}$ inside $\mathbb M_{m_i}(M)$. This contradicts the moreover assertion of Proposition \ref{step1}.

Next, we define $\varphi=\varphi_1\otimes\varphi_2:L(\Gamma)=P_1\overline{\otimes}P_2\rightarrow\mathbb M_{m_1}(M)\overline{\otimes}\mathbb M_{m_2}(M)$. Then $\varphi(L(\Gamma))\subset\mathcal M$, where
 $\mathcal M=\mathcal M_1\overline{\otimes}\mathcal M_2$. 
 We let $v=v_1\otimes v_2$ and note that $\varphi(x)v=vx$, for every $x\in L(\Gamma)$.
 We denote $e=\varphi(1)\in\mathcal M$ and
 $\mathcal B=\mathcal B_1\overline{\otimes}\mathcal B_2$. Then $\mathcal M=\mathcal B\rtimes\Lambda$, where we consider the product action of $\Lambda=\Lambda_{T_1}\times\Lambda_{T_2}$ on $\mathcal B$. The rest of the proof is split between three claims.

\vskip 0.05in
{\bf Claim 1.} If a von Neumann subalgebra $Q\subset L(\Gamma)$ satisfies $\varphi(Q)\prec_{\mathcal M}\mathcal B\rtimes\Lambda_{T_1}$, then $Q\prec_{L(\Gamma)}P_1$. 

\vskip 0.05in
{\it Proof of Claim 1.} Assuming that $Q\nprec_{L(\Gamma)}P_1$, we will prove that $\varphi(Q)\nprec_{\mathcal M}\mathcal B\rtimes\Lambda_{T_1}$. By applying Theorem \ref{corner}, we can find a sequence $u_n\in\mathcal U(Q)$ such that $\|E_{P_1}(u_na)\|_2\rightarrow 0$, for all $a\in L(\Gamma)$. 

For every $i\in\{1,2\}$, let $\psi_i:M\rightarrow\mathbb M_{m_i}(M)$ be the embedding given by $\psi_i(x)=x\otimes e_{11}$. Let $\psi=\psi_1\otimes\psi_2:L(\Gamma)=P_1\overline{\otimes}P_2\rightarrow\mathbb M_{m_1}(M)\overline{\otimes}\mathbb M_{m_2}(M)$. 
We claim that \begin{equation}\label{u_n}\|E_{\mathbb M_{m_1}(M)\overline{\otimes}\mathcal B_2}(a\psi(u_n)b)\|_2\rightarrow 0,\;\;\;\text{for all $a,b\in \mathbb M_{m_1}(M)\overline{\otimes}\mathbb M_{m_2}(M)$}.\end{equation}
By using that $\mathcal B_2=\mathbb M_{m_2}(B)$ and the position of $A\subset B$,  we find $\alpha_1,...,\alpha_{D},\beta_1,...,\beta_{D}\in\mathbb M_{m_2}(M)$ such that  $E_{\mathcal B_2}(x)=\sum_{d=1}^{D}\alpha_dE_{\psi_2(A)}(\alpha_d^*x\beta_d)\beta_d^*$, for every $x\in  \mathbb M_{m_2}(M)$. This allows us to reduce \ref{u_n} to showing that $\|E_{\mathbb M_{m_1}(M)\overline{\otimes}\psi_2(A)}(a\psi(u_n)b)\|_2\rightarrow 0$, for all $a,b\in \mathbb M_{m_1}(M)\overline{\otimes}\psi_2(p)\mathbb M_{m_2}(M)\psi_2(p)$.
Since $\psi_2(p)\mathbb M_{m_2}(M)\psi_2(p)=\psi_2(pMp)=\psi_2(A\rtimes\Gamma)$ and $E_{\mathbb M_{m_1}(M)\overline{\otimes}\psi_2(A)}$ is $\mathbb M_{m_1}(M)\overline{\otimes}\psi_2(A)$-bimodular, it is enough to treat the case when $a=1\otimes\psi_2(\xi)$, $b=1\otimes\psi_2(\zeta)$, for some $\xi,\zeta\in L(\Gamma)$. 

In this case we have $a\psi(u_n)b=(\psi_1\otimes\psi_2)((1\otimes\xi)u_n(1\otimes\zeta))\in P_1\overline{\otimes}\psi_2(L(\Gamma))$. Since  $E_{\psi_2(A)}(\psi_2(x))=\tau(x)\psi_2(p)$, for every $x\in P_2$, we get that $E_{\mathbb M_{m_1}(M)\overline{\otimes}\psi_2(A)}((\psi_1\otimes\psi_2)(x))=(\psi_1\otimes\tau)(x)\psi_2(p)$, for all $x\in L(\Gamma)$. Also, note that $(\psi_1\otimes\tau)(x)=(\psi_1\otimes\tau)(1\otimes E_{P_2})(E_{P_1}\otimes 1)(x)=((\psi_1\otimes\tau)\circ E_{P_1\overline{\otimes}P_2})(x)$, for every $x\in P_1\overline{\otimes}L(\Gamma)$, and that $(1\otimes\xi)u_n(1\otimes\zeta)\in P_1\overline{\otimes}L(\Gamma)$. By combining these fact we get that
\begin{align*}E_{\mathbb M_{m_1}(M)\overline{\otimes}\psi_2(A)}(a\psi(u_n)b)&=  (\psi_1\otimes\tau)((1\otimes\xi)u_n(1\otimes\zeta))\psi_2(p)\\&=((\psi_1\otimes\tau)\circ E_{P_1\overline{\otimes}P_2})((1\otimes\xi)u_n(1\otimes\zeta))\psi_2(p) \\&=\psi_1(E_{P_1}(u_nE_{P_2}(\zeta\xi)))\psi_2(p),\end{align*}
where in order to get that the last equality we used the fact that for all $\alpha\in P_1, \beta\in P_2$ we have \begin{align*}(1\otimes\tau)(E_{P_1}{\otimes}E_{P_2})((1\otimes\xi)(\alpha\otimes\beta)(1\otimes\zeta))&=E_{P_1}(\alpha)\tau(\beta\zeta\xi)\\&=E_{P_1}(\alpha)\tau(\beta E_{P_2}(\zeta\xi))\\&=E_{P_1}((\alpha\otimes\beta)E_{P_2}(\zeta\xi)).\end{align*}

Since $\|E_{P_1}(u_nE_{P_2}(\zeta\xi))\|_2\rightarrow 0$, equation \ref{u_n} follows.

Let $c=c_1c_2>0$. Since $E_{\mathcal M}(vv^*)\geq c\varphi(1)=ce$, if $w=v^*E_{\mathcal M}(vv^*)^{-1}$, then $\varphi(x)=E_{\mathcal M}(vxw)$, for any $x\in L(\Gamma)$.
Let $a, b\in\mathcal M$. Since $\mathcal B\rtimes\Lambda_{T_1}\subset\mathcal M$, we have $E_{\mathcal B\rtimes\Lambda_{T_1}}(a\varphi(u_n)b)=E_{\mathcal B\rtimes\Lambda_{T_1}}(avu_nwb)$. Using that $\mathcal B\rtimes\Lambda_{T_1}=\mathcal M_1\overline{\otimes}\mathcal B_2\subset 
\mathbb M_{m_1}(M)\overline{\otimes}\mathcal B_2$ in combination with \ref{u_n}, the claim follows.
\hfill$\square$

To finish the proof, it suffices to show that Claim 1 applies to $Q=L(\Omega_k)$, for $k$ large enough. This will be achieved  by combining Claims 2 and 3 below.
 We fix $j\in T_2$ and denote $T=\{1,...,n\}\setminus\{j\}$. 
  For $k\geq 1$, we put $N_k=\varphi(L(C_{\Gamma}(\Omega_k)))$, and
  let $f_k\in \mathcal Z(N_k'\cap e\mathcal Me)$ be the maximal projection such that  $N_kf_k$ is amenable relative to $\mathcal B\rtimes\Lambda_T$ inside $\mathcal M$. 

 \vskip 0.05in
 
 {\bf Claim 2.}   $\varphi(L(\Omega_k))(e-f_k)\prec^s_{\mathcal M} \mathcal B\rtimes\Lambda_T$, for any $k\geq 1$.
 
 \vskip 0.05in
 
 {\it Proof of Claim 2.} 
Since $\varphi(L(\Omega_k))\subset N_k'\cap e\mathcal Me$,  by parts (1) and (2) of Lemma \ref{facts1}, it suffices to show that $(N_k'\cap e\mathcal Me)z\prec_{\mathcal M}\mathcal B\rtimes\Lambda_T$, whenever $z\in\mathcal Z((N_k'\cap e\mathcal Me)'\cap e\mathcal Me)(e-f_k)$ is a non-zero projection. 
Since $\mathcal Z((N_k'\cap e\mathcal Me)'\cap e\mathcal Me)\subset\mathcal Z(N_k'\cap e\mathcal Me)$, we get $z\in\mathcal Z(N_k'\cap e\mathcal Me)$. Since $z\leq e-f_k$, the maximality of $f_k$ implies that $N_kz$ is not amenable relative to $\mathcal B\rtimes\Lambda_T$.  Since $(N_k'\cap e\mathcal Me)z$ and $N_kz$ commute, and we can decompose $\mathcal M=(\mathcal B\rtimes\Lambda_T)\rtimes\Lambda_j$, where $\Lambda_j\in\mathcal C_{\text{rss}}$, Lemma \ref{KV15} implies that $(N_k'\cap e\mathcal Me)z\prec_{\mathcal M}\mathcal B\rtimes\Lambda_T$.
 This proves the claim. \hfill$\square$
 
 \vskip 0.05in
 
 Next, put $N=\varphi(L(\cup_{k\geq 1}C_{\Gamma}(\Omega_k)))$. Since $P_2\prec_{L(\Gamma)}L(\cup_{k\geq 1}C_{\Gamma}(\Omega_k))$ by Proposition \ref{ultra}, and $P_2$ is regular in $L(\Gamma)$, Lemma \ref{facts1}(3) implies that $P_2\prec_{L(\Gamma)}^sL(\cup_{k\geq 1}C_{\Gamma}(\Omega_k))$.
Thus $\varphi(P_2)\prec^s_{\varphi(L(\Gamma))}N$, hence $\varphi(P_2)\prec^s_{\mathbb M_{m_1}(M)\overline{\otimes}\mathbb M_{m_2}(M)}N$, so in particular 
 $\varphi(P_2)vv^*\prec_{\mathbb M_{m_1}(M)\overline{\otimes}\mathbb M_{m_2}(M)}N$.
Using Lemma \ref{facts1}(4) we find a non-zero projection $e'\in \mathcal Z(N'\cap e(\mathbb M_{m_1}(M)\overline{\otimes}\mathbb M_{m_2}(M))e)$ such that $\varphi(P_2)vv^*\prec_{\mathbb M_{m_1}(M)\overline{\otimes}\mathbb M_{m_2}(M)} Nf$, for any non-zero projection $f\in N'\cap e(\mathbb M_{m_1}(M)\overline{\otimes}\mathbb M_{m_2}(M))e$ with $f\leq e'$.

 We continue with the following:
 
 \vskip 0.05in

 {\bf Claim 3.}  $\tau(f_ke')\rightarrow 0$, as $k\rightarrow\infty$.

 \vskip 0.05in
 
 {\it Proof of Claim 3.} Assume that the claim is false. Since $N_k\subset N_{k+1}$, we have $f_{k+1}\leq f_k$, for any $k\geq 1$. If $f=\bigwedge_{k}f_k$, then  $f\in\mathcal Z(N'\cap e\mathcal Me)$. Since $f\leq f_k$, we get that $N_kf$ is amenable relative to $\mathcal B\rtimes\Lambda_T$ inside $\mathcal M$, for all $k\geq 1$.
 By  Lemma \ref{union} we get that $Nf=(\cup_{k\geq 1}N_kf)''$ is amenable relative to $\mathcal B\rtimes\Lambda_T$ inside $\mathcal M$. Lemma \ref{facts2}(1) then gives that $Nfe'$ is amenable relative to $\mathcal B\rtimes\Lambda_T$ inside $\mathcal M$.
 

Since $\tau(fe')=\lim\limits_k\tau(f_ke')$ and the claim is assumed false,  $fe'\not=0$. Since  $fe'\leq e'$ belongs to $N'\cap e(\mathbb M_{m_1}(M)\overline{\otimes}\mathbb M_{m_2}(M))e$, the discussion before the claim gives $\varphi(P_2)vv^*\prec_{\mathbb M_{m_1}(M)\overline{\otimes}\mathbb M_{m_2}(M)} Nfe'$.
By Lemma \ref{facts1}(2) there is a non-zero projection $p'\in\mathcal Z((\varphi(P_2)vv^*)'\cap vv^*(\mathbb M_{m_1}(M)\overline{\otimes}\mathbb M_{m_2}(M))vv^*)$ with $\varphi(P_2)p'\prec^s_{\mathbb M_{m_1}(M)\overline{\otimes}\mathbb M_{m_2}(M)} Nfe'$.
 Lemma \ref{facts2}(3) further  gives that $\varphi(P_2)p'$ is amenable relative to $Nfe'$ inside $\mathbb M_{m_1}(M)\overline{\otimes}\mathbb M_{m_2}(M)$. Since $\varphi(P_2)vv^*=v_1v_1^*\otimes\varphi_2(P_2)v_2v_2^*$ and $M$ is a factor, we get that $p'=v_1v_1^*\otimes p''$, for some projection $p''\in\varphi_2(P_2)'\cap\varphi_2(1)\mathbb M_{m_2}(M)\varphi_2(1)$ with $p''\leq v_2v_2^*$.
It follows that $\varphi_2(P_2)p''$ is amenable relative to $Nfe'$ inside $\mathbb M_{m_1}(M)\overline{\otimes}\mathbb M_{m_2}(M)$.

 By combining the conclusions of the last two paragraphs with \cite[Proposition 2.4(3)]{OP07}, we deduce that $\varphi_2(P_2)p''\subset\mathbb M_{m_2}(M)$ is amenable relative to $\mathcal B\rtimes\Lambda_T$ inside $\mathbb M_{m_1}(M)\overline{\otimes}\mathbb M_{m_2}(M)$. Since $\mathcal B\rtimes\Lambda_T$ and $\mathbb M_{m_2}(M)$ are in a commuting square position and regular,  by Lemma \ref{PV11}(2),  $\varphi_2(P_2)p''$ is amenable relative to their intersection, $\mathcal B_2\rtimes\Lambda_{T_2\setminus\{j\}}$, inside $\mathbb M_{m_1}(M)\overline{\otimes}\mathbb M_{m_2}(M)$. 
As $\varphi_2(P_2)p''$ and $\mathcal B_2\rtimes\Lambda_{T_2\setminus\{j\}}$ are  subalgebras of $\mathbb M_{m_2}(M)$, it follows that $\varphi_2(P_2)p''$ is amenable relative to $\mathcal B_2\rtimes\Lambda_{T_2\setminus\{j\}}$ inside $\mathbb M_{m_2}(M)$.
This contradicts the second paragraph of the proof of the lemma. \hfill$\square$
 
 \vskip 0.05in
 
 Next, by combining claims 2 and 3,
  for every $j\in T_2$, we can find projections $f_{k,j}\in \mathcal Z(N_k'\cap e\mathcal Me)$ such that $\varphi(L(\Omega_k))(e-f_{k,j})\prec^s_{\mathcal M} \mathcal B\rtimes\Lambda_{\{1,...,n\}\setminus\{j\}}$, for any $k\geq 1$, and $\tau(f_{k,j}e')\rightarrow 0$, as $k\rightarrow\infty$. 
  
  For $k\geq 1$, let $r_k=\bigvee_{j\in T_2}f_{k,j}$. Then $r_{k}\in \mathcal Z(N_k'\cap e\mathcal Me)$ and 
 since $\tau(r_ke')\leq\sum_{j\in T_2}\tau(f_{k,j}e')$, we get that 
$\tau(r_ke')\rightarrow 0,$ as $k\rightarrow\infty$.
In particular, since $0\not=e'\leq e$, we get that $e-r_k\not=0$, for $k$ large enough.
 On the other hand, since $\varphi(L(\Omega_k))(e-r_k)\prec^s_{\mathcal M}\mathcal B\rtimes\Lambda_{\{1,...,n\}\setminus\{j\}}$, for every $j\in T_2$, and the algebras $\mathcal B\rtimes\Lambda_{\{1,...,n\}\setminus\{j\}}$, with $j\in T_2$, are in a commuting square position and regular in $\mathcal M$, Lemma \ref{PV11}(2) implies that
$\varphi(L(\Omega_k))(e-r_k)\prec^s_{\mathcal M}\mathcal B\rtimes\Lambda_{\{1,...,n\}\setminus T_2}=\mathcal B\rtimes\Lambda_{T_1}$, for any $k\geq 1$.
 
Thus, if $k$ is large enough then $\varphi(L(\Omega_k))\prec_{\mathcal M}\mathcal B\rtimes\Lambda_{T_1}$, hence $L(\Omega_k)\prec_{L(\Gamma)}P_1$, by Claim 1.
 \hfill$\blacksquare$

We are now ready to complete the proof of {\bf Step 4}.

\begin{lemma}\label{Sigma_i} For every $i\in\{1,2\}$ 
we can find a subgroup $\Sigma_i<\Gamma$ such that 
\begin{enumerate}
\item $B\rtimes\Lambda_{S_i}\prec^s_{M} A\rtimes\Sigma_i$.
\item $A\rtimes\Sigma_i\prec^s_{M}B\rtimes\Lambda_{S_i}$.
\item $P_i\prec^s_{L(\Gamma)}L(\Sigma_i)$.
\item $L(\Sigma_i)\prec^s_{L(\Gamma)}P_i$.
\end{enumerate}
\end{lemma}
 
 {\it Proof.}  Assume for simplicity $i=1$.
By Lemma \ref{ultra} we can find a decreasing sequence of subgroups $\Omega_k<\Gamma$ such that
$B\rtimes\Lambda_{S_1}\prec_M A\rtimes\Omega_k$, for all $k\geq 1$, and 
 $P_2\prec_{L(\Gamma)} L(\cup_{k\geq 1}C_{\Gamma}(\Omega_k)).$ By Lemma \ref{Omega_k}, for any $k\geq 1$ large enough $\Sigma_1:=\Omega_k$ satisfies $L(\Sigma_1)\prec_{L(\Gamma)}P_1$ in addition to $B\rtimes\Lambda_{S_1}\prec_MA\rtimes\Sigma_1$. Since $B\rtimes\Lambda_{S_1}$ is regular in the II$_1$ factor $M$, by Lemma \ref{facts1}(3) we get that $B\rtimes\Lambda_{S_1}\prec^s_{M}A\rtimes\Sigma_1$. This proves (1).

By Lemma \ref{facts1}(3), we can find a non-zero projection $e\in L(\Sigma_1)'\cap L(\Gamma)$ with $L(\Sigma_1)e\prec_{L(\Gamma)}^sP_1$. By Proposition \ref{step1} we have that $P_1\prec_M^{s}B\rtimes\Lambda_{T_1}$. By combining these facts with Lemma \ref{facts1}(1) we derive that $L(\Sigma_1)e\prec_{M}^sB\rtimes\Lambda_{T_1}$.
Our next goal is to upgrade this to the following conclusion:

\vskip 0.05in
{\bf Claim 1.} $A\rtimes\Sigma_1\prec_M^sB\rtimes\Lambda_{T_1}$.

\vskip 0.05in

{\it Proof of Claim 1.}  For  $F\subset\Lambda$, let $\mathcal K_F\subset L^2(M)$ be the closed linear span of $\{(B\rtimes\Lambda_{T_1})v_g|g\in F\}$. We denote by
$P_F$ be the orthogonal projection onto $\mathcal K_F$. The proof relies on the following fact: let $R\subset rMr$  be a von Neumann subalgebra and $\mathcal U\subset\mathcal U(R)$ a subgroup with $\mathcal U''=R$. Then $R\prec^s_M B\rtimes\Lambda_{T_1}$ iff for any $\varepsilon>0$, there is  $F\subset\Lambda$ finite such that $\|u-P_F(u)\|_2\leq\varepsilon$, for all $u\in\mathcal U$. This fact follows from \cite[Lemma 2.5]{Va10b} by using that
 $\Lambda_{T_1}<\Lambda$ is a normal subgroup.
 
 Let $\varepsilon>0$.
Since $A\subset pMp$ is maximal abelian and $e\in L(\Gamma)$, we have $E_{A'\cap pMp}(e)=E_A(e)=\frac{\tau(e)}{\tau(p)}p$. On the other hand, $E_{A'\cap pMp}(e)$ belongs to the closed convex hull of $\{vev^*|v\in\mathcal U(A)\}$ (being precisely its element of minimal $\|.\|_2$).
 We can therefore find $v_1,...,v_D,w_1,...,w_D\in\mathcal U(A)$ such that $\|p-\sum_{d=1}^Dv_dew_d\|_2\leq\frac{\varepsilon}{2}$.
Since $L(\Sigma_1)e\prec_{M}^sB\rtimes\Lambda_{T_1}$, by using the above fact, we can find $F\subset\Lambda$ finite such that  $\|u_ge-P_F(u_ge)\|_2\leq\frac{\varepsilon}{2D}$, for any $g\in\Sigma_1$.
 
 By combining the last two inequalities, for every $a\in\mathcal U(A)$ and $g\in\Sigma_1$ we have
 
 $$\|au_g-\sum_{d=1}^Da(u_gv_du_g^*)P_F(u_ge)w_d\|_2\leq\|p-\sum_{d=1}^Dv_dew_d\|_2+D\;\|u_ge-P_F(u_ge)\|_2\leq\varepsilon.$$
 
 Since $\mathcal K_F$ is an $A$-$A$-bimodule, we derive that $\sum_{d=1}^Da(u_gv_du_g^*)P_F(u_ge)w_d\in\mathcal K_F$. Hence, we have $\|au_g-P_F(au_g)\|_2\leq\varepsilon$, for every $a\in\mathcal U(A)$ and $g\in\Sigma_1$. Since $\varepsilon>0$ is arbitrary and the group $\mathcal U=\{au_g|a\in\mathcal U(A), g\in\Sigma_1\}$ generates $A\rtimes\Sigma_1$, the above fact gives the claim.
 \hfill$\square$

By combining the claim with $B\rtimes\Lambda_{S_1}\prec_M A\rtimes\Sigma_1$ and with Lemma \ref{facts1}(1) we conclude that $B\rtimes\Lambda_{S_1}\prec_{M}B\rtimes\Lambda_{T_1}$. This readily implies that $S_1\subset T_1$. By symmetry, we also get that $S_2\subset T_2$. Since $\{S_1, S_2\}$ and $\{T_1,T_2\}$ are partitions of $\{1,...,n\}$ we must have that $S_1=T_1$ and $S_2=T_2$. Thus, Claim 2 reads $A\rtimes\Sigma_1\prec_M^sB\rtimes\Lambda_{S_1}$, which proves (2). 

We are left with proving (3) and (4), which is done in the following two claims.

\vskip 0.05in
{\bf Claim 2.} $P_1\prec^s_{L(\Gamma)}L(\Sigma_1)$.

\vskip 0.05in
{\it Proof of Claim 2.}
Since $P_1$ is regular in $L(\Gamma)$ and $L(\Gamma)$ is a II$_1$ factor, by Lemma \ref{facts1}(3) it suffices to show that $P_1\prec_{L(\Gamma)}L(\Sigma_1)$.
 By Proposition \ref{step1}, $P_1\prec_{M}B\rtimes\Lambda_{T_1}=B\rtimes\Lambda_{S_1}$. 
 By combining this with (1) and  Lemma \ref{facts1}(1), it follows that $P_1\prec_{M}A\rtimes\Sigma_1$. 

Assume by contradiction that $P_1\nprec_{L(\Gamma)}L(\Sigma_1)$.
By Theorem \ref{corner} we can find  $u_n\in\mathcal U(P_1)$ such that $\|E_{L(\Sigma_1)}(au_nb)\|_2\rightarrow 0$, for every $a,b\in L(\Gamma)$.
We claim that $\|E_{A\rtimes\Sigma_1}(au_nb)\|_2\rightarrow 0$, for every $a,b\in pMp=A\rtimes\Gamma$. Since $E_{A\rtimes\Sigma_1}$ is $A$-$A$-bimodular, it suffices to verify this for every $a,b\in L(\Gamma)$. But, since $au_nb\in L(\Gamma)$, we have that $\|E_{A\rtimes\Sigma_1}(au_nb)\|_2=\|E_{L(\Sigma_1)}(au_nb)\|_2\rightarrow 0$. Since the claim implies that $P_1\nprec_{M}A\rtimes\Sigma_1$, we get the desired contradiction.
\hfill$\square$

\vskip 0.05in
{\bf Claim 3.} $L(\Sigma_1)\prec^s_{L(\Gamma)}P_1$.

\vskip 0.05in
{\it Proof of Claim 3.} By Proposition \ref{step2}  we have $\Delta(B\rtimes\Lambda_{S_1})\prec_{M\overline{\otimes}L(\Gamma)}M\overline{\otimes}P_1$. Since $\Gamma$ is icc, we get that $\Delta(M)'\cap M\overline{\otimes}L(\Gamma)=\mathbb C1$. Therefore, by applying Lemma \ref{facts1}(3) we conclude that $\Delta(B\rtimes\Lambda_{S_1})\prec^s_{M\overline{\otimes}L(\Gamma)}M\overline{\otimes}P_1$.
On the other hand, since $L(\Sigma_1)\subset A\rtimes\Sigma_1$, Claim 1 implies that $L(\Sigma_1)\prec_{M}^sB\rtimes\Lambda_{S_1}$, and therefore $\Delta(L(\Sigma_1))\prec_{M\overline{\otimes}L(\Gamma)}^s\Delta(B\rtimes\Lambda_{S_1})$. 
By combining these facts with Lemma \ref{facts1}(1), we derive that $\Delta(L(\Sigma_1))\prec^s_{M\overline{\otimes}L(\Gamma)}M\overline{\otimes}P_1$.

Let $p'\in L(\Sigma_1)'\cap L(\Gamma)$ be a non-zero projection. Assuming that 
 $L(\Sigma_1)p'\nprec_{L(\Gamma)}P_1$, we will reach a contradiction, which will prove the claim.
 By Theorem \ref{corner} we can find a sequence $g_n\in\Sigma_1$ such that $\|E_{P_1}(au_{g_n}p'b)\|_2\rightarrow 0$, for every $a,b\in L(\Gamma)$. We claim that $\|E_{M\overline{\otimes}P_1}(a\Delta(u_{g_n})(1\otimes p')b)\|_2\rightarrow 0$, for every $a,b\in M\overline{\otimes}L(\Gamma)$. Since $\Delta(u_{g_n})\in\mathcal U(\Delta(L(\Sigma_1)))$ and $1\otimes p'\in\Delta(L(\Sigma_1))'\cap M\overline{\otimes}L(\Gamma)$ is  non-zero projection (recall that $\Delta(u_g)=u_g\otimes u_g$, for all $g\in\Gamma$), we get that $\Delta(L(\Sigma_1))(1\otimes p')\nprec_{M\overline{\otimes}L(\Gamma)}M\overline{\otimes}P_1$, which contradicts the conclusion of the previous paragraph. Thus, it remains to prove the claim. 
 
 Since $E_{M\overline{\otimes}P_1}$ is $M\otimes 1$-$M\otimes 1$-bimodular, we may assume that $a,b\in 1\otimes L(\Gamma)$. But in this case we have $\|E_{M\overline{\otimes}P_1}(a\Delta(u_{g_n})(1\otimes p')b)\|_2=\|E_{P_1}(au_{g_n}p'b)\|_2\rightarrow 0$, which finishes the proof.
 \hfill$\blacksquare$

\subsection{Step 5: completion of the proof of Theorem \ref{tensor}}
Let $i\in\{1,2\}$. By Lemma \ref{Sigma_i} we have that $B\rtimes\Lambda_{S_i}\prec^s_MA\rtimes\Sigma_i$ and $A\rtimes\Sigma_i\prec_{M}^sB\rtimes\Lambda_{S_i}$.
Recalling that $A=L^{\infty}(X)$ and $B\rtimes\Lambda_{S_i}=(L^{\infty}(Y)\rtimes\Lambda_{S_i})\otimes\mathbb M_{\ell}(\mathbb C)$, we get that $L^{\infty}(Y)\rtimes\Lambda_{S_i}\prec_{M}L^{\infty}(X)\rtimes\Sigma_i$ and  also that $L^{\infty}(X)\rtimes\Sigma_i\prec^s_{M}L^{\infty}(Y)\rtimes\Lambda_{S_i}$. Since $\Lambda=\Lambda_{S_1}\times\Lambda_{S_2}$, Proposition \ref{P: twine2ME} and implies that $\Sigma_i$ is measure equivalent to $\Lambda_{S_i}$. Together with Lemma \ref{Sigma_i}, this finishes the proof of Theorem \ref{tensor}. \hfill$\blacksquare$

\section{From tensor decompositions to product decompositions}

The goal of this section is prove the following result that we will need in the proof of Theorem \ref{C}. We say that two subgroups $\Sigma,\Omega$ of a countable group $\Gamma$ are called {\it commensurable} if we have that $[\Sigma:\Sigma\cap\Omega]<\infty$ and $[\Omega:\Sigma\cap\Omega]<\infty$.

\begin{theorem}\label{product}
Let $\Gamma$ be a countable icc group, denote $M=L(\Gamma)$, and assume that $M=P_1\overline{\otimes}P_2$. For every $i\in\{1,2\}$, let $\Sigma_i<\Gamma$ be a subgroup such that $P_i\prec_M^sL(\Sigma_i)$ and $L(\Sigma_i)\prec_M^sP_i$.

Then we can find a decomposition $\Gamma=\Gamma_1\times\Gamma_2$, a decomposition $M=P_1^s\overline{\otimes}P_2^{1/s}$, for some $s>0$, and a unitary $u\in\mathcal U(M)$ such that
\begin{itemize}
\item $\Gamma_1$ is commensurable to $k\Sigma_1k^{-1}$, for some $k\in\Gamma$, $\Gamma_2$ is commensurable to $\Sigma_2$, 
\item $P_1^s=uL(\Gamma_1)u^*$ and $P_2^{1/s}=uL(\Gamma_2)u^*$.
\end{itemize}
\end{theorem}

The proof of Theorem \ref{product} 
relies on several results. Before continuing, we introduce some terminology. Let $\Gamma$ be a countable group and $\Sigma<\Gamma$ be a subgroup. Following \cite{CdSS15},  we denote by $\mathcal O_{\Sigma}(g)=\{hgh^{-1}|h\in\Sigma\}$ the orbit of $g\in\Gamma$ under the conjugation action of $\Sigma$. 
Note that $\mathcal O_{\Sigma}(g_1g_2)\subset\mathcal O_{\Sigma}(g_1)\mathcal O_{\Sigma}(g_2)$, thus $|\mathcal O_{\Sigma}(g_1g_2)|\leq |\mathcal O_{\Sigma}(g_1)| |\mathcal O_{\Sigma}(g_2)|$, for all $g_1,g_2\in\Gamma$. 
Therefore, the set $\Delta=\{g\in\Gamma\;|\;\text{$\mathcal O_{\Sigma}(g)$ is finite}\}$ is a subgroup of $\Gamma$. Moreover, we note that $L(\Sigma)'\cap L(\Gamma)\subset L(\Delta)$.

\subsection{From commuting subalgebras to almost commuting subgroups} The first step towards proving Theorem \ref{product} is to show the existence of conjugates of finite index subgroups of $\Sigma_1,\Sigma_2$ that ``almost'' commute, in the sense that they have finite commutator.

\begin{theorem}\label{commute}
Let $\Gamma$ be a countable group and $\Sigma_1,\Sigma_2<\Gamma$ be two subgroups. Assume that we have $L(\Sigma_1)\prec_{L(\Gamma)}L(\Sigma_2)'\cap L(\Gamma)$.

Then we can find finite index subgroups $\Omega_1<k\Sigma_1k^{-1}$ and $\Omega_2<\Sigma_2$, for some $k\in\Gamma$, such that the group $[\Omega_1,\Omega_2]$ generated by all commutators $[g,h]=g^{-1}h^{-1}gh$ with $g\in\Omega_1,h\in\Omega_2$, is finite and satisfies $[\Omega_1,\Omega_2]\subset C_{\Gamma}(\Omega_1)\cap C_{\Gamma}(\Omega_2).$
\end{theorem}

\begin{remark}
We do not know whether the following more natural, stronger conclusion holds: there exist finite index commuting subgroups $\Omega_1<k\Sigma_1k^{-1}$ and $\Omega_2<\Sigma_2$, for some $k\in\Gamma$. Note, however, that Lemma \ref{commute2} below implies that this is the case if $\Sigma_1$ is finitely generated.
\end{remark}

The proof of Theorem \ref{commute} relies on the following lemma inspired by \cite[Claims 4.9-4.11]{CdSS15}.

\begin{lemma}\label{commute2} Assume the setting of Theorem \ref{commute}. Let $\Delta=\{g\in\Gamma\;|\;\text{$\mathcal O_{\Sigma_2}(g)$ is finite}\}$.

Then we can find a finite index subgroup $\Omega_1<\Sigma_1$ and $k\in\Gamma$ such that $k\Omega_1 k^{-1}\subset\Delta$ and $L(k\Omega_1k^{-1})\prec_{L(\Delta)}L(\Sigma_2)'\cap L(\Gamma)$.
\end{lemma}

{\it Proof.} Since $L(\Sigma_1)\prec_{L(\Gamma)}L(\Sigma_2)'\cap L(\Gamma)$, by Theorem \ref{corner} we can find $k_1,...,k_n,l_1,...,l_n\in\Gamma$ and a constant $\delta>0$ such that 
\begin{equation}\label{embed}\sum_{i=1}^n\|E_{L(\Sigma_2)'\cap L(\Gamma)}(u_{k_i}u_gu_{l_i})\|_2^2\geq\delta,\;\;\;\;\text{for every $g\in\Sigma_1$.}\end{equation}
If $g\in\Gamma$, then $E_{L(\Sigma_2)'\cap L(\Gamma)}(u_g)$ is equal to $\frac{1}{|\mathcal O_{\Sigma_2}(g)|}\sum_{h\in\mathcal O_{\Sigma_2}(g)}u_h$, if $g\in\Delta$, and to $0$, otherwise. Thus, we have $\|E_{L(\Sigma_2)'\cap L(\Gamma)}(u_g)\|_2^2=\frac{1}{|\mathcal O_{\Sigma_2}(g)|}$, where we make the convention that $\frac{1}{\infty}=0$. 

Let $c=\frac{n}{\delta}$ and define $S=\{g\in\Gamma| \;|\mathcal O_{\Sigma_2}(g)|\leq c\}$. 
By using  \ref{embed} we get that for any $g\in\Sigma_1$, there is $i\in\{1,...,n\}$ such that $|\mathcal O_{\Sigma_2}(k_igl_i)|\leq c$.
Hence, we have $\Sigma_1\subset\cup_{i=1}^nk_i^{-1}Sl_i^{-1}$.
For $i\in\{1,...,n\}$, let $a_i\in \Sigma_1\cap k_i^{-1}Sl_i^{-1}$,
if $\Sigma_1\cap k_i^{-1}Sl_i^{-1}$ is non-empty, and let $a_i=e$, otherwise.

Since $S\subset\Delta$, we get that 
$\Sigma_1\subset\cup_{i=1}^n(k_i^{-1}\Delta k_i)a_i. $
This implies that at least one of the groups $\Sigma_1\cap k_i^{-1}\Delta k_i$, with $1\leq i\leq n$, has finite index in $\Sigma_1$. 
After renumbering, we find $m\in\{1,...,n\}$ such that the index $[\Sigma_1:\Sigma_1\cap k_i^{-1}\Delta k_i]$ is finite, for all $1\leq i\leq m$, and infinite, for all $m+1\leq i\leq n$. 

Define $\Omega_1=\cap_{i=1}^m(\Sigma_1\cap k_i^{-1}\Delta k_i)$. Then $\Omega_1$ has finite index in $\Sigma_1$, and  $\Omega_1=\cup_{i=1}^n(\Omega_1\cap k_i^{-1}Sl_i^{-1})$.
For $1\leq i\leq n$, let $b_i\in \Omega_1\cap k_i^{-1}Sl_i^{-1}$,
if $\Omega_1\cap k_i^{-1}Sl_i^{-1}$ is non-empty, and $b_i=e$, otherwise. If $i\leq m$, then since $b_i\in\Omega_1\subset k_i^{-1}\Delta k_i$, we get $k_ib_ik_i^{-1}\in\Delta$, or equivalently $|\mathcal O_{\Sigma_2}(k_ib_ik_i^{-1})|<\infty$.  Let $0<d\leq 1$ be a constant such that $d\leq\frac{1}{c^2|\mathcal O_{\Sigma_2}(k_ib_ik_i^{-1})|}$, for every $1\leq i\leq m$.

Next, fix $g\in\Omega_1$. Then $g\in\Omega_1\cap k_i^{-1}Sl_i^{-1}$, for some $1\leq i\leq n$. 
Thus, $gb_i^{-1}\in k_i^{-1}SS^{-1}k_i$ and hence $gb_i^{-1}\in\Omega_1\cap k_i^{-1}\Delta k_i$. 
Moreover, since $k_igb_i^{-1}k_i^{-1}\in SS^{-1}$, we get that $k_igb_i^{-1}k_i^{-1}\in\Delta$ and that $|\mathcal O_{\Sigma_2}(k_igb_i^{-1}k_i^{-1})|\leq c^2$. 
Now, if $i\leq m$, then  $|\mathcal O_{\Sigma_2}(k_igk_i^{-1})|\leq c^2|\mathcal O_{\Sigma_2}(k_ib_ik_i^{-1})|\leq\frac{1}{d}$, hence $\|E_{L(\Sigma_2)'\cap L(\Gamma)}(u_{k_igk_i^{-1}})\|_2^2\geq d$. Altogether, since $d\leq 1$, we conclude that 
 \begin{equation}\sum_{i=1}^m\|E_{L(\Sigma_2)'\cap L(\Gamma)}(u_{k_igk_i^{-1}})\|_2^2+\sum_{i=m+1}^n\|E_{L(\Omega_1\cap k_i^{-1}\Delta k_i)}(u_gu_{b_i}^*)\|_2^2 \geq d,\;\;\;\;\text{for every $g\in\Omega_1$.}\end{equation}
Since $k_i\Omega_1 k_i^{-1}\subset\Delta$ and $b_i\in\Omega_1$, Remark \ref{corner2} implies that either $L(k_i\Omega_1 k_i^{-1})\prec_{L(\Delta)}L(\Sigma_2)'\cap L(\Gamma)$, for some $1\leq i\leq m$, or that $L(\Omega_1)\prec_{L(\Omega_1)}L(\Omega_1\cap k_i^{-1}\Delta k_i)$, for some $m+1\leq i\leq n$. The latter is however impossible by Lemma \ref{groups}(1) since the inclusion $\Omega_1<\Sigma_1$ has finite index and thus the inclusion $\Omega_1\cap k_i^{-1}\Delta k_i<\Omega_1$ has infinite index, for every $m+1\leq i\leq n$. This proves the lemma.
\hfill$\blacksquare$

{\bf Proof of Theorem \ref{commute}}. Let $\Delta=\{g\in\Gamma\;|\;\text{$\mathcal O_{\Sigma_2}(g)$ is finite}\}$.
By Lemma \ref{commute2}, we can find a finite index subgroup $\Omega_1<k\Sigma_1k^{-1}$, for some $k\in\Gamma$, such that $\Omega_1\subset\Delta$ and $L(\Omega_1)\prec_{L(\Delta)}L(\Sigma_2)'\cap L(\Gamma)$. 
We continue with the following claim. If $A\subset pL(\Gamma)p$ and $B\subset L(\Gamma)$ are von Neumann subalgebras, then we write $A\subset_{\varepsilon}B$ if $\|a-E_B(a)\|_2\leq\varepsilon$, for every $a\in A$ with $\|a\|\leq 1$.

{\bf Claim.} There exists a non-zero projection $z\in L(\Omega_1)'\cap L(\Delta)$ with the following property: for every $\varepsilon>0$ we can find a finite index  subgroup $\Omega_2<\Sigma_2$ such that $L(\Omega_1)z\subset_{\varepsilon}L(\Omega_2)'\cap L(\Gamma)$.

{\it Proof of the claim.} By Theorem \ref{corner} we can find projections $p\in L(\Omega_1), q\in L(\Sigma_2)'\cap L(\Gamma)$, a non-zero partial isometry $v\in qL(\Delta)p$, and a $*$-homomorphism $\theta:pL(\Omega_1)p\rightarrow q(L(\Sigma_2)'\cap L(\Gamma))q$  such that $vx=\theta(x)v$, for every $x\in pL(\Omega_1)p$.
Since $v^*v\in (pL(\Omega_1)p)'\cap pL(\Delta)p$, we can find a projection $p'\in L(\Omega_1)'\cap L(\Delta)$ such that $v^*v=pp'$. Let $p''\in\mathcal Z(L(\Omega_1))$ be the central support of $p$. 

We will prove that $z=p''p'$ satisfies the claim.
To this end, fix $\varepsilon>0$ and $x\in L(\Omega_1)$ with $\|x\|\leq 1$. Let $v_i\in L(\Omega_1)$ be partial isometries such that $p''=\sum_{i\geq 1}v_iv_i^*$ and $v_i^*v_i\leq p$, for every $i\geq 1$.
Let $n\geq 1$ such that $\|p''-\sum_{i=1}^nv_iv_i^*\|_2\leq\frac{\varepsilon}{4}$. Then
 \begin{equation}\label{p'}\|xp''p'-\sum_{i,j=1}^nv_iv_i^*xv_jv_j^*p'\|_2\leq \|xp''-\sum_{i,j=1}^nv_iv_i^*xv_jv_j^*\|_2\leq\frac{\varepsilon}{2}.\end{equation}
On the other hand, using that  $v_j$ and $p'$ commute, for every $j$,  that $v_i^*xv_j\in pL(\Omega_1)p$, for every $i,j$, and that $yp'=v^*\theta(y)v$, for every $y\in pL(\Omega_1)p$, we derive that $$\sum_{i,j=1}^nv_iv_i^*xv_jv_j^*p'=\sum_{i,j=1}^nv_iv_i^*xv_jp'v_j^*=\sum_{i,j=1}^nv_iv^*\theta(v_i^*xv_j)vv_j^*.$$ Now, if $g\in\Delta$, then $\mathcal O_{\Sigma_2}(g)$ is finite, hence $g$ commutes with a finite index subgroup of $\Sigma_2$. Therefore,  any finite subset of $\Delta$ commutes with some finite index subgroup of $\Sigma_2$. This implies that for every $y\in L(\Delta)$ and $\delta>0$, we can find a finite index subgroup $\Omega_2<\Sigma_2$ such that $\|y-E_{L(\Omega_2)'\cap L(\Delta)}(y)\|_2\leq \delta$. 

Thus, there is a finite index subgroup $\Omega_2<\Sigma_2$ such that $\|v_iv^*-E_{L(\Omega_2)'\cap L(\Delta)}(v_iv^*)\|_2\leq\frac{\varepsilon}{4n^2}$, for all $1\leq i\leq n$. Using these inequalities and the last displayed formula, it follows that 
\begin{equation}\label{p''}\|\sum_{i,j=1}^nv_iv_i^*xv_jv_j^*p'-\sum_{i,j=1}^nE_{L(\Omega_2)'\cap L(\Delta)}(v_iv^*)\theta(v_i^*xv_j)E_{L(\Omega_2)'\cap L(\Delta)}(v_jv^*)^*\|_2\leq\frac{\varepsilon}{2}. \end{equation}
Since $\sum_{i,j=1}^nE_{L(\Omega_2)'\cap L(\Delta)}(v_iv^*)\theta(v_i^*xv_j)E_{L(\Omega_2)'\cap L(\Delta)}(v_jv^*)^*$ belongs to $L(\Omega_2)'\cap L(\Delta)$, by combining \ref{p'} and \ref{p''} we deduce that $\|xp''p'-E_{L(\Omega_2)'\cap L(\Delta)}(xp''p')\|_2\leq\varepsilon$. Since $x\in L(\Omega_1)$ with $\|x\|\leq 1$ is arbitrary, the claim follows.
\hfill$\square$

Now, write $z=\sum_{g\in\Delta}c_gu_g$, where $c_g\in\mathbb C$. Let $\alpha=\max_{g\in\Delta}|c_g|$ and put $F=\{g\in\Delta| \;|c_g|=\alpha\}$. Then $F$ is a finite set, and there is $\varepsilon>0$ such that if $k\in\Gamma$ satisfies $\|u_kz-z\|_2<\varepsilon$, then $kF=F$. Indeed, one can check that $\varepsilon=\alpha-\beta$, where $\beta=\max_{g\in\Delta\setminus F}|c_g|$, works. 

The claim gives a finite index subgroup $\Omega_2<\Sigma_2$ such that $L(\Omega_1)z\subset_{\frac{\varepsilon}{4}}L(\Omega_2)'\cap L(\Gamma)$. As $z\in L(\Delta)$, after replacing $\Omega_2$ with a finite index subgroup, we may assume that $\|z-E_{L(\Omega_2)'\cap L(\Delta)}(z)\|_2<\frac{\varepsilon}{4}$. Let $g\in\Omega_1$ and $h\in\Omega_2$. Since $\|u_gz-E_{L(\Omega_2)'\cap L(\Delta)}(u_gz)\|_2\leq\frac{\varepsilon}{4}$, we get that $\|u_gz-u_h(u_gz)u_h^*\|_2\leq\frac{\varepsilon}{2}$. Since $\|z-E_{L(\Omega_2)'\cap L(\Delta)}(z)\|_2<\frac{\varepsilon}{4}$, we also have that $\|zu_h^*-u_h^*z\|_2<\frac{\varepsilon}{2}$. Altogether, we deduce that $\|u_gz-u_hu_gu_h^*z\|_2<\varepsilon$, hence $\|z-u_{g^{-1}hgh^{-1}}z\|_2<\varepsilon$. By the previous paragraph, this implies that $g^{-1}hgh^{-1}F=F$, for every $g\in\Omega_1$ and $h\in\Omega_2$. 

Therefore, $[\Omega_1,\Omega_2]$ is finite and contained in the group  $\langle F\rangle$ generated by $F$.
 Since $z\in L(\Omega_1)'\cap L(\Delta)$ and $F\subset\Delta$, after replacing $\Omega_1,\Omega_2$ with finite index subgroups, we may assume that they commute with $F$.  Thus, $[\Omega_1,\Omega_2]$ is finite and $[\Omega_1,\Omega_2]\subset\langle F\rangle\subset C_{\Gamma}(\Omega_1)\cap C_{\Gamma}(\Omega_2)$. This finishes the proof.
 \hfill$\blacksquare$

\subsection{Finite index commensurator} The next step towards proving Theorem \ref{product} is to show that $\Sigma_i$ is commensurated by a finite index subgroup of $\Gamma$, for every $i\in\{1,2\}$.

\begin{lemma}\label{finiteindex}
Let $\Gamma$ be a countable icc group, denote $M=L(\Gamma)$, and assume that $M=P_1\overline{\otimes}P_2$. Let $\Sigma<\Gamma$ be a subgroup such that $P_1\prec_M^sL(\Sigma)$ and $L(\Sigma)\prec_M^sP_1$.  Let $\Gamma_0<\Gamma$ be the subgroup of $g\in\Gamma$ such that $\Sigma$ and $g\Sigma g^{-1}$ are commensurable.

Then $[\Gamma:\Gamma_0]<\infty$.
\end{lemma}

{\it Proof.} The proof is inspired by \cite[Claims 4.5 and 4.6]{CdSS15}.
Let $\Delta=\{g\in\Gamma|\;\text{$\mathcal O_{\Sigma}(g)$ is finite}\}$. Then $\Delta\subset\Gamma_0$, hence $\Sigma\Delta\subset\Gamma_0$. Indeed, if $k\in\Delta$, then $k$ commutes with a finite index subgroup of $\Sigma$, hence $k\in\Gamma_0$.
Since $L(\Sigma)'\cap M\subset L(\Delta)$, we have $L(\Sigma)\vee (L(\Sigma)'\cap M)\subset L(\Sigma\Delta)\subset L(\Gamma_0)$.
Lemma \ref{groups}(1) implies that in order to reach the conclusion it is sufficient to prove that
\begin{equation}\label{vee}M\prec L(\Sigma)\vee (L(\Sigma)'\cap M).\end{equation}
Towards proving \ref{vee}, we denote $Q_1=L(\Sigma)$.  Then the hypothesis gives that $P_1\prec^sQ_1$ and $Q_1\prec^sP_1$.  By Lemma \ref{facts1}(4), there is a non-zero projection $z\in\mathcal Z(Q_1'\cap M)$ such that $P_1\prec Q_1q'$ for every non-zero projection $q'\in (Q_1'\cap M)z$.
 We claim that 
  $(Q_1'\cap M)z\prec^sP_2$.
By Lemma \ref{facts1}(2), it suffices to show that $(Q_1'\cap M)z'z\prec P_2$, for any projection $z'\in\mathcal Z((Q_1'\cap M)'\cap M)$ such that $z'z\not=0$.
  But $z'z\in (Q_1'\cap M)z$, and thus by the above $P_1\prec Q_1z'z$. By \cite[Lemma 3.5]{Va08} we derive that $(Q_1'\cap M)z'z\prec P_2=P_1'\cap M$, which proves our claim.
  
 Next, we denote $Q_2=Q_1'\cap M$. 
Then  $z\in\mathcal Z(Q_2)$ is a non-zero projection such that $Q_2z\prec^sP_2$ and $P_1\prec Q_1q'$, for every non-zero projection $q'\in Q_2z$. Since we also have that $Q_1z\prec^sP_1$, we get that $\mathcal Z(Q_1)z=Q_1z\cap Q_2z$ satisfies $\mathcal Z(Q_1)z\prec^sP_1$ and $\mathcal Z(Q_1)z\prec^sP_2$. By Lemma \ref{PV11}(2) we deduce that $\mathcal Z(Q_1)z\prec^s P_1\cap P_2=\mathbb C1$, hence $\mathcal Z(Q_1)z$ is completely atomic.

Further, since $Q_1z\prec P_1$, by \cite[Lemma 3.5]{Va08} we get that $P_2=P_1'\cap M\prec Q_2z$. By arguing as in the second paragraph, we can find a non-zero projection $z'\in\mathcal Z((Q_2z)'\cap zMz)=\mathcal Z(Q_2'\cap M)z$ such that $(Q_2'\cap M)z'\prec^s P_1$. Since $\mathcal Z(Q_2'\cap M)\subset\mathcal Z(Q_2)$, we have that $z'\in\mathcal Z(Q_2)z$. Since $Q_2z'\prec^sP_2$, by arguing as in the previous paragraph, we get that $\mathcal Z(Q_2)z'$ is completely atomic.

Thus, $z'\in\mathcal Z(Q_2)$ is a non-zero projection such that $\mathcal Z(Q_1)z'$ and $\mathcal Z(Q_2)z'$ are completely atomic. By shrinking $z'$ we may assume that in fact $\mathcal Z(Q_2)z'=\mathbb Cz'$. Since $\mathcal Z(Q_1)z'$ is completely atomic we can find a non-zero projection $f\in\mathcal Z(Q_1)z'$ such that $\mathcal Z(Q_1)f=\mathbb Cf$. But then also $\mathcal Z(Q_2)f=\mathbb Cf$. Therefore, $f\in Q_2=Q_1'\cap M$ is a projection such that both $Q_1f$ and $fQ_2f=(Q_1f)'\cap fMf$ are II$_1$ factors. Since $Q_1f\prec P_1$, by \cite[Proposition 12]{OP03}, we can find a decomposition $fMf=P_1^{t_1}\overline{\otimes}P_2^{t_2}$, for some $t_1,t_2>0$ satisfying $t_1t_2=\tau(f)$,  and a unitary $u\in fMf$ such that $uQ_1fu^*\subset P_1^{t_1}$.

Since $f\in Q_2z$ is a non-zero projection, we have that $P_1\prec Q_1f$, hence $P_1^{t_1}\prec_{fMf}uQ_1fu^*$. We claim that $P_1^{t_1}\prec_{P_1^{t_1}}uQ_1fu^*$.  Otherwise we can find a sequence $u_n\in\mathcal U(P_1^{t_1})$ such that $\|E_{uQ_1fu^*}(au_nb)\|_2\rightarrow 0$, for every $a,b\in P_1^{t_1}$. We will show that $\|E_{uQ_1fu^*}(a_0u_nb_0)\|_2\rightarrow 0$, for every $a_0,b_0\in fMf$, contradicting the fact that $P_1^{t_1}\prec_{fMf}uQ_1fu^*$. Since $fMf=P_1^{t_1}\overline{\otimes}P_2^{t_2}$, we may assume that $a_0=a_1\otimes a_2$ and $b_0=b_1\otimes b_2$, for 
$a_1,a_2\in P_1^{t_1}$ and $b_1,b_2\in P_2^{t_2}$.
Using that $uQ_1fu^*\subset P_1^{t_1}$, we get that  $\|E_{uQ_1fu^*}(a_0u_nb_0)\|_2=\|E_{uQ_1fu^*}(a_1u_nb_1\otimes a_2b_2)\|_2=\|E_{uQ_1fu^*}(a_1u_nb_1)\|_2|\tau(a_2b_2)|\rightarrow 0$. This altogether proves that $P_1^{t_1}\prec_{P_1^{t_1}}uQ_1fu^*$.

This implies that $fMf=P_1^{t_1}\overline{\otimes}P_2^{t_2}\prec_{fMf}uQ_1fu^*\overline{\otimes}P_2^{t_2}$. Since $P_2^{t_2}\subset (uQ_1fu^*)'\cap fMf$, we get that $fMf\prec_{fMf}Q_1f\vee ((Q_1f)'\cap fMf)=f(Q_1\vee (Q_1'\cap M))f$, which proves \ref{vee} and the lemma.
\hfill$\blacksquare$

\subsection{Proof of Theorem \ref{product}} The proof Theorem \ref{product} has two main parts.

In the first part of the proof, we construct two commuting icc subgroups $\Omega_1,\Omega_2<\Gamma$ which are conjugates of finite index subgroups of $\Sigma_1,\Sigma_2$, and satisfy $[\Gamma:\Omega_1\Omega_2]<\infty$ (compare with \cite[Theorem 4.3]{CdSS15}).

Since $L(\Sigma_2)\prec P_2$, \cite[Lemma 3.5]{Va08} implies that $P_1\prec L(\Sigma_2)'\cap M$. Since $P_1$ is regular in $M$ and $M$ is a II$_1$ factor, Lemma \ref{facts1}(3) implies that $P_1\prec^sL(\Sigma_2)'\cap M$. Since $L(\Sigma_1)\prec P_1$, by combining this with Lemma \ref{facts1}(1) we deduce that $L(\Sigma_1)\prec L(\Sigma_2)'\cap M$. 

By applying Theorem \ref{commute}, we find finite index subgroups $\Omega_1<k\Sigma_1k^{-1}$, $\Omega_2<\Sigma_2$, for some $k\in\Gamma$, such that $[\Omega_1,\Omega_2]$ is finite and contained in $C_{\Gamma}(\Omega_1)\cap C_{\Gamma}(\Omega_2)$.
 If $i\in\{1,2\}$, 
then  Lemma \ref{groups}(2) implies that $L(\Omega_i)\prec^sL(\Sigma_i)$ and $L(\Sigma_i)\prec^sL(\Omega_i)$. Since $L(\Sigma_i)\prec^sP_i$ and $P_i\prec^sL(\Sigma_i)$, we conclude that $L(\Omega_i)\prec^sP_i$ and $P_i\prec^s L(\Omega_i)$.

By applying Lemma \ref{finiteindex} to $\Omega_1$ we deduce that $[\Gamma:\Gamma_0]<\infty$, where $\Gamma_0<\Gamma$ is the subgroup of $g\in\Gamma$ such that $\Omega_1$ and $g\Omega_1g^{-1}$ are commensurable. Since $[\Gamma:\Gamma_0]<\infty$ and $\Gamma$ is icc, it follows that $\mathcal O_{\Gamma_0}(g)$ is infinite, for every $g\in\Gamma\setminus\{e\}$. From this we deduce that $L(\Gamma_0)'\cap M=\mathbb C1$. 
Using that $P_1\prec L(\Omega_1)$ and $P_2\prec L(\Omega_2)$, we find non-zero elements $v,v_1,...,v_m,w,w_1,...,w_m\in M$  such that 
\begin{equation}\label{module} (P_1)_1v\subset\sum_{i=1}^mv_i(L(\Omega_1))_1\;\;\;\;\text{and}\;\;\;\;w(P_2)_1\subset\sum_{i=1}^m (L(\Omega_2))_1w_i. \end{equation}
We claim that we can find $g\in\Gamma_0$ such that $vu_gw\not=0$. Indeed, otherwise we would get that $u_g^*v^*vu_gww^*=0$, for every $g\in\Gamma_0$. Thus,  if $K$ denotes the $\|.\|_2$-closure of the convex hull of $\{u_g^*v^*vu_g|g\in\Gamma_0\}$, then $\xi ww^*=0$, for all $\xi\in K$. Let $\eta\in K$ be the unique element of minimal $\|.\|_2$. Since the map $K\ni \xi\mapsto u_h^*\xi u_h\in K$ preserves $\|.\|_2$, we get that $u_h^*\eta u_h=\eta$, for all $h\in\Gamma_0$. Thus, $\eta\in L(\Gamma_0)'\cap M = \mathbb C1$ and since $\tau(\eta)=\tau(v^*v)$, we deduce that $\eta=\tau(v^*v)1$. But this implies that $0=\eta ww^*=\tau(v^*v)ww^*$, contradicting that both $v$ and $w$ are non-zero. This proves the claim.

Next, since $[\Omega_1:\Omega_1\cap g\Omega_1g^{-1}]<\infty$, we can find $g_1,...,g_n\in\Gamma$ such that $\Omega_1g\subset\cup_{j=1}^ng_j\Omega_1$ and thus $(L(\Omega_1))_1u_g\subset\sum_{j=1}^nu_{g_j}(L(\Omega_1))_1$. By combining this inclusion with equation \ref{module} we get that
\begin{equation}\label{module2}(P_1)_1(vu_gw)(P_2)_1\subset\sum_{i=1}^m\sum_{j=1}^nv_iu_{g_j}(L(\Omega_1))_1(L(\Omega_2))_1w_i. \end{equation}
Thus, if we denote by $\Omega<\Gamma$ the subgroup generated by $\Omega_1$ and $\Omega_2$, then \ref{module2} implies that 
\begin{equation}\label{module3}\mathcal U(P_1)\;(vu_gw)\;\mathcal U(P_2)\subset\sum_{i=1}^m\sum_{j=1}^nv_iu_{g_j}(L(\Omega))_1w_i. \end{equation}
Let us show that $[\Gamma:\Omega]<\infty$. Otherwise, if $[\Gamma:\Omega]=\infty$,  Lemma \ref{groups}(1) implies that $M\nprec L(\Omega)$. Since the group of unitaries $\{u_1\otimes u_2|u_1\in\mathcal U(P_1),u_2\in\mathcal U(P_2)\}$ generates $M$, by Theorem \ref{corner} we can find a sequence $u_n=u_{n,1}\otimes u_{n,2}$, with $u_{n,1}\in\mathcal U(P_1)$ and $u_{n,2}\in\mathcal U(P_2)$, such that $\|E_{L(\Omega)}(au_nb)\|_2\rightarrow 0$, for every $a,b\in M$. We claim that $\|E_{L(\Omega)}(au_{n,1}bu_{n,2}c)\|_2\rightarrow 0$, for every $a,b,c\in M$. 
Since this claim contradicts equation \ref{module3}, we conclude that the assumption $[\Gamma:\Omega]=\infty$ is false.
To prove this claim, we may assume that $a=a_1\otimes a_2, b=b_1\otimes b_2, c=c_1\otimes c_2$, where $a_1,b_1,c_1\in P_1$ and $a_2,b_2,c_2\in P_2$. But in this case $au_{n,1}bu_{n,2}c=a_1u_{n,1}b_1c_1\otimes a_2b_2u_{n,2}c_2=(a_1\otimes a_2b_2)u_n(b_1c_1\otimes c_2)$, and therefore $\|E_{L(\Omega)}(au_{n,1}bu_{n,2}c)\|_2\rightarrow 0$ by the above. 

Since $\Gamma$ is icc and $[\Gamma:\Omega]<\infty$, we get that $\Omega$ is icc. On the other hand, $[\Omega_1,\Omega_2]$ is a finite central subgroup of $\Omega$. Thus, we must have $[\Omega_1,\Omega_2]=\{e\}$, or, in other words, $\Omega_1$ and $\Omega_2$ commute. Moreover, since $\Gamma$ is icc, it follows that both $\Omega_1$ and $\Omega_2$ are icc.

In the second part of the proof, we derive the conclusion by repeating almost verbatim part of the proof of \cite[Theorem 4.14]{CdSS15}. Nevertheless, we include details for the reader's convenience.

Since $L(\Omega_1)$ is a II$_1$ factor and $L(\Omega_1)\prec P_1$, by applying \cite[Proposition 12]{OP03}, we can find a decomposition $M=P_1^{t}\overline{\otimes}P_2^{1/t}$, for some $t>0$, and a non-zero partial isometry $v\in M$ such that $vv^*\in P_2^{1/t}, v^*v\in L(\Omega_1)'\cap M$, and 
\begin{equation}\label{omega_1}vL(\Omega_1)v^*\subset P_1^tvv^*. \end{equation}
Next, let $H_2\subset\Gamma$ be the subgroup of $g\in\Gamma$ for which $\mathcal O_{\Omega_1}(g)$ is finite. 
Then $H_2\supset\Omega_2$ and since $\Omega_1$ is icc, we get that $H_2\cap\Omega_1\Omega_2=\Omega_2$. Using that $[\Gamma:\Omega_1\Omega_2]<\infty$, we deduce that $[H_2:\Omega_2]<\infty$. Let $g_1,...,g_n\in H_2$ such that $H_2=\cup_{i=1}^n\Omega_2g_i$. Since $C_{\Omega_1}(g_i)<\Omega_1$ is a finite index subgroup, for every $i\in\{1,...,n\}$, we derive that $H_1:=C_{\Omega_1}(H_2)=\cap_{i=1}^nC_{\Omega_1}(g_i)$ is a finite index subgroup of $\Omega_1$.   Since $[\Omega_1\Omega_2:H_1\Omega_2]\leq [\Omega_1:H_1]<\infty$ and $H_1\Omega_2\subset H_1H_2$, we get that $[\Gamma:H_1H_2]<\infty$. In particular, it follows that the commuting subgroups $H_1,H_2<\Gamma$ are icc.

Since $H_1\subset\Omega_1$, by equation \ref{omega_1} we get that $vL(H_1)v^*\subset P_1^tvv^*$. Since $L(\Omega_1)'\cap M\subset L(H_2)$, we also get that $v^*v\in L(H_2)$. Note that $L(H_2)$ is a II$_1$ factor and $L(H_2)\subset L(H_1)'\cap M$. By combining these facts and proceeding as in the last paragraph of the 
 proof of \cite[Proposition 12]{OP03} (see also the proof of \cite[Theorem 4.14]{CdSS15}), we find a unitary $u\in M$ such that
  \begin{equation}\label{H_1}uL(H_1)u^*\subset P_1^t.\end{equation}
Let $\Gamma_2<\Gamma$ be the subgroup of $g\in\Gamma$ for which $\mathcal O_{H_1}(g)$ is finite. By  repeating the argument from above it follows that $\Gamma_2$ is icc,  $[\Gamma_2:H_2]<\infty$, $[H_1:C_{H_1}(\Gamma_2)]<\infty$, and $[\Gamma:C_{H_1}(\Gamma_2)\Gamma_2]<\infty$. Since $L(H_1)'\cap M\subset L(\Gamma_2)$, equation \ref{H_1} implies that 
\begin{equation}\label{inclusion}  uL(\Gamma_2)u^*\supset P_2^{1/t}.\end{equation}
 Since $L(\Gamma_2)$ is a II$_1$ factor, by using \ref{inclusion} and applying \cite[Theorem A]{Ge95}, we find a factor $A\subset P_1^t$ such that $uL(\Gamma_2)u^*=A\overline{\otimes}P_2^{1/t}$.  Since $[\Gamma_2:H_2]<\infty$ and $[H_2:\Omega_2]<\infty$, we have that $[\Gamma_2:\Omega_2]<\infty$. In particular, we conclude that $\Gamma_2$ and $\Sigma_2$ are commensurable. Using that $L(\Omega_2)\prec P_2$, we get that $L(\Gamma_2)\prec P_2$, hence $A\prec P_2$. In combination with $A\subset P_1^t$, this implies that $A$ is not diffuse.   Since $A$ is a factor, it must be finite dimensional, hence $A=\mathbb M_k(\mathbb C)$, for some $k\geq 1$. Denoting $s=t/k$, we obtain a decomposition $M =P_1^s\overline{\otimes}P_2^{1/s}$ such that 
 \begin{equation}\label{equal}uL(\Gamma_2)u^*=P_2^{1/s}.\end{equation}
 Finally, let $\Gamma_1<\Gamma$ be the subgroup of $g\in\Gamma$ for which $\mathcal O_{\Gamma_2}(g)$ is finite. Then $C_{H_1}(\Gamma_2)\subset\Gamma_1$, and since $\Gamma_2$ is icc we have that $\Gamma_1\cap C_{H_1}(\Gamma_2)\Gamma_2\subset C_{H_1}(\Gamma_2)$. Using that $[\Gamma:C_{H_1}(\Gamma_2)\Gamma_2]<\infty$, we get that $[\Gamma_1:C_{H_1}(\Gamma_2)]<\infty$. In combination with $[k\Sigma_1k^{-1}:\Omega_1]<\infty$, $[\Omega_1:H_1]<\infty$ and $[H_1:C_{H_1}(\Gamma_2)]<\infty$, this implies that $\Gamma_1$ and $k\Sigma_1 k^{-1}$ are commensurable. 
 
 Using \ref{equal}, we get that $P_1^s= u(L(\Gamma_2)'\cap M)u^*\subset uL(\Gamma_1)u^*$. By applying \cite[Theorem A]{Ge95} again, we find a von Neumann subalgebra $B\subset P_2^{1/s}$ such that $uL(\Gamma_1)u^*=P_1^{s}\overline{\otimes}B$. Since $\Gamma_2$ is icc, we get that $B=uL(\Gamma_1)u^*\cap uL(
 \Gamma_2)u^*=uL(\Gamma_1\cap\Gamma_2)u^*=\mathbb C1$.
 Therefore, we have that 
 \begin{equation}\label{gamma_1}uL(\Gamma_1)u^*=P_1^{s}.\end{equation}
 It is now clear that \ref{equal} and \ref{gamma_1} imply that $\Gamma=\Gamma_1\times\Gamma_2$, which finishes the proof.
\hfill$\blacksquare$

\section{Proofs of main results}

In this section we prove Theorems \ref{A} and \ref{C}, and Corollary \ref{D}.

\subsection{A strengthening of Theorem \ref{C}}
We establish the following strengthening of Theorem \ref{C}. This result will also be used to derive Theorem \ref{A}.

\begin{theorem}\label{general}
Let $\Gamma$ be a countable icc group and assume that $\Gamma$ is measure equivalent to a product $\Lambda=\Lambda_1\times...\times\Lambda_n$ of $n\geq 1$  groups $\Lambda_1,...,\Lambda_n$ which belong to $\mathcal C_{\text{rss}}$. Assume the notation from \ref{ME}. Suppose that $L(\Gamma)=P_1\overline{\otimes}P_2$, for some II$_1$ factors $P_1$ and $P_2$.

Then there exist a decomposition  $\Gamma=\Gamma_1\times\Gamma_2$, a partition $S_1\sqcup S_2=\{1,...,n\}$, a decomposition $L(\Gamma)=P_1^{t}\overline{\otimes}P_2^{1/t}$, for some $t>0$, and a unitary $u\in L(\Gamma)$  such that 
\begin{enumerate}
\item  $P_1^t=uL(\Gamma_1)u^*$ and $P_2^{1/t}=uL(\Gamma_2)u^*$,
\item $A\rtimes\Gamma_i\prec^s_{M}B\rtimes\Lambda_{S_i}$, $B\rtimes\Lambda_{S_i}\prec^s_{M} A\rtimes\Gamma_i$ for every $i\in\{1,2\}$, and
\item $\Gamma_i$ is measure equivalent to $\Lambda_{S_i}$, for every $i\in\{1,2\}$.

\end{enumerate}
\end{theorem}

 {\it Proof}. By applying Theorem \ref{tensor} we find subgroups $\Sigma_1,\Sigma_2<\Gamma$ and a partition $S_1\sqcup S_2=\{1,...,n\}$ such that the following conditions hold for all $i\in\{1,2\}$:
\begin{enumerate}[label=(\alph*)]
\item \text{$P_i\prec_{L(\Gamma)}^sL(\Sigma_i)$, $L(\Sigma_i)\prec_{L(\Gamma)}^sP_i$,} 
 \item \text{$A\rtimes\Sigma_i\prec^s_{M}B\rtimes\Lambda_{S_i}, B\rtimes\Lambda_{S_i}\prec^s_{M} A\rtimes\Sigma_i$, and}
 \item \text{$\Sigma_i$ is measure equivalent to $\Lambda_{S_i}$.}
 \end{enumerate}
Further, by using (a), Theorem \ref{product} provides decompositions $\Gamma=\Gamma_1\times\Gamma_2$ and $L(\Gamma)=P_1^s\overline{\otimes}P_2^{1/s}$, for some $s>0$, and a unitary $u \in L(\Gamma)$ such that $\Gamma_1$ is commensurable to $k\Sigma_1k^{-1}$, for some $k\in\Gamma$, $\Gamma_2$ is commensurable to $\Sigma_2$, and condition (1) is satisfied.
It is clear that (b) implies (2).
Finally, since commensurable groups are measure equivalent, we deduce that $\Gamma_i$ is measure equivalent to $\Sigma_i$ hence to $\Lambda_{S_i}$, for all $i\in\{1,2\}$. 
This shows that condition (3) also holds and finishes the proof. \hfill$\blacksquare$

\subsection{Proof of Theorem \ref{C}} Since non-elementary hyperbolic groups belong to $\mathcal C_{\text{rss}}$ by \cite{PV12}, Theorem \ref{C} follows from Theorem \ref{general}. \hfill$\blacksquare$

\subsection{Proof of Theorem \ref{A}} By Remark \ref{L}(1), any irreducible lattice  in a product of connected non-compact rank one simple Lie groups with finite center belongs to $\mathscr L$. Thus, it suffices to prove the second assertion of Theorem \ref{A}.

Let $\Gamma\in\mathscr L$ be an icc group and assume by contradiction that the II$_1$ factor $L(\Gamma)$ is not prime.
Then $\Gamma$ is an irreducible lattice in a product $G=G_1\times...\times G_n$ of $n\geq 1$ locally compact groups, each admitting a non-elementary hyperbolic lattice $\Lambda_j<G_j$, and not all admitting an open normal compact subgroup. Moreover, $\Gamma$ does not contain a non-trivial element which commutes with an open subgroup of $G$.
Denote $\Lambda=\Lambda_1\times...\times\Lambda_n$. Then $\Lambda<G$ is also a lattice, and hence $\Gamma$ and $\Lambda$ are measure equivalent. Since non-elementary hyperbolic groups belong to $\mathcal C_{\text{rss}}$ by \cite{PV12}, we deduce that $\Gamma$ satisfies the hypothesis of Theorem \ref{general}. 

To get a contradiction we will apply Theorem \ref{general}. 
We begin by defining a concrete stable orbit equivalence between certain actions of $\Gamma$ and $\Lambda$.
 Let $m$ be a fixed Haar measure of $G$,  consider the left-right translation action $\Gamma\times\Lambda\curvearrowright (G,m)$ given by $(g,h)\cdot x=gxh^{-1}$, and put $\mathcal R=\mathcal R(\Gamma\times\Lambda\curvearrowright G)$.

Let $X=G/\Lambda$ and $Y=\Gamma\backslash G$, endowed with left and right translation actions of $G$, and the unique $G$-invariant probability measures $m_X$ and $m_Y$. Let $p:X\rightarrow G$ and $q:Y\rightarrow G$ be Borel maps such that $p(x)\in x\Lambda$ and $q(y)\in\Gamma y$, for all $x\in X,y\in Y$. 
Let $\ell\geq 1$ such that $\ell \;m(q(Y))\geq m(p(X))$.  Let $\{X_j\}_{1\leq j\leq\ell}$ be a measurable partition of $X$ such that $m(p(X_j))\leq m(q(Y))$, for every $1\leq j\leq\ell$. Since $\mathcal R$ is ergodic, we can find $\{\theta_j\}_{1\leq j\leq\ell}\subset [\mathcal R]$ such that 
$\theta_j(p(X_j))\subset q(Y)$, for every $1\leq j\leq\ell$.  Let $\alpha_j:G\rightarrow\Gamma, \beta_j:G\rightarrow\Lambda$ be Borel maps such that $\theta_j(x)=\alpha_j(x)x\beta_j(x)$, for almost every $x\in G$.

 We define $\iota:X\rightarrow Y\times\mathbb Z/\ell\mathbb Z$ by letting $$\text{$\iota(x)=(\Gamma\theta_j(p(x)),j+\ell\mathbb Z)$, if $x\in X_j$, for some $1\leq j\leq\ell$.}$$

 We view $X$ as a subset of $Y\times\mathbb Z/\ell\mathbb Z$ by identifying it with $\iota(X)$. Fix $x_1,x_2\in X$ and let $1\leq j_1,j_2\leq\ell$ such that $x_1\in X_{j_1},x_2\in X_{j_2}$.  Then $x_1\in\Gamma x_2$ iff $p(x_1)\in\Gamma p(x_2)\Lambda$ iff  $\theta_{j_1}(p(x_1))\in\Gamma\theta_{j_2}(p(x_2))\Lambda$ iff $\Gamma\theta_{j_1}(p(x_1))\in (\Gamma\theta_{j_2}(p(x_2)))\Lambda.$  Thus, if $\mathbb Z/\ell\mathbb Z$ acts on itself by addition, then $$\mathcal R(\Gamma\curvearrowright X)=\mathcal R(\Lambda\times\mathbb Z/\ell\mathbb Z\curvearrowright Y\times\mathbb Z/\ell\mathbb Z)|_{X}.$$ 
 
Since $\Gamma$ does not contain a non-trivial element which commutes with an open subgroup of $G$, it is easy to see that the actions $\Gamma\curvearrowright (X,\mu_X)$ and $\Lambda\curvearrowright (Y,\mu_Y)$ are free.
 
 We are therefore in the situation from \ref{ME}, so we may assume the notation introduced therein: $A=L^{\infty}(X), B=L^{\infty}(Y)\otimes\mathbb M_{\ell}(\mathbb C)$, $M=L^{\infty}(Y\times\mathbb Z/\ell\mathbb Z)\rtimes(\Lambda\times\mathbb Z/\ell\mathbb Z)=B\rtimes\Lambda$. 
 We denote by $\{u_g\}_{g\in\Gamma}\subset A\rtimes\Gamma$ and $\{v_{h}\}_{h\in\Lambda}\subset M$ the canonical unitaries.
 Additionally, we let $\Lambda_{S}=\times_{i\in S}\Lambda_i$, $G_S=\times_{i\in S}G_i$, and $\pi_S:G\rightarrow G_S$ denote the canonical projection, for every subset $S\subset\{1,...,n\}$.

 Since $L(\Gamma)$ is  not prime, Theorem \ref{general}  implies that we can find a  decomposition $\Gamma=\Gamma_1\times\Gamma_2$, with $\Gamma_1$ and $\Gamma_2$ icc, and a partition $S_1\sqcup S_2=\{1,...,n\}$ such that $A\rtimes\Gamma_i\prec^s_{M}B\rtimes\Lambda_{S_i}$, for all $i\in\{1,2\}$. The rest of the proof relies on the following:
 
 {\bf Claim}. The subgroups $\overline{\pi_{S_1}(\Gamma_2)}\subset G_{S_1}$ and $\overline{\pi_{S_2}(\Gamma_1)}\subset G_{S_2}$ are compact.

{\it Proof of the claim.} By symmetry, it suffices to prove the first assertion. Assume by contradiction that $\overline{\pi_{S_1}(\Gamma_2)}$ is not compact. Then we can find a sequence $g_n\in \Gamma_2$ such that $\pi_{S_1}(g_n)\rightarrow\infty$, as $n\rightarrow\infty$, in $G_{S_1}$. We claim that 
\begin{equation}\label{noncpt}\text{$\|E_{B\rtimes\Lambda_{S_2}}(u_{g_n}v_k^*)\|_2\rightarrow 0$, for every $k\in\Lambda_{S_1}$.}\end{equation}

Since $E_{B\rtimes\Lambda_{S_2}}$ is $B\rtimes\Lambda_{S_2}$-bimodular and $M$ is generated by $B\rtimes\Lambda_{S_2}$ together with the unitaries $\{v_k \;|\; k \in \Lambda_{S_1}\}$ that normalize it, claim \ref{noncpt} readily implies that $\|E_{B\rtimes\Lambda_{S_2}}(au_{g_n}b)\|_2\rightarrow 0$, for every $a,b\in M$, which contradicts that $A\rtimes\Gamma_2\prec_{M}B\rtimes\Lambda_{S_2}$. 

For $1\leq j\leq\ell$, let $e_j\in L^{\infty}(X)$ denote the characteristic function of $X_j$. Since $\sum_{1\leq j\leq\ell}e_{j}=1_{X}$, claim \ref{noncpt} reduces to proving \begin{equation}\label{noncpt2}\text{$\|E_{B\rtimes\Lambda_{S_2}}(e_{j_1}u_{g_n}e_{j_2}v_k^*)\|_2\rightarrow 0$, for every $k\in\Lambda_{S_1}$ and $1\leq j_1,j_2\leq\ell$.}\end{equation}

To prove \ref{noncpt2}, fix $k\in\Lambda_{S_1}$ and $1\leq j_1,j_2\leq\ell$. For $g\in\Gamma$, the Fourier expansion of $e_{j_1}u_ge_{j_2}$ in $M=B\rtimes\Lambda$ is given by $e_{j_1}u_ge_{j_2}=\sum_{h\in\Lambda\times\mathbb Z/\ell\mathbb Z}{1}_{\{x\in X_{j_1}\cap gX_{j_2}|g^{-1}x=h^{-1}x\}}v_h$.
If $x\in X_{j_1}\cap gX_{j_2}$, then $\iota(x)=(\Gamma g^{-1}p(x)\beta_{j_1}(p(x)),j_1+\ell\mathbb Z)$ and $\iota(g^{-1}x)=(\Gamma p(g^{-1}x)\beta_{j_2}(p(g^{-1}x)),j_2+\ell\mathbb Z)$.
Thus, denoting $$w(x)=\beta_{j_1}(p(x))^{-1}p(x)^{-1}gp(g^{-1}x)\beta_{j_2}(p(g^{-1}x))\in\Lambda,$$ and recalling that the action $\Lambda\curvearrowright (Y,\mu_Y)$ is free, we get that
\begin{equation}\label{noncpt3} e_{j_1}u_ge_{j_2}=\sum_{h\in\Lambda}1_{\{x\in X_{j_1}\cap gX_{j_2}|w(x)=h\}}\;v_{(h,j_1-j_2+\ell\mathbb Z)}.
\end{equation}

From this it follows that \begin{equation}\label{noncpt4}\text{$\|E_{B\rtimes\Lambda_{S_2}}(e_{j_1}u_ge_{j_2}v_k^*)\|_2^2\leq m_X(\{x\in X|w(x)\in\Lambda_{S_2}k\})$, for every $g\in\Gamma$}.\end{equation}

Now, let $\varepsilon>0$. Then we can find a compact set $C\subset G_{S_1}$ such that we have
\begin{itemize}
\item $\mu_X(\{x\in X|\pi_{S_1}(p(x))\notin C\}\leq\frac{\varepsilon}{4}$, and
\item $\mu_X(\{x\in X|\pi_{S_1}(\beta_j(p(x))\notin C\}\leq\frac{\varepsilon}{4}$, for $j\in\{j_1,j_2\}$.
\end{itemize}

If $x\in X$ satisfies $w(x)\in\Lambda_{S_2}k$, then $\pi_{S_1}(w(x))=k$. By using the definition of $w(x)$, the fact that  the action of $\Gamma$ on $X$ is measure preserving, and the last two inequalities one obtains that 
\begin{equation}\label{noncpt5}\mu_{X}(\{x\in X|w(x)\in\Lambda_{S_2}k\})\leq \varepsilon+\mu_{X}(\{x\in X|k\in (C^{-1})^2\pi_{S_1}(g)C^2\}).\end{equation} 

By combining \ref{noncpt4} and \ref{noncpt5} we derive  that 

$$\text{$\|E_{B\rtimes\Lambda_{S_2}}(e_{j_1}u_ge_{j_2}v_k^*)\|_2^2\leq \varepsilon+\mu_{X}(\{x\in X|k\in (C^{-1})^2\pi_{S_1}(g)C^2\})$, for every $g\in\Gamma$.}$$

Since $\pi_{S_1}(g_n)\rightarrow\infty$, we have that $k\notin(C^{-1})^2\pi_{S_1}(g_n)C^2$, for large enough $n$. Therefore, the last inequality implies that $\limsup_{n\rightarrow\infty}\|E_{B\rtimes\Lambda_{S_2}}(e_{j_1}u_{g_n}e_{j_2}v_k^*)\|_2^2\leq\varepsilon$. Since $\varepsilon>0$ is arbitrary, this proves \ref{noncpt2} and thus the claim.
\hfill$\square$

To finish the proof, let $i\in\{1,2\}$. Since $\Gamma_i$ is infinite and $A\rtimes\Gamma_i\prec^s_{M}B\rtimes\Lambda_{S_i}$, we get that $\Lambda_{S_i}$ is infinite, hence $S_i$ is nonempty. Thus, $S_i$ is a proper subset of $\{1,...,n\}$. Therefore, since $\Gamma$ is an irreducible lattice in $G$, we derive that $\pi_{S_i}(\Gamma)<G_{S_i}$ is dense. In combination with the claim, this implies that $K_1=\overline{\pi_{S_1}(\Gamma_2)}$ and $K_2=\overline{\pi_{S_2}(\Gamma_1)}$ are normal compact subgroups of $G_{S_1}$ and $G_{S_2}$, respectively.
 Thus, $K=K_1\times K_2$ is a normal compact subgroup of $G=G_{S_1}\times G_{S_2}$. 
 
 Let $\rho_i:G_{S_i}\rightarrow G_{S_i}/K_i$, for $i\in\{1,2\}$, and $\rho=(\rho_1,\rho_2):G\rightarrow G/K$ be the canonical projections.  If  $g_1\in\Gamma_1$ and $g_2\in\Gamma_2$, then $\rho_1(\pi_{S_1}(g_2))=\text{id}$ and $\rho_2(\pi_{S_2}(g_1))=\text{id}$. Thus, we derive that $\rho(g_1g_2)=(\rho_1(\pi_{S_1}(g_1g_2)),\rho_2(\pi_{S_2}(g_1g_2)))=(\rho_1(\pi_{S_1}(g_1)),\rho_2(\pi_{S_2}(g_2)))$, which implies that \begin{equation}\label{noncpt6} \rho(\Gamma)=\rho_1(\pi_{S_1}(\Gamma_1))\times\rho_2(\pi_{S_2}(\Gamma_2)).\end{equation}
 
 If $i\in\{1,2\}$, then $\pi_{S_i}(\Gamma)<G_{S_i}$ is dense, hence $\rho_i(\pi_{S_i}(\Gamma))=\rho_i(\pi_{S_i}(\Gamma_i))$ is dense in $G_{S_i}/K_i$. 
 In combination with \ref{noncpt6}, we conclude that $\rho(\Gamma)<G/K$ is dense.
 On the other hand, since $\Gamma<G$ is discrete and $K<G$ is compact, we get that $\rho(\Gamma)<G/K$ is discrete hence closed. Altogether, we deduce that $\rho(\Gamma)=G/K$ and thus $K<G$ is an open normal compact subgroup. 
 This implies that $\pi_{\{j\}}(K)<G_j$ is an open normal compact subgroup, for every $1\leq j\leq n$, a contradiction.
\hfill$\blacksquare$

\subsection{Proof of Corollary \ref{D}}
Let $k\geq 1$ be the largest integer for which there are a decomposition $\Gamma=\Gamma_1\times...\times\Gamma_k$ and a partition $T_1\sqcup...\sqcup T_k=\{1,...,n\}$ such that $T_i$ is non-empty and $\Gamma_i$ is measure equivalent to $\underset{j\in T_i}{\times}\Lambda_j$, for all $1\leq i\leq k$.
Theorem \ref{C} implies that $L(\Gamma_i)$ is a prime II$_1$ factor, for all $1\leq i\leq k$. This proves the existence of a decomposition with the desired property.

In order to prove the uniqueness of the decomposition, we  establish the following fact:
if $\Gamma=\Sigma_1\times\Sigma_2$, then there is a partition $I_1\sqcup I_2=\{1,...,k\}$ such that $\Sigma_1=\underset{i\in I_1}{\times}\Gamma_i$ and $\Sigma_2=\underset{i\in I_2}{\times}\Gamma_i$.  
 To see this, for $1\leq i\leq k$, let $\pi_i:\Gamma\rightarrow\Gamma_i$ be the canonical projection. Then $\Gamma_i$ is generated by the commuting subgroups $\pi_i(\Sigma_1)$ and $\pi_i(\Sigma_2)$. Since $\Gamma_i$ has trivial center, we have that $\pi_i(\Sigma_1)\cap\pi_i(\Sigma_2)=\{e\}$, which implies that $\Gamma_i=\pi_i(\Sigma_1)\times\pi_i(\Sigma_2)$. Since $L(\Gamma_i)$ is prime, we deduce that either $\pi_i(\Sigma_1)=\{e\}$ or $\pi_i(\Sigma_2)=\{e\}$. Since this holds for every $1\leq i\leq k$, the fact follows.

 Now, if $\Gamma=\Sigma_1\times...\times\Sigma_l$ is another decomposition such that $L(\Sigma_j)$ is a prime II$_1$ factor, for every $1\leq j\leq l$, then the fact implies that $l=k$ and that, after a permutation of indices, we have $\Sigma_i=\Gamma_i$, for every $1\leq i\leq k$.

(1) Assume that $M=P_1\overline{\otimes}P_2$, for some II$_1$ factors $P_1$ and $P_2$. By applying Theorem \ref{C} we find a decomposition $\Gamma=\Sigma_1\times\Sigma_2$, a decomposition $M=P_1^t\overline{\otimes}P_2^{1/t}$, for some $t>0$, and a unitary $u\in M$, such that $P_1^t=uL(\Sigma_1)u^*$ and $P_2^{1/t}=uL(\Sigma_2)u^*$. The above fact now clearly implies the conclusion.

(2) $\&$ (3) Assume that $M=P_1\overline{\otimes}\dots\overline{\otimes}P_m$, where $P_1,\dots,P_m$ are II$_1$ factors.
Then by induction, part (1) implies that $m\leq k$ and there are a partition $I_1\sqcup\dots\sqcup I_m=\{1,\dots,k\}$, a decomposition $M=P_1^{t_1}\overline{\otimes}\dots\overline{\otimes}P_m^{t_m}$, for some $t_1,\dots,t_m>0$ with $t_1\dots t_m=1$, and a unitary $u\in M$ such that $P_j^{t_j}=u(\overline{\otimes}_{i\in I_j}L(\Gamma_i))u^*$, for every $1\leq j\leq m$. 

If $m\geq k$, then we get that $m=k$. Since $I_j$ is nonempty, it follows that $I_j$ consists of one element, for every $1\leq j\leq m$. This implies part (2).
If $P_j$ is prime, for every $1\leq j\leq m$, then again it follows that $I_j$ consists of one element, for every $1\leq j\leq m$. This implies part (3).
\hfill$\blacksquare$

\end{document}